\documentclass[10pt]{article} 
\usepackage[preprint]{tmlr}

\usepackage{hyperref}       
\usepackage{url}            
\usepackage{booktabs}       
\usepackage{amsfonts}       
\usepackage{nicefrac}       
\usepackage{microtype}      
\usepackage{bbold}
\usepackage{capt-of}

\usepackage{amsmath, amsthm, amssymb, graphicx, mathabx, epsfig,epstopdf,color,soul}

\usepackage{algorithm,tabularx}

\let\classAND\AND
\let\AND\relax
\usepackage{algorithmic}

\let\AND\classAND
\AtBeginEnvironment{algorithmic}{\let\AND\algoAND}

\allowdisplaybreaks

\usepackage{multicol}
\usepackage{empheq} 
\usepackage{caption} 
\usepackage{float}

\newcommand{\Acal}{{\cal A}}

\newcommand{\Ecal}{{\cal E}}
\newcommand{\Fcal}{{\cal F}}
\newcommand{\Gcal}{{\cal G}}

\newcommand{\Ncal}{{\cal N}}
\newcommand{\Ocal}{{\cal O}}
\newcommand{\Pcal}{{\cal P}}
\newcommand{\Qcal}{{\cal Q}}

\newcommand{\Scal}{{\cal S}}

\newcommand{\Vcal}{{\cal V}}

\newcommand{\Xcal}{{\cal X}}

\newcommand{\1}{{\mathbf{1}}}

\newcommand{\argmin}{\mathop{\rm argmin}}
\newcommand{\argmax}{\mathop{\rm argmax}}

\newtheorem{prop}{Proposition}
\newtheorem{lem}{Lemma}
\newtheorem{thm}{Theorem}
\newtheorem{cor}{Corollary}
\newtheorem{definition}{Definition}

\newtheorem{assump}{Assumption}
\newtheorem{remark}{Remark}

\newcommand{\vct}[1]{\boldsymbol{#1}}

\newcommand{\vpi}{\vct{\pi}}

\newcommand{\diag}{\operatorname{diag}}

\title{Natural Policy Gradient and Actor Critic Methods for Constrained Multi-Task Reinforcement Learning}

\author{\name Sihan Zeng \email sihan.zeng@jpmchase.com \\
      \addr JPMorgan AI Research
      \AND
      \name Thinh T. Doan \email thinhdoan@vt.edu \\
      \addr Department of Electrical and Computer Engineering\\
      Virginia Tech
      \AND
      \name Justin Romberg \email jrom@ece.gatech.edu\\
      \addr Department of Electrical and Computer Engineering\\
      Georgia Tech}

\def\month{MM}  
\def\year{YYYY} 
\def\openreview{\url{https://openreview.net/forum?id=XXXX}} 

\begin{document}

\maketitle
\thispagestyle{empty}

\begin{abstract}
Multi-task reinforcement learning (RL) aims to find a single policy that effectively solves multiple tasks at the same time. This paper presents a constrained formulation for multi-task RL where the goal is to maximize the average performance of the policy across tasks subject to bounds on the performance in each task. We consider solving this problem both in the centralized setting, where information for all tasks is accessible to a single server, and in the decentralized setting, where a network of agents, each given one task and observing local information, cooperate to find the solution of the globally constrained objective using local communication.

We first propose a primal-dual algorithm that provably converges to the globally optimal solution of this constrained formulation under exact gradient evaluations. When the gradient is unknown, we further develop a sampled-based actor-critic algorithm that finds the optimal policy using online samples of state, action, and reward. Finally, we study the extension of the algorithm to the linear function approximation setting.
\end{abstract}

\section{Introduction}

Multi-task reinforcement learning (RL) aims to find a common policy that effectively solves a range of tasks simultaneously, where each task is the policy optimization problem defined over a Markov decision process (MDP). 
The MDPs can have different state spaces, reward functions, and transition kernels in general, but may be implicitly or explicitly correlated.

The most common mathematical formulation for multi-task RL is to maximize the average cumulative rewards collected by a single policy across all MDPs \cite{zeng2021decentralized,jiang2022mdpgt,junru2022decentralized}. In this paper, we study a generalized formulation in which we maximize the average cumulative rewards subject to constraints on the performance of the policy for each task. This formulation is a special case of the policy optimization problem for a constrained Markov decision process (CMDP) \cite{altman1999constrained} and is a flexible framework that allows more fine-grained specification of the performance of the optimal policy in each task.
In applications where the tasks exhibit major conflicts of interest and/or the magnitude of the rewards varies significantly across tasks \cite{kalashnikov2021mt,guo2022learning}, the optimal policy under the average-cumulative-reward formulation may excel in some tasks at the cost of compromised performance in others \cite{hayes2022practical}. The constrained formulation provides a way to mitigate this task imbalance.  Illustrative numerical simulations are given in Section~\ref{sec:simulations}.

Under the constrained multi-task formulation, we consider centralized and decentralized learning paradigms. ``Centralized'' in this context means that information of all tasks is available at a single server, while ``decentralized'' is a scenario where a group of agents, each deployed to one local environment/task, work together to solve the global constrained optimization program without any central coordination.
The centralized setting can be regarded as an easier special case of the decentralized problem with a fully connected communication graph. 
We will initiate our algorithmic development in the centralized setting and extend them to the decentralized scenario.


To solve the constrained multi-task policy optimization problem, we propose a multi-task primal-dual natural policy gradient algorithm (MT-PDNPG) that allows the policy at each local agent to achieve global optimality. 
%
%
The updates in MT-PDNPG require the computation of the exact gradients of the value function, which is impractical in environments with large state spaces and/or unknown transition probability kernel.
To extend our algorithm to the case where we do not have perfect knowledge of the environments and can only obtain samples of the state transitions, we propose a sampled-based multi-task primal-dual natural actor-critic algorithm (MT-PDNAC) and study its finite-sample performance. 
%
Finally, to tackle problems where the state space is enormous or even infinitely large, we extend MT-PDNAC to the case where the policy and value functions are linearly approximated with pre-determined lower-dimensional feature vectors. Despite the complications introduced to the algorithm and analysis by the linear function approximation, we show that a modified version of the MT-PDNAC algorithm in this setting achieves the same convergence rate as in the tabular case. We note that our proposed sampled-based algorithms are completely online in the sense that they use a single trajectory of continuously generated samples, which makes them convenient to implement in practice.

\subsection{Main Contribution}



The first contribution of our work is to study the constrained multi-task RL formulation and to propose and analyze a primal-dual natural actor-critic algorithm provably converges to the globally optimal policy in expectation. The algorithm is completely data-driven, single-loop, and relies on a single trajectory of samples. After $K$ iterations, the policy parameter converges in both objective function and constraint violation up to the precision $\Ocal(1/K^{1/6})$. This matches the best-known time and sample complexity of the natural actor-critic algorithm for single-agent (non-constrained) MDPs under comparable assumptions.

Our second contribution is to extend the primal-dual algorithms to the decentralized learning paradigm. The extended algorithm makes each agent compute and update in the direction of a locally observable component of the policy gradient, followed by a parameter averaging step. We show that the decentralized algorithms enjoy finite time and sample complexity matching their centralized counterparts, differing only by a factor that scales inversely with the connectivity of the communication graph.

Finally, we study the setting where both the policy and the critic variables are approximated using linear features. This on-policy linear function approximation setting presents peculiar technical challenges. Specifically, the TD learning target under linear function approximation becomes ill-defined when the policy to evaluate is not completely mixed (i.e. does not have uniformly lower bounded entries). While a well-defined TD learning target can be ensured through careful control of the policy iterates, it is not sufficiently Lipschitz continuous for the direct extension of our algorithm for the tabular setting to converge. 
We overcome the challenge by dynamically adjusting the number of TD learning iterations devoted to chasing a fixed TD target, which results in a modified algorithm with a sample complexity of $\Ocal(\epsilon^{-6})$, matching that in the tabular case.
To our best knowledge, a finite-time and finite-sample analysis for the on-policy natural actor-critic algorithm (even for standard non-constrained MDPs) has been missing from the existing literature, and our work fills in this gap.

\subsection{Related Work}

This paper presents reliably convergent decentralized algorithm for finding the optimal solution of the new constrained multi-task RL objective. It closely relates to the literature on multi-task RL, decentralized optimization, CMDP, and actor-critic algorithms in RL, which we discuss in this section to give context to our novelty.

\noindent\textbf{Multi-Task Reinforcement Learning.} Multi-task RL in general studies efficiently solving the policy optimization tasks for multiple RL environments at the same time by leveraging connections between the tasks. Its most common mathematical formulation is to find a single policy that maximizes the (weighted) average of the cumulative returns collected across all environments, and \cite{zeng2021decentralized,jiang2022mdpgt,junru2022decentralized,chen2022decentralized} study various gradient-based algorithms that provably converge to global or local solutions of this objective. 
However, as pointed out in \cite{hessel2019multi}, this average return formulation can be inadequate when modelling practical problems where the tasks have strong conflicting or imbalanced interests. While the authors in \cite{hessel2019multi} address this issue by dynamically addressing the weight of each task in the policy updates, we are motivated to propose the constrained multi-task formulation that allows fine-grained control of the performance of the policy in each environment.

It is worth pointing out the large volume of literature that approaches multi-task RL from a more empirical perspective. Some important lines of work include (but are certainly not limited to) policy distillation \cite{rusu2015policy,traore2019discorl,wadhwania2019policy}, transfer learning \cite{gupta2017learning,dsharing}, and innovated design of the policy representation \cite{yang2020multi,hong2021structure}. We also note that there exist multi-task RL formulations where rather than learning a single policy, task-specific adaptation is allowed \cite{finn2017model,raghurapid}.

\noindent\textbf{Decentralized Optimization.} Closely connected to distributed and decentralized optimization, our problem formulation considers maximizing a global objective function composed of local rewards under local constraints in the case where each agent only has access to its local information. In the unconstrained setting, it has been well-known that decentralized gradient, sub-gradient, and Newton's methods converge as fast as their centralized counterparts up to a constant that describes the connectivity of the communication graph \cite{nedic2009distributed,yuan2016convergence,nedic2020distributed,bullins2021stochastic,islamov2021distributed}.
For constrained convex optimization programs, \cite{chang2014distributed,lei2016primal} show that primal-dual algorithms provably converge to the globally optimal solution. Our work is inspired by these existing results but focuses on a non-convex optimization program where we establish global convergence using the specific problem structure.

\noindent\textbf{Constrained Markov Decision Process.}
Our constrained multi-task formulation can be regarded as a special case of policy optimization under a CMDP (if the information for all tasks is centrally available). 
A common approach to finding the optimal policy for a CMDP is to search for a saddle point of the Lagrangian using primal-dual gradient descent ascent. Variants of this approach are considered in \cite{prashanth2016variance,chow2017risk,tessler2018reward}. However, they only establish the asymptotic convergence to a stationary point or locally optimal solution; global optimality is not achieved and the exact finite-time complexity is unknown due to the non-convexity of the underlying optimzation program.

For a convex constrained optimization problem under Slater's condition, it is well known that strong duality holds and the primal-dual gradient descent ascent algorithm efficiently converges to the globally optimal solution.
The policy optimization problem for a CMDP is non-convex, but by leveraging the structure of the CMDP, \cite{altman1999constrained, paternain2019constrained} show that the strong duality holds despite the non-convexity.
A number of works \cite{ding2020natural,liu2021policy,ding2022convergence,bai2022achieving} take advantage of this property to design primal-dual natural policy gradient algorithms and establish their finite-time convergence to the globally optimal policy in the tabular case. 
Some other works focus on deriving regret bounds (rather than finite-time convergence) \cite{zheng2020constrained,ding2021provably,agarwal2022regret}.
A limit of these studies is that the policy and dual variable updates rely on an oracle that always returns the exact gradient (or its highly accurate unbiased estimate). In real-life problems where the transition probability kernel is not fully known and/or the state space is large, computing such gradients can become computationally prohibitive.
We improve these prior works by presenting a completely sample-based algorithm that solves the CMDP optimization using a single continuously-generated trajectory, both in the tabular setting and under linear function approximation. 

There exist non-primal-dual algorithms for finding the optimal policy for a CMDP. For example, \cite{yu2019convergent} constructs a succession of surrogate convex relaxation to the non-convex CMDP optimization problem and shows that the solutions to these surrogate programs converge to a stationary point of the CMDP optimization problem. The work \cite{chow2018lyapunov} extends value-based methods including value iteration, policy iteration, and Q learning to the context of CMDP. The paper \cite{hasanzadezonuzy2021learning} builds an empirical estimate of the probability transition kernel, on which planning is carried out.
The authors of \cite{liu2021learning} propose an algorithm driven by the principle of optimistic pessimism that improves the state-of-the-art analysis on the constraint violation. In \cite{ying2022dual}, a dual-only approach is studied which solves a regularized version of the CMDP.

\noindent\textbf{Actor-Critic Algorithms.} 
The sample-based algorithms presented in our paper fall under the category of actor-critic algorithms, which can be considered as a variant of policy gradient methods where the unknown value function is estimated by an auxiliary variable updated with stochastic approximation. This class of algorithms has been analyzed in various settings of dynamical systems (such as standard MDP, entropy regularized MDP, and linear-quadratic regulator), on/off-policy sample collection, natural/standard policy gradient, and function approximation \cite{yang2019provably,wu2020finite,xu2020non,zeng2021two,ju2022model,khodadadian2022finite,chen2022finite,barakat2022analysis}. 
However, as discussed in the previous subsection, analyzing the natural actor-critic algorithm incurs unique technicality that has not been treated in the previous literature under the combined effect of on-policy samples and linear function approximation. Our work thoroughly describes the technical challenges and presents our solution.

Finally, we note that a preliminary version of the work has been presented in \cite{zeng2022finite}, which analyzes an online actor-critic algorithm for a single-agent CMDP. This current paper significantly generalizes \cite{zeng2022finite} along several axes by introducing the multi-task formulation, decentralized learning paradigm, and linear function approximation.

\section{Constrained Multi-Task Formulation}
Consider a collection of $N$ infinite-horizon discounted-reward MDPs characterized by $\{(\Scal_i,\Acal,\Pcal_i,\gamma, r_i)\}_{i=1}^{N}$, where $\Scal_i$ is the finite state space, $\Acal$ is the finite and common action space, $P_i:\Scal_i\times\Acal\rightarrow\Delta_{\Scal_i}$ is the transition probability kernel, $\gamma\in(0,1)$ is the discount factor, and $r_i:\Scal_i\times\Acal\rightarrow[0,1]$ is the reward function.

Given a policy $\pi\in\Delta_{\Acal}^{\Scal_i}$, we define the local value functions for each task $i=1,\cdots,N$
\begin{gather*}
    V_{i}^{\pi}(s)\triangleq\mathbb{E}_{\pi}\Big[\sum_{k=0}^{\infty} \gamma^k r_i\left(s_k, a_k\right) \mid s_{0}=s\Big],\quad
    Q_{i}^{\pi}(s,a)\triangleq\mathbb{E}_{\pi}\Big[\sum_{t=0}^{\infty} \gamma^k r_i\left(s_k, a_k\right) \mid s_{0}=s, a_{0}=a\Big],\notag\\
    A_{i}^{\pi}(s,a)\triangleq Q_{i}^{\pi}(s,a)-V_{i}^{\pi}(s),\quad\forall s\in\Scal,a\in\Acal.
\end{gather*}

With some abuse of notation, we use $V_{i}^{\pi}(\rho)$ to denote the expected cumulative reward collected by policy $\pi$ under the initial distribution $\rho$
\begin{align}
    V_{i}^{\pi}(\rho)\triangleq\mathbb{E}_{s_{0} \sim \rho}\left[V_{i}^{\pi}\left(s_{0}\right)\right].
\end{align}

For the simplicity of notation, we define $V_0^{\pi}(\rho)$ to be the averaged cumulative reward over tasks
\[V_0^{\pi}(\rho)\triangleq\frac{1}{N}\sum_{i=1}^{N}V_i^{\pi}(\rho).\]

The common multi-task formulation considered in the literature is to maximize this average value function
\begin{align}
\max_{\pi\in\Delta_{\Acal}^{\Scal}} V_{0}^{\pi}(\rho).\label{eq:unconstrained_MTRL}
\end{align}

Given local performance upper and lower bounds $\{\ell_i\in\mathbb{R},u_i\in\mathbb{R}\}_{i=1}^{N}$, our multi-task policy optimization objective is to solve the following constrained program
\begin{align}
    \pi^{\star}\in\argmax_{\pi\in\Delta_{\Acal}^{\Scal}} \quad &V_{0}^{\pi}(\rho)\notag\\
    \text{subject to} \quad &\ell_i\leq V_{i}^{\pi}(\rho) \leq u_i\quad\forall i = 1,2,...,N.
    \label{eq:objective}
\end{align}

We note that Eq.~\eqref{eq:objective} obviously subsumes the non-constraint multi-task formulation in \cite{zeng2021decentralized,jiang2022mdpgt,junru2022decentralized} by properly choosing $\{\ell_i,u_i\}$. It can be shown that the multi-task policy optimization problem (even without constraints) does not observe the gradient domination condition in general, which makes it difficult for any gradient-based algorithm to find the globally optimal policy. In this paper, we make the simplifying assumption that the local MDPs have the same state space and transition probability kernel and can only be different in the reward functions (i.e. $\Scal=\Scal_1=\Scal_2=\cdots$ and $\Pcal=\Pcal_1=\Pcal_2=\cdots$), under which the gradient domination condition is recovered \cite{zeng2021decentralized}.

We consider the following assumption on Eq.~\eqref{eq:objective}, which essentially states that the constraint set must have a non-empty interior. This is a standard assumption in the study of constrained MDPs \cite{paternain2019constrained,ding2020natural,zeng2022finite2} and ensures that strong duality holds for the constrained program despite the lack of convexity.

\begin{assump}[Slater's Condition]\label{assump:Slater}
There exists a constant $0<\xi\leq 1$ and a policy $\pi$ such that $\ell_i+\xi\leq V_i^{\pi}(\rho)\leq u_i-\xi$ for all $i=1,\cdots,N$.
\end{assump}

\section{Preliminary -- Centralized Computation Setting}\label{sec:MT-PDNPG}

We start by discussing how Eq.~\eqref{eq:objective} can be solved in the ``centralized'' setting where information of all tasks are available at a single server. 
For now, we work with deterministic gradients where we assume having perfect information of the environments to compute the exact value functions. Sampled-based algorithms will be built upon these results in Section~\ref{sec:AC}.
Under this simplification, the problem formulation and setting become a special case of those in \cite{ding2020natural} where the reward to be optimized is $\frac{1}{N}\sum_{i=1}^{N}r_i$ and the reward for the definition of the $i_{\text{th}}$ constraint is $r_i$. 

Inspired by \cite{ding2020natural}, we follow a primal-dual approach to solve the constrained program in Eq.~\eqref{eq:objective}. The first step is to form the Lagrangian

\begin{align}
V_{L}^{\pi,\lambda,\nu}(\rho)\triangleq V_{0}^{\pi}(\rho)+\sum_{i=1}^{N}\left(\lambda_i\left(V_{i}^{\pi}(\rho)-\ell_i\right)-\nu_i\left(V_{i}^{\pi}(\rho)-u_i\right)\right),
\label{eq:Lagrangian}
\end{align}
where $\lambda=[\lambda_1,\dots,\lambda_N]\in\mathbb{R}_{+}^N$ and $\nu=[\nu_1,\dots,\nu_N]\in\mathbb{R}_{+}^N$ are the dual variables associated with the lower and upper bounds. The dual function $V_{D}^{\lambda,\nu}$ is defined as
\begin{align}
    V_D^{\lambda,\nu}(\rho)=\max_{\pi\in\Delta_{\Acal}^{\Scal}} V_{L}^{\pi,\lambda,\nu}(\rho),
    \label{eq:def_dualfunc}
\end{align}
and the dual problem is to solve
\begin{align}
    \lambda^{\star},\nu^{\star}\in\argmin_{\lambda,\nu\,\in\,\mathbb{R}_{+}^{N}} V_D^{\lambda,\nu}(\rho).
    \label{eq:def_dualprob}
\end{align}

Although the program in Eq.~\eqref{eq:objective} is non-convex, it is known that strong duality holds under Slater's condition in Assumption~\ref{assump:Slater} (see \cite{altman1999constrained}[Theorem 3.6])
\begin{align}
    V_{D}^{\lambda^{\star},\nu^{\star}}(\rho)=V_{0}^{\pi^{\star}}(\rho),\label{lem:strong_duality}
\end{align}
where $\pi^{\star}$, $\lambda^{\star}$, and $\nu^{\star}$ are the (not necessarily unique) optimal solutions to Eq.~\eqref{eq:objective} and Eq.~\eqref{eq:def_dualprob}, respectively. 
The optimal dual variables are known to be bounded also as a consequence of Assumption~\ref{assump:Slater}, which we present in the following lemma and is a simple extension of \cite{ding2020natural}[Lemma 1] to the case of multiple constraints. 

\begin{lem}\label{lem:bounded_lambdastar}
Under Assumption \ref{assump:Slater}, we have
\begin{align*}
    \|\lambda^\star\|_{\infty}\leq \frac{B_{\lambda}}{2} \quad\text{and}\quad \|\nu^\star\|_{\infty}\leq \frac{B_{\lambda}}{2},
\end{align*}
where $B_{\lambda}=\frac{1}{\xi(1-\gamma)}$
\end{lem}

Motivated by the existence of strong duality, we solve Eq.~\eqref{eq:objective} by finding the minimax saddle point of the Lagrangian. We represent the policy through the softmax parameterization, as that introduces more favorable structure to the optimization landscape \cite{agarwal2020optimality}.
Specifically, using a policy parameter $\theta\in \mathbb{R}^{|\mathcal{S}||\mathcal{A}|}$ that encodes the policy $\pi_{\theta}$ through the softmax function as follows
\begin{align*}
    \pi_\theta(a \mid s)=\frac{\exp \left(\theta(s,a)\right)}{\sum_{a^{\prime} \in \mathcal{A}} \exp \left(\theta(s,a')\right)}, \quad \text { for all } \, \theta \in \mathbb{R}^{|\mathcal{S}||\mathcal{A}|},
\end{align*}
we take a gradient-descent-ascent approach to find the saddle point $(\theta^{\star},\lambda^{\star},\nu^{\star})$ such that
\begin{align}
    \theta^{\star},(\lambda^{\star},\nu^{\star})=\argmax_{\theta\in\mathbb{R}^{|\Scal||\Acal|}}\argmin_{\lambda,\nu\,\in\,\mathbb{R}_{+}^{N}}V_{L}^{\pi_{\theta},\lambda,\nu}.\label{eq:maximin}
\end{align}

Our approach requires computing the gradients of the Lagrangian with respect to both the policy parameter and dual variables, which we now derive.



\noindent\textbf{Gradient of the dual variable.}
The Lagrangian in Eq.~\eqref{eq:Lagrangian} is obviously a linear function of $\lambda$ and $\nu$, and the gradients have simple closed-form expressions.
\begin{align*}
    \nabla_{\lambda_i}V_L^{\pi,\lambda,\nu}(\rho)&=V_i^{\pi}(\rho)-\ell_i=\sum_{s:\rho(s)>0,a}\rho(s) \pi(a\mid s) Q_i^{\pi}(s,a)-\ell_i,\notag\\
    \nabla_{\nu_i}V_L^{\pi,\lambda,\nu}(\rho)&=-V_i^{\pi}(\rho)+u_i=-\sum_{s:\rho(s)>0,a}\rho(s) \pi(a\mid s) Q_i^{\pi}(s,a)+u_i.
\end{align*}
This naturally leads to the update in Eq.~\eqref{Alg:MT-PDNPG:dual_update}, in which the operator $\Pi_{[0,B_{\lambda}]}:\mathbb{R}^{N}\rightarrow\mathbb{R}^{N}$ denotes the element-wise projection of a vector to the interval $[0,B_{\lambda}]$. The projection guarantees the stability of the dual variables. We note that the optimal dual variables are in the range $[0,B_{\lambda}]$ according to Lemma~\ref{lem:bounded_lambdastar}.

\noindent\textbf{Gradient of the primal variable.}
It is known that the natural gradient of the value function under reward $r_i$ with respect to $\theta$, denoted by $\widetilde{\nabla}_{\theta}V_i^{\pi_{\theta}}(\rho)$, is the advantage function $A_i^{\pi_{\theta}}$ weighted by $1/(1-\gamma)$ \cite{agarwal2020optimality}. This means that $\widetilde{\nabla}_{\theta}V_L^{\pi_{\theta},\lambda,\nu}(\rho)$, the natural policy gradient of the Lagrangian can be expressed as\looseness=-1
\begin{align}
&\widetilde{\nabla}_{\theta(s,a)} V_L^{\pi_\theta, \lambda,\nu}(\rho)\notag\\
&\hspace{40pt}=\frac{1}{1-\gamma}\big(A_0^{\pi_\theta}(s, a)+\sum_{i=1}^N (\lambda_i-\nu_i) A_i^{\pi_\theta}(s, a)\big)=\frac{1}{1-\gamma}\sum_{i=1}^N (\frac{1}{N}+\lambda_i-\nu_i) A_i^{\pi_\theta}(s, a).\label{eq:NPG_Lagrangian}
\end{align}

Formally presented in Algorithm~\ref{Alg:MT-PDNPG-C}, the multi-task primal-dual natural policy gradient (MT-PDNPG) algorithm proposed for the centralized setting ascends the primal variable according to Eq.~\eqref{Alg:MT-PDNPG-C:policy_update} in the direction of the natural gradient, with dual iterates $\lambda_i^k,\nu_i^k$ plugged in as $\lambda_i,\nu_i$ and the Q function as a proxy of the advantage function. On the other hand, dual variables $\lambda_i^k,\nu_i^k$ are iteratively refined according to Eq.~\eqref{Alg:MT-PDNPG-C:dual_update} with gradient descent.

Despite its simplicity, the iterates of Algorithm~\ref{Alg:MT-PDNPG-C} efficiently converge to the globally optimal solution of Eq.~\eqref{eq:objective}. As this algorithm can be regarded as a special case of that in \cite{ding2020natural} under a specific choice of reward functions, its analysis directly follows as a corollary of \cite{ding2020natural}[Theorem 1], which we state below. 
In the following sections, we will generalize Algorithm~\ref{Alg:MT-PDNPG-C} along three dimensions to make it 1) sample-based, 2) compatible with decentralized learning, and 3) effective under linear function approximation.

\begin{algorithm}[!h]
\caption{Multi-Task Primal-Dual Natural Policy Gradient Algorithm (Centralized)}
\label{Alg:MT-PDNPG-C}
\begin{algorithmic}[1]
\STATE{\textbf{Initialization:} Initialize $\theta^0\in\mathbb{R}^{|\Scal||\Acal|}=0$ and for each task $i$ dual variables $\lambda_i^0,\nu_i^0\in\mathbb{R}_+=0$}\\
\FOR{$k=0,1,\cdots,K-1$}
\STATE{Policy update: 
    \begin{align}
        \theta^{k+1} = \alpha\sum_{j=1}^{N}(\frac{1}{N}+\lambda_j^k-\nu_j^k)Q_j^{\pi_{\theta^k}},
        \quad \pi^{k+1}(a\mid s)=\frac{\exp \left(\theta^{k+1}(s,a)\right)}{\sum_{a' \in \Acal} \exp\left(\theta^{k+1}(s,a')\right)}\label{Alg:MT-PDNPG-C:policy_update}
    \end{align}
}
\FOR{Each task $i=1,\cdots,N$}
\STATE{Dual variable update:
    \begin{align}
        \begin{aligned}
        &\lambda_{i}^{k+1}\hspace{-2pt}=\hspace{-2pt}\Pi_{[0,B_{\lambda}]}\Big(\hspace{-1pt}\lambda_{i}^{k}\hspace{-1pt}-\hspace{-1pt}\eta\big(V_i^{\pi_{\theta_i^k}}\hspace{-1pt}(\rho)\hspace{-1pt}-\ell_i\big)\hspace{-2pt}\Big)\\
        &\nu_{i}^{k+1}\hspace{-2pt}=\hspace{-2pt}\Pi_{[0,B_{\lambda}]}\Big(\hspace{-1pt}\lambda_{i}^{k}\hspace{-1pt}+\hspace{-1pt}\eta\big(V_i^{\pi_{\theta_i^k}}\hspace{-1pt}(\rho)\hspace{-1pt}-\hspace{-1pt}u_i\big)\hspace{-2pt}\Big)\end{aligned}\label{Alg:MT-PDNPG-C:dual_update}
    \end{align}}
\ENDFOR
\ENDFOR
\end{algorithmic}
\end{algorithm}

\begin{cor}\label{cor:pg}
Consider the iterates $\{\pi^k\}$ obtained from $K$ iterations of Algorithm \ref{Alg:MT-PDNPG-C} in the centralized setting. Let the step size sequences be $\alpha=\frac{\alpha_0}{K^{1/2}}, \eta=\frac{\eta_0}{K^{1/2}}$,
with $\alpha_0,\eta_0>0$.
Under Assumption~\ref{assump:Slater}, we have 
\begin{align*}
    \max\Big\{\frac{1}{K}\hspace{-2pt}\sum_{k=0}^{K-1}(V_{0}^{\pi^{\star}}(\rho)-V_{0}^{\pi^k}(\rho)),\frac{1}{K}\hspace{-2pt}\sum_{k=0}^{K-1}\sum_{i=1}^{N}\left(\big[\ell_i-V_{i}^{\pi^k}(\rho)\big]_{+}\hspace{-2pt}+\hspace{-2pt}\big[V_{i}^{\pi^k}(\rho)-u_i\big]_{+}\right)\Big\}\leq\Ocal\Big(K^{-1/2}\Big).
\end{align*}
\end{cor}
Corollary~\ref{cor:pg} guarantees that the policy iterate converges in both objective function value and constraint violation with a rate of $\Ocal(1/\sqrt{K})$. As the centralized setting is a special decentralized case where every agent may communicate with every other agent, this result can also be regarded as a corollary of Theorem~\ref{thm:main_pg} to be presented in Section~\ref{sec:decentralized:pg}, which analyzes the MT-PDNPG algorithm under decentralized computation.

\section{Sample-Based Setting: Online Primal-Dual Natural Actor-Critic Algorithm}\label{sec:AC}

The policy gradient algorithm presented in the previous section employs deterministic gradient updates, which requires evaluating the Q function. In large and/or unknown environments, the Q function cannot be exactly computed instantaneously with samples, which makes Algorithm~\ref{Alg:MT-PDNPG-C} inapplicable in practice. In Algorithm~\ref{Alg:MT-PDNAC-C}, we introduce a sample-based extension of Algorithm~\ref{Alg:MT-PDNPG-C}, where the key extension is to maintain a variable $\widehat{Q}_i^k\in\mathbb{R}^{|\Scal||\Acal|}$ as an estimate of $Q_i^{\pi^k}$ in each task $i$, refined over time using samples of the state transition. 
Instead of relying on the true value function, the updates of the policy and dual variable in Eqs.~\eqref{Alg:MT-PDNAC:policy_update} and \eqref{Alg:MT-PDNAC:dual_update} use approximate gradients obtained by plugging in the most up-to-date $\widehat{Q}_i^k$ estimate.
Algorithms of this flavor are usually categorized as actor-critic methods, where the actor refers to the policy iterates and the critic is the value function estimator. We stress that our actor-critic algorithm is single-loop and truly online in the sense that the samples are generated continuously from a single trajectory and updates both actor and critic variables on the run, which makes it very convenient to implement in practice.


While we can generate samples according to the current policy iterate $\pi^k$, doing so may cause certain actions to be very infrequently selected, which leads to exploration issues.
To guarantee sufficient exploration (visitation of all state-action pairs), we introduce a behavior policy $\widehat{\pi}^k$ in Eq.~\eqref{Alg:MT-PDNAC:behaviorpolicy_update}, which is the $\epsilon$-exploration version of $\pi^k$. We note that $\epsilon$ needs to be properly selected with respect to the desired precision and controls an important trade-off: an excessively large $\epsilon$ facilitates the exploration of all state-action pairs but lead to a substantial gap between $\widehat{\pi}^k$ and $\pi^k$, and vice versa. We note that similar behavior policies are also adopted in \cite{borkar2005actor,khodadadian2022finite}.

In the following subsection, we show that this simple, intuitive, completely sample-based algorithm is guaranteed to converge efficiently.

\begin{algorithm}[!h]
\caption{Multi-Task Primal-Dual Natural Actor-Critic Algorithm (Centralized)}
\label{Alg:MT-PDNAC-C}
\begin{algorithmic}[1]
\STATE{\textbf{Initialization:} Initialize $\theta^0\in\mathbb{R}^{|\Scal||\Acal|}=0$ and uniform behavior policy $\widehat{\pi}^0$. For each task $i$, initialize dual variables $\lambda_i^0,\nu_i^0\in\mathbb{R}_+=0$ and critic parameters $\widehat{Q}_i^0 = 0\in\mathbb{R}^{d}$. For each task $i$, draw the initial sample $s_i^{0}$ and $a_i^{0}\sim\widehat{\pi}_i^{0}(\cdot\mid s_i^{0})$}\\
\FOR{$k=0,1,\cdots,K-1$}
\STATE{1) Policy (actor) update: $\forall s\in\Scal, a\in\Acal$
    \begin{align}
        \theta^{k+1} = \alpha\sum_{j=1}^{N}(\frac{1}{N}+\lambda_j^k-\nu_j^k)\widehat{Q}_j^{k},\quad\pi^{k+1}(a\mid s)=\frac{\exp \left(\theta^{k+1}(s,a)\right)}{\sum_{a' \in \Acal} \exp\left(\theta^{k+1}(s,a')\right)}\label{Alg:MT-PDNAC-C:policy_update}
    \end{align}}

\STATE{2) Behavior policy update: $\forall s\in\Scal,a\in\Acal$
    \begin{align}
        \widehat{\pi}^{k+1}(a\mid s)=\frac{\epsilon}{|\Acal|}+\left(1-\epsilon\right) \pi^{k+1}(a\mid s).
        \label{Alg:MT-PDNAC-C:behaviorpolicy_update}
    \end{align}}\\
\FOR{Each task $i=1,\cdots,N$}
\STATE{1) Observe $s_i^{k+1}\sim P(\cdot\mid s_i^{k}, a_i^{k})$ and take action $a_i^{k+1}\sim\widehat{\pi}^{k}(\cdot\mid s_i^{k+1})$}\\
\STATE{2) Value function estimator (critic) update:
    \begin{align}
        \widehat{Q}_i^{k+1}(s_i^k, a_i^k)&=(1-\beta)\widehat{Q}_{i}^{k}(s_i^k, a_i^k) + \beta\left(r_i\left(s_i^k, a_i^k\right)+\gamma \widehat{Q}_i^{k}\left(s_i^{k+1}, a_i^{k+1}\right)\right)\label{Alg:MT-PDNAC-C:critic_update} 
    \end{align}}\\
\STATE{3) Local dual variable update:
    \begin{align}
        \begin{aligned}
        &\lambda_i^{k+1}=\Pi_{[0,B_{\lambda}]}\Big(\lambda_i^k-\eta\big(\sum_{s,a}\rho(s)\pi^k(a\mid s)\widehat{Q}_i^k(s,a)-\ell_i\big)\hspace{-1pt}\Big)\\
        &\nu_i^{k+1}=\Pi_{[0,B_{\lambda}]}\Big(\nu_i^k + \eta\big(\sum_{s,a}\rho(s)\pi^k(a\mid s)\widehat{Q}_i^k(s,a)-u_i\big)\hspace{-1pt}\Big)
        \end{aligned}\label{Alg:MT-PDNAC-C:dual_update}
    \end{align}}
\ENDFOR
\ENDFOR
\end{algorithmic}
\end{algorithm}


\subsection{Finite-Sample Complexity}\label{sec:thm_ac}

We start by introducing some notations and stating our main assumptions.

Given a policy $\pi:\Scal\rightarrow\Delta_{\Acal}$, we use $P^{\pi}:\Scal\rightarrow\Delta_{\Scal}$ to denote the state transition probability under $\pi$
\[P^{\pi}(s'\mid s) = \sum_{a\in\Acal}\Pcal(s'\mid s,a)\pi(a\mid s),\]
and $\mu_{\pi}\in\Delta_{\Scal}$ to denote the stationary distribution of the Markov chain induced by $P^{\pi}$.
In addition, we denote by $\widetilde{\mu}_{\pi}\in\Delta_{\Scal\times\Acal}$ the stationary distribution of state-action pairs under $\pi$, which relates to $\mu_{\pi}$ as follows
\begin{align*}
    \widetilde{\mu}_{\pi}(s,a) = \mu_{\pi}(s)\pi(a\mid s).
\end{align*}

To measure the amount of time a Markov chain takes to approach its stationary distribution, we define the mixing time as follows.
\begin{definition}\label{def:mixing_time}
Given policy $\pi$, consider the Markov chain $\{s^k\}$ generated according to $s^{k+1}\sim P^{\pi}(\cdot\mid s^k)$.
For any scalar $c>0$, the mixing time of $\{s^{k}\}$ associated with $c$ is
\begin{align}
    \tau_{\pi}(c) &\triangleq\min\{k\geq0:\sup_{s\in\Scal}d_{\text{TV}}(P(s^k=\cdot\,|s^0=s),\mu_{\pi}(\cdot))\leq c\},\label{def:mixing_time:eq1}
\end{align}
where given two probability distributions $u_1$ and $u_2$, $d_{TV}$ denotes their total variation distance
\begin{equation}
    d_{\text{TV}}(u_1,u_2)=\frac{1}{2} \sup _{\nu: \Xcal \rightarrow[-1,1]}\left|\int \nu d u_1-\int \nu d u_2\right|.
    \label{eq:TV_def}
\end{equation}
\end{definition}

We consider the following assumption on the the transition probability kernel $\Pcal$ which states that the Markov chain induced by any fixed policy approaches its stationary distribution geometrically fast. 
\begin{assump}[Uniform Ergodicity]\label{assump:markov-chain}
Given any $\pi$, the Markov chain $\{s^k\}$ generated by $P^{\pi}$ according to $s^{k+1}\sim P^{\pi}(\cdot\mid s^k)$ has a unique stationary distribution $\mu_{\pi}$ and is uniformly geometrically ergodic, i.e., there exist $C_{0}\geq 1$ and $\ell\in (0,1)$ such that
\begin{equation*}
\sup_{s}d_{\text{TV}}(P(s^k=\cdot \mid s^0=s),\mu_{\pi}(\cdot))\leq C_{0}\ell^k,\, \forall k\geq 0.
\end{equation*}
\end{assump}

This assumption is commonly made in RL and optimization papers that study gradient-based algorithms under samples collected from a (time-varying) Markovian chain \cite{wu2020finite,khodadadian2022finite,zeng2021two}. Recall the mixing time $\tau_\pi(c)$ defined in Eq.~\eqref{def:mixing_time:eq1}. As an obvious result of Assumption~\ref{assump:markov-chain}, there exists a constant $D>0$ such that
\begin{equation}
\tau_{\pi}(c) \leq D\log(1/c),\quad\forall c\in(0,1)\,\, \text{and}\,\,\pi.\label{eq:mixing:tau}
\end{equation} 
In this work, we denote
\begin{align}
    \tau(c)\triangleq\max_{\pi}\tau_{\pi}(c)\leq D\log(1/c).\label{eq:tau_def}
\end{align}

Another consequence of the uniform ergodicity assumption is that the stationary distribution $\mu_{\pi}$ is uniformly bounded away from 0 for any policy $\pi$, and we denote
$\underline{\mu}\triangleq\min_{\pi,s}\mu_{\pi}(s)>0$.

Algorithm~\ref{Alg:MT-PDNAC-C} is guaranteed to converge to the optimal multi-task policy $\pi^{\star}$ under proper choice of the step sizes. We state the result below, which we take to be a corollary of the analysis of the decentralized MT-PDNAC algorithm to be presented in the next section.

\begin{cor}\label{cor:ac}
Consider the iterates $\{\pi^k\}$ obtained from $K$ iterations of Algorithm \ref{Alg:MT-PDNAC-C} in the centralized setting. Let the step size sequences be
\begin{align*}
   \alpha=\frac{\alpha_0}{K^{5/6}}, \,\,\beta=\frac{\beta_0}{K^{1/2}},\,\,\eta=\frac{\eta_0}{K^{5/6}},\,\,\epsilon=\frac{\epsilon_0}{K^{1/6}},
\end{align*}
with $\frac{(1-\gamma)\underline{\mu}\epsilon_0\beta_0}{|\Acal|}\leq 1$ and $\alpha_0=\Ocal(N^{-1/4})$.
Then, under Assumptions~\ref{assump:Slater} and \ref{assump:markov-chain}, we have
\begin{align*}
    \max\Big\{\frac{1}{K}\sum_{k=0}^{K-1}(V_{0}^{\pi^{\star}}\hspace{-2pt}(\rho)\hspace{-2pt}-\hspace{-2pt}V_{0}^{\pi^{k}}\hspace{-1pt}(\rho)),\frac{1}{K}\sum_{k=0}^{K-1}\sum_{i=1}^{N}\left(\big[\ell_i-V_{i}^{\pi^k}(\rho)\big]_{+}+\big[V_{i}^{\pi^k}(\rho)-u_i\big]_{+}\right)\hspace{-1pt}\Big\}&\\
    \leq\Ocal\Big(\frac{N^{5/4}\log(K)}{K^{1/6}}\Big).&
\end{align*}
\end{cor}

Corollary~\ref{cor:ac} shows a finite-time complexity of $\widetilde{\Ocal}(K^{-1/6})$, where $\widetilde{\Ocal}$ hides structural constants and logarithm factors.
As Algorithm~\ref{Alg:MT-PDNAC-C} draws exactly $N$ samples in each iteration, it implies that the algorithm will converge in both objective function and constraint violation up to a precision $\delta$ with at most $\widetilde{\Ocal}(\delta^{-6})$ samples. This result matches the best-known rate of actor-critic algorithms for solving the policy optimization problem under a single-task, unconstrained MDP \cite{khodadadian2022finite}.
Compared with the complexity of Algorithm~\ref{Alg:MT-PDNPG-C}, we have a slightly inferior rate for not exactly knowing the Q function of the current policy iterates.

\section{Algorithms Under Decentralized Computation}

In many practical applications (for example, each task is associated with an environment physically separated from each other), the information of the multiple tasks may not all be available on a central server. A more realistic learning paradigm in such scenarios is to employ a network of $N$ agents, each placed to explore and learn in one different environment. For generality, we do not assume the existence of a central coordinator that exchanges information with all agents. Rather, the agents are connected according to an undirected graph $\Gcal=(\Vcal,\Ecal)$; agent $i\in\Vcal$ can communicate with agent $j\in\Ncal_i$, where $\Ncal_i=\{j\in\Vcal :(i,j)\in\Ecal\}$ denotes the neighbors of agent $i$. Since each agent can only access local information, collaboration between the agents is necessary to solve the global objective in Eq.~\eqref{eq:objective}.

\subsection{Deterministic Gradient Setting}\label{sec:decentralized:pg}
We first present the development of the decentralized MT-PDNPG method, which extends Algorithm~\ref{Alg:MT-PDNPG-C}. While we would still like to perform alternating gradient descent ascent on the dual and primal variables as in Eqs.~\eqref{Alg:MT-PDNPG-C:dual_update} and \eqref{Alg:MT-PDNPG-C:policy_update}, the challenge is that the natural gradient in Eq.~\eqref{Alg:MT-PDNPG-C:policy_update} involves information across all tasks and obviously cannot be computed by any single agent locally. 

\begin{algorithm}[!h]
\caption{Multi-Task Primal-Dual Natural Policy Gradient Algorithm (Decentralized)}
\label{Alg:MT-PDNPG}
\begin{algorithmic}[1]
\STATE{\textbf{Initialization:} Each agent $i$ initializes $\theta_i^0\in\mathbb{R}^{|\Scal||\Acal|}=0$ and dual variables $\lambda_i^0,\nu_i^0\in\mathbb{R}_+=0$}\\
\FOR{$k=0,1,\cdots,K-1$}
\FOR{Each task $i=1,\cdots,N$}
\STATE{1) Policy update: 
    \begin{gather}
        \theta_i^{k+1} =  \sum_{j\in\mathcal{N}_i}W_{i,j}\theta_{j}^{k} + \alpha (\frac{1}{N}+\lambda_i^k-\nu_i^k)Q_i^{\pi_{\theta_i^k}}\label{Alg:MT-PDNPG:policy_update}\\
        \pi_i^{k+1}(a\mid s)\hspace{-2pt}=\hspace{-2pt}\frac{\exp \left(\theta_i^{k+1}(s,a)\right)}{\sum_{a' \in \Acal} \exp\left(\theta_i^{k+1}(s,a')\right)}\notag
    \end{gather}
    }
\STATE{2) Dual variable update:
    \begin{align}
        \begin{aligned}
        &\lambda_{i}^{k+1}\hspace{-2pt}=\hspace{-2pt}\Pi_{[0,B_{\lambda}]}\Big(\hspace{-1pt}\lambda_{i}^{k}\hspace{-1pt}-\hspace{-1pt}\eta\big(V_i^{\pi_{\theta_i^k}}\hspace{-1pt}(\rho)\hspace{-1pt}-\ell_i\big)\hspace{-2pt}\Big)\\
        &\nu_{i}^{k+1}\hspace{-2pt}=\hspace{-2pt}\Pi_{[0,B_{\lambda}]}\Big(\hspace{-1pt}\lambda_{i}^{k}\hspace{-1pt}+\hspace{-1pt}\eta\big(V_i^{\pi_{\theta_i^k}}\hspace{-1pt}(\rho)\hspace{-1pt}-\hspace{-1pt}u_i\big)\hspace{-2pt}\Big)\end{aligned}\label{Alg:MT-PDNPG:dual_update}
    \end{align}}
\ENDFOR
\ENDFOR
\end{algorithmic}
\end{algorithm}

Our solution 
to the problem is to make each agent move in the direction of a locally computable portion of this gradient, followed by an averaging step (weighted according to matrix $W\in\mathbb{R}^{N\times N}$) that mixes the agents' policy parameters to achieve consensus. 
This specific update is shown in Algorithm~\ref{Alg:MT-PDNPG} Eq.~\eqref{Alg:MT-PDNPG:policy_update}.
In the long run, the local policy parameter obtained by following Eq.~\eqref{Alg:MT-PDNPG:policy_update} behaves almost as if each agent is updated using the global aggregate gradient in Eq.~\eqref{eq:NPG_Lagrangian}.

To analyze the complexity of MT-PDNPG, we make the following assumption on the matrix $W$, which specifies the averaging weight in Eq.~\eqref{Alg:MT-PDNPG:policy_update}.

\begin{assump}\label{assump:W_matrix}
    The matrix $W$ is doubly stochastic, i.e. $\sum_{i=1}^{N}W_{i,j}=\sum_{j=1}^{N}W_{i,j}=1$. In addition, $W_{i,j}>0$ if and only if $(i,j)\in\Ecal$ and $W_{i,j}=0$ otherwise.
\end{assump}

This assumption is standard in the literature of consensus optimization \cite{yuan2016convergence,zhang2018fully,zeng2021decentralized,zeng2022finite2}. For any connected communication graph $\Gcal$, a weight matrix $W$ that satisfies Assumption~\ref{assump:W_matrix} can be simply found using the lazy Metropolis method \cite{olshevsky2015linear}.
The largest singular value of $W$ is always 1, and we use $\sigma_2(W)\in[0,1)$ to denote its second largest singular value. In general, a more densely connected graph $\Gcal$ leads to a smaller $\sigma_2(W)$.
We now present the first main theoretical result of the paper, which guarantees the finite-time convergence of Algorithm~\ref{Alg:MT-PDNPG}.

\begin{thm}\label{thm:main_pg}
Consider the iterates $\{\pi_i^k\}$ obtained from $K$ iterations of Algorithm \ref{Alg:MT-PDNPG}. Let the step size sequences be
\begin{align*}
   \alpha=\frac{\alpha_0}{K^{1/2}}, \quad\eta=\frac{\eta_0}{K^{1/2}},
\end{align*}
with $\alpha_0=\Ocal(\sqrt{1-\sigma_2(W)})$.
Then, under Assumption~\ref{assump:Slater} and \ref{assump:W_matrix}, we have for any $j=1,\cdots,N$
\begin{align*}
    \max\Big\{\frac{1}{K}\sum_{k=0}^{K-1}(V_{0}^{\pi^{\star}}(\rho)-V_{0}^{\pi_j^k}(\rho)),\frac{1}{K}\sum_{k=0}^{K-1}\sum_{i=1}^{N}\left(\big[\ell_i-V_{i}^{\pi_j^k}(\rho)\big]_{+}\hspace{-2pt}+\hspace{-2pt}\big[V_{i}^{\pi_j^k}(\rho)-u_i\big]_{+}\right)\Big\}&\\
    \leq\Ocal\Big(\frac{N^{5/4}}{\sqrt{1-\sigma_2(W)}K^{1/2}}\Big).&
\end{align*}
\end{thm}

Theorem~\ref{thm:main_pg} states that the policy iterate at every local agent $i$ converges to the globally optimal multi-task policy $\pi^{\star}$ in both objective function and constraint violation with rate $\Ocal(1/\sqrt{K})$, which matches the complexity of the algorithm in the centralized setting. 
The dependency of the bound on $N$ reflects the increasing difficulty as the number of tasks scales up, while the inverse dependency on $\sqrt{1-\sigma_2(W)}$ captures the impact of the communication graph $\Gcal$ and matrix $W$.

\subsection{Sample-Based Setting}

The same principle of distributing computation across the network can be applied to extend Algorithm~\ref{Alg:MT-PDNAC-C} to the decentralized setting. Policy averaging is performed by each agent with its neighbors, while the critic and dual variables are updated completely locally. We present the updates formally in Algorithm~\ref{Alg:MT-PDNAC} and study its finite-time complexity below.

\begin{algorithm}[!h]
\caption{Multi-Task Primal-Dual Natural Actor-Critic Algorithm (Decentralized)}
\label{Alg:MT-PDNAC}
\begin{algorithmic}[1]
\STATE{\textbf{Initialization:} Each agent $i$ initializes $\theta_i^0\in\mathbb{R}^{|\Scal||\Acal|}=0$, dual variables $\lambda_i^0,\nu_i^0\in\mathbb{R}_+=0$, critic parameters $\widehat{Q}_i^0 = 0\in\mathbb{R}^{d}$, and uniform behavior policy $\widehat{\pi}_i^0$. For each task $i$, draw the initial sample $s_i^{0}$ and $a_i^{0}\sim\widehat{\pi}_i^{0}(\cdot\mid s_i^{0})$}\\
\FOR{$k=0,1,\cdots,K-1$}
\FOR{Each task $i=1,\cdots,N$}
\STATE{1) Observe $s_i^{k+1}\sim P(\cdot\mid s_i^{k}, a_i^{k})$ and take action $a_i^{k+1}\sim\widehat{\pi}_i^{k}(\cdot\mid s_i^{k+1})$}\\
\STATE{2) Value function estimator (critic) update:
    \begin{align}
        \widehat{Q}_i^{k+1}(s_i^k, a_i^k)&=(1-\beta)\widehat{Q}_{i}^{k}(s_i^k, a_i^k) + \beta\left(r_i\left(s_i^k, a_i^k\right)+\gamma \widehat{Q}_i^{k}\left(s_i^{k+1}, a_i^{k+1}\right)\right)\label{Alg:MT-PDNAC:critic_update} 
    \end{align}}\\
\STATE{3) Policy (actor) update: $\forall s\in\Scal, a\in\Acal$
    \begin{gather}
        \theta_i^{k+1} = \sum_{j\in\mathcal{N}_i}W_{i,j}\theta_{j}^{k} + \alpha (\frac{1}{N}+\lambda_i^k-\nu_i^k)\widehat{Q}_i^{k}\label{Alg:MT-PDNAC:policy_update}\\
        \pi_i^{k+1}(a\mid s)=\frac{\exp \left(\theta_i^{k+1}(s,a)\right)}{\sum_{a' \in \Acal} \exp\left(\theta_i^{k+1}(s,a')\right)}\notag
    \end{gather}}
\STATE{4) Behavior policy update: $\forall s\in\Scal,a\in\Acal$
    \begin{align}
        \widehat{\pi}_i^{k+1}(a\mid s)=\frac{\epsilon}{|\Acal|}+\left(1-\epsilon\right) \pi_i^{k+1}(a\mid s).
        \label{Alg:MT-PDNAC:behaviorpolicy_update}
    \end{align}}\\
\STATE{5) Local dual variable update:
    \begin{align}
        \begin{aligned}
        &\lambda_i^{k+1}=\Pi_{[0,B_{\lambda}]}\Big(\lambda_i^k-\eta\big(\sum_{s,a}\rho(s)\pi_i^k(a\mid s)\widehat{Q}_i^k(s,a)-\ell_i\big)\hspace{-1pt}\Big)\\
        &\nu_i^{k+1}=\Pi_{[0,B_{\lambda}]}\Big(\nu_i^k + \eta\big(\sum_{s,a}\rho(s)\pi_i^k(a\mid s)\widehat{Q}_i^k(s,a)-u_i\big)\hspace{-1pt}\Big)
        \end{aligned}\label{Alg:MT-PDNAC:dual_update}
    \end{align}}
\ENDFOR
\ENDFOR
\end{algorithmic}
\end{algorithm}

\begin{thm}\label{thm:main_ac}
Consider the iterates $\{\pi_i^k\}$ from $K$ iterations of Algorithm \ref{Alg:MT-PDNAC}. Let the step size sequences be
\begin{align*}
   \alpha=\frac{\alpha_0}{K^{5/6}}, \,\,\beta=\frac{\beta_0}{K^{1/2}},\,\,\eta=\frac{\eta_0}{K^{5/6}},\,\,\epsilon=\frac{\epsilon_0}{K^{1/6}},
\end{align*}
with $\frac{(1-\gamma)\underline{\mu}\epsilon_0\beta_0}{|\Acal|}\leq 1$ and $\alpha_0=\Ocal(N^{-1/4}\sqrt{1-\sigma_2(W)})$.
Then, under Assumptions~\ref{assump:Slater}-\ref{assump:W_matrix}, we have for any local agent $j$
\begin{align*}
    \max\Big\{\frac{1}{K}\sum_{k=0}^{K-1}(V_{0}^{\pi^{\star}}\hspace{-2pt}(\rho)\hspace{-2pt}-\hspace{-2pt}V_{0}^{\pi_j^{k}}\hspace{-1pt}(\rho)),\frac{1}{K}\sum_{k=0}^{K-1}\sum_{i=1}^{N}\left(\big[\ell_i-V_{i}^{\pi_j^k}(\rho)\big]_{+}+\big[V_{i}^{\pi_j^k}(\rho)-u_i\big]_{+}\right)\hspace{-1pt}\Big\}&\\
    \leq\Ocal\Big(\frac{N^{5/4}\log(K)}{\sqrt{1-\sigma_2(W)}K^{1/6}}\Big).&
\end{align*}
\end{thm}

Theorem~\ref{thm:main_ac} generalizes the result in Section~\ref{sec:AC} and guarantees that the decentralized MT-PDNAC algorithm finds the globally optimally policy with a convergence rate of $\widetilde{O}(K^{-1/6})$. As Algorithm~\ref{Alg:MT-PDNAC} again draws exactly $N$ sample in each iteration, Theorem~\ref{thm:main_ac} implies a sample complexity of $\widetilde{\Ocal}(\delta^{-6})$ for converging to an optimal solution up to a precision $\delta$.

A main challenge to the analysis of Algorithm~\ref{Alg:MT-PDNAC} lies the coupling between the actor, critic, and dual variables, further complicated by the access to only local information. The actor controls the behavior policy which generates the samples for the update of the critic, which in turn is involved in the dual variable update. In addition, the accuracy of the critic, along with that of the dual variable, affects the update of the actor. To handle this intertwined system of variables, we leverage the fact that the effect of the samples on the actor (and dual variable) is indirect through the critic. This enables us to isolate the analysis of the actor from the evolution of the Markov chain conditioning on the critic. Specifically, as key components of the analysis, we show that the sub-optimality gap converge up to the cumulative error of the critic while the critic converges linearly fast to find the value function of the behavior policy under evolving Markovian samples. 
Further details of the analysis can be found in Section~\ref{sec:thm:main_pg} of the appendix.

\section{Sample-Based Algorithm Under Linear Function Approximation}

When we take an optimization approach to solve Eq.~\eqref{eq:objective} in the previous sections, the policy and Q function are both $|\Scal||\Acal|$-dimensional objects. This makes optimizing Eq.~\eqref{eq:objective} (or even evaluating the objective function for a given policy) prohibitively expensive in real-life problems where the state space is huge or even infinitely large.

To reduce the dimensionality of the problem, we consider parameterizing the policy and Q function with linear function approximation. Suppose that the optimal policy parameter and Q functions of all tasks are be (approximately) represented using a given set of $d$ basis vectors $\{\phi_1,\dots,\phi_d\in\mathbb{R}^{|\Scal||\Acal|}\}$ where $d$ can be much smaller than $|\Scal||\Acal|$. Each state and action pair $(s,a)$ is associated with the feature vector $\phi(s,a)=[\phi_1(s,a),\dots,\phi_d(s,a)]\in\mathbb{R}^d$. We define the feature matrix $\Phi\in\mathbb{R}^{|\Scal||\Acal|\times d}$ as
\begin{align*}
\Phi=\left[\begin{array}{ccc}
\mid & & \mid \\
\phi_1 & \cdots & \phi_d \\
\mid & & \mid
\end{array}\right]=\left[\begin{array}{ccc}
- & \phi(s_0,a_0)^{\top} & - \\
\cdots & \cdots & \cdots \\
- & \phi(s_{|\Scal|},a_{|\Acal|})^{\top} & -
\end{array}\right].
\end{align*}
Without loss of generality, we assume that $\|\phi(s,a)\|\leq1$ for all $s\in\Scal,a\in\Acal$ and $\Phi$ is full column rank, and use $\sigma_{\min}(\Phi)$ and $\sigma_{\max}(\Phi)$ to denote the smallest and largest singular value of $\Phi$, respectively.

Given the feature vectors, we maintain a parameter $\theta\in\mathbb{R}^d$ which represents the policy through the log-linear parameterization
\begin{align}
    \pi_{\theta}(a \mid s)=\frac{\exp(\phi(s, a)^{\top} \theta)}{\sum_{a'\in \Acal} \exp(\phi(s, a')^{\top} \theta)}.\label{eq:policy_loglinear}
\end{align}

We use the same features to approximate the Q function and restrict it to lie in the $d$-dimensional subspace $\widehat{\Qcal}=\{\Phi\omega:\omega\in\mathbb{R}^d\}$. Given any policy $\pi$, its Q function in task $i$ satisfies the Bellman equation
\[Q_i^{\pi}=T_i^{\pi}Q_i^{\pi},\]
where $T_i^{\pi}:\mathbb{R}^{|\Scal||\Acal|}\rightarrow\mathbb{R}^{|\Scal||\Acal|}$ is the Bellman backup operator
\[(T_i^{\pi}Q_i^{\pi})(s,a)=r_i(s,a)+\gamma\hspace{-2pt}\sum_{a'\in\Acal}\hspace{-2pt}P(s'\mid s,a)\pi(a'\mid s')Q_i^{\pi}(s,a).\]
Under linear function approximation the feature matrix $\Phi$ may not exactly span $Q_i^{\pi}$. In this work, we solve the projected Bellman equation. Letting $M^{\pi}\in\mathbb{R}^{|\Scal||\Acal|\times|\Scal||\Acal|}$ denote the diagonal matrix such that
\begin{align}
    M_{(s,a),(s,a)}=\widetilde{\mu}_{\pi}(s,a)=\mu_{\pi}(s)\pi(a\mid s),\label{eq:def_M}
\end{align}
we define the operator $\Pi_{\Phi}^{\pi}:\mathbb{R}^{|\Scal||\Acal|}\rightarrow\mathbb{R}^{|\Scal||\Acal|}$ as the projection to the linear subspace spanned by $\Phi$ under the weighted $\ell_2$ norm $\|\cdot\|_{\pi}^2=\cdot^{\top}M^{\pi}\cdot$. If the policy $\pi$ is in the interior of the probability simplex (has uniformly lower bounded entries),
the projected Bellman equation takes the following form
\begin{align}
\Phi\omega=\Pi_{\Phi}^{\pi}T_i^{\pi}\Phi\omega.\label{eq:projected_bellman_eq}
\end{align}
A straightforward extension of \cite{tsitsiklis1997analysis}[Theorem 1] guarantees that a unique solution to Eq.~\eqref{eq:projected_bellman_eq} exists, and we use $\omega_i^{\star}:\Delta_{\Acal}^{\Scal}\rightarrow d$ to denote the mapping from a policy to its optimal value function parameter, which is the solution of Eq.~\eqref{eq:projected_bellman_eq}. Note that $\omega_i^{\star}(\pi)$ satisfies
\begin{align}
    \widebar{H}^{\pi}\omega_i^{\star}(\pi)+\widebar{b}_i^{\pi}=0,\label{eq:Aomega=b}
\end{align}
where
\begin{align}
\begin{aligned}
    &\widebar{H}^{\pi}\hspace{-3pt}=\hspace{-2pt}\mathbb{E}_{s,a\sim\widetilde{\mu}_{\pi},s'\hspace{-1pt}\sim\hspace{-1pt} P(\cdot\mid s,a),a'\hspace{-1pt}\sim\hspace{-1pt}\pi(\cdot\mid s')}\hspace{-1pt}[\phi(s,\hspace{-1pt}a)(\gamma\phi(s',a')\hspace{-2pt}-\hspace{-2pt}\phi(s,a))^{\top}]\\
    &\widebar{b}_i^{\pi}=\mathbb{E}_{s,a\sim\widetilde{\mu}_{\pi}}[r_i(s,a)\phi(s,a)].
\end{aligned}
\label{eq:def_Ab}
\end{align}

\subsection{Online Nested-Loop Algorithm}

To extend Algorithm~\ref{Alg:MT-PDNAC} to the linear function approximation setting, a most straightforward approach simply replaces the actor and critic updates in Eqs.~\eqref{Alg:MT-PDNAC:policy_update} and \eqref{Alg:MT-PDNAC:critic_update} with their properly generalized versions (see Eqs.\eqref{Alg:MT-PDNAC_nestedloop:policy_update} and \eqref{Alg:MT-PDNAC_nestedloop:critic_update} of Algorithm~\ref{Alg:MT-PDNAC_nestedloop}). However, we note that the analysis of the single-loop actor-critic algorithm relies critically on the bounded variation of the TD learning target over iterations. In the tabular setting, the TD learning target is $Q_i^{\pi_i^{k}}$, the Q function of the current policy iterate $\pi_i^k$. The shift $\|Q_i^{\pi_i^{k+1}}-Q_i^{\pi_i^k}\|$ can be easily controlled by $\|\pi_i^{k+1}-\pi_i^k\|$ and eventually by the step size $\alpha$, due to the following Lipschitz condition (established in Lemma~\ref{lem:misc})
\begin{align*}
    \|Q_i^{\pi}-Q_i^{\pi'}\|\leq\frac{\gamma|\Scal||\Acal|}{(1-\gamma)^2}\|\pi-\pi'\|,\quad\forall\pi,\pi'.
\end{align*}
In the case of linear function approximation, the TD learning target is $\omega_i^{\star}$, which is nevertheless not a (sufficiently) Lipschitz operator. It can be shown that there exists $\pi,\pi'$ and constant $L>0$ such that
\begin{align*}
    \|\omega_i^{\star}(\pi)-\omega_i^{\star}(\pi')\|\geq\frac{L}{\min_{s,a}\{\pi(a\mid s),\pi'(a\mid s)\}}\|\pi-\pi'\|.
\end{align*}
However, to maintain the sample complexity of $\widetilde{\Ocal}(\delta^{-6})$ of Algorithm~\ref{Alg:MT-PDNAC}, we need at least
\begin{align*}
    \|\omega_i^{\star}(\pi_i^k)-\omega_i^{\star}(\pi_i^{k+1})\|\leq\frac{L}{\left(\min_{s,a}\{\pi_i^k(a\mid s),\pi_i^{k+1}(a\mid s)\}\right)^{1/2}}\|\pi_i^k-\pi_i^{k+1}\|,\quad\forall k.
\end{align*}

This apparent gap in the Lipschitz continuity of $\omega_i^{\star}$ means that we cannot apply the analysis established earlier in the paper. To work around the challenge without degrading the sample complexity, we employ nested-loop updates in the linear function approximation setting. Stated formally in Algorithm~\ref{Alg:MT-PDNAC_nestedloop}, our method updates the policy, dual variable, and behavior policy in the outer loop and allows more iterations for the critic in the inner loop for it to chase the moving target $\omega_i^{\star}(\pi_i^{k})$.

\begin{remark}
We can still use a single-loop algorithm under linear function approximation despite the degraded Lipschitz condition on the TD learning target. By choosing the step sizes slightly differently from the tabular case, the single-loop algorithm will converge in finite time, but the complexity is slightly worse than $\widetilde{\Ocal}(K^{-1/6})$ established in Theorem~\ref{thm:main_ac}.
\end{remark}

Note that Eq.~\eqref{Alg:MT-PDNAC_nestedloop} is equivalent to the following direct update on $\pi_i^k$ \cite{chen2022finite}[Lemma 3.1]
\begin{align}
    &\pi_i^{k+1}(a\mid s)\propto\pi_i^k(a\mid s)\exp(\alpha(1/N+\lambda_i^k-\nu_i^k)\phi(s,a)^{\top})\omega_i^{k,T}.
\end{align}
We present the expressions of $\pi_i^{k+1}$ and $\widehat{\pi}_i^{k+1}$ in Eqs.~\eqref{Alg:MT-PDNAC_nestedloop:policy_update} and \eqref{Alg:MT-PDNAC_nestedloop:behaviorpolicy_update} only for the purpose of clarity. As the policies are large $|\Scal||\Acal|$-dimensional objects, the policies are never explicitly maintained or updated when the algorithm is deployed; we only need to track the policy parameter $\theta_i^k$.

\subsection{Finite-Sample Complexity}

This section presents the complexity of the nested-loop actor-critic algorithm under linear function approximation. As the function approximation is not necessarily perfect (i.e. the true Q functions may not lie in the column space of the feature matrix $\Phi$), the optimality error does not converge to 0 asymptotically but rather depends on the approximation error associated with $\Phi$, which is the distance between $Q_i^{\pi}$ and its approximation in $\widehat{\Qcal}$. In the next assumption, we assume that there is a uniform upper bound on the approximation error.
\begin{assump}\label{assump:epsilon_max}
There exists a constant $\varepsilon_{\max}>0$ such that for any $\pi$, we have
\begin{align}
    \|\Phi\omega_i^{\star}(\pi)-Q_{i}^{\pi}\|\leq\varepsilon_{\max}.
\end{align}
\end{assump}

A consequence of the assumption is that the target critic parameter always has a bounded norm. More specifically, we have $\|\omega_i^{\star}(\widehat{\pi}_i^k)\|\leq B_{\omega}$ for all $k\geq 0$, where
\begin{align}
B_{\omega}=\sigma_{\min}^{-1}(\hspace{-1pt}\Phi\hspace{-1pt})(\sqrt{\frac{|\Scal||\Acal|}{1-\gamma}}+\varepsilon_{\max}).\label{eq:def_Bomega}
\end{align}
This result is stated and proved in Lemma~\ref{lem:omega_bound} of the appendix. We need to confine the possibly infinite growth of the critic parameter $\omega_i^k$ through projection. Knowing the boundedness of $\omega_i^{\star}(\widehat{\pi}_i^k)$ allows us to use the operator $\Pi_{B_{\omega}}:\mathbb{R}^{d}\rightarrow\mathbb{R}^{d}$, which projects a vector to the $\ell_2$ ball with radius $B_{\omega}$. 

\begin{thm}\label{thm:LFA}
Let $\delta>0$ be a desired precision. Under Assumptions~\ref{assump:Slater}-\ref{assump:epsilon_max} and properly selected step sizes, the iterates of Algorithm \ref{Alg:MT-PDNAC_nestedloop} satisfies for any agent $j$
\begin{align*}
    \max\Big\{\frac{1}{K}\sum_{k=0}^{K-1}(V_{0}^{\pi^{\star}}\hspace{-2pt}(\rho)\hspace{-2pt}-\hspace{-2pt}V_{0}^{\pi_j^{k}}\hspace{-1pt}(\rho)),\frac{1}{K}\sum_{k=0}^{K-1}\sum_{i=1}^{N}\left(\big[\ell_i-V_{i}^{\pi_j^k}(\rho)\big]_{+}+\big[V_{i}^{\pi_j^k}(\rho)-u_i\big]_{+}\right)\hspace{-1pt}\Big\}&\notag\\
    \leq\frac{N^{5/4}\delta}{\sqrt{1-\sigma_2(W)}}+\Ocal(N\varepsilon_{\max})&
\end{align*}
with at most $K=\Ocal(\delta^{-2})$ outer loop iterations and at most $\Ocal(\frac{\log(1/\delta)}{\delta^6})$ total samples.
\end{thm}

Theorem~\ref{thm:LFA} again establishes a $\widetilde{\Ocal}(\delta^{-6})$ convergence for every agent's local policy parameter, which scales inversely with $1-\sigma_2(W)$, the spectral gap of the weight matrix. Due to the presence of approximation error, the optimality gap in objective function and constraint violation does not converge to 0, but to a quantity proportional to $\varepsilon_{\max}$. Prior to our work, there are few results on the finite-sample complexity of data-driven algorithms for solving the CMDP optimization problem under linear function approximation, even in the single task setting. 
The recent work \cite{ding2022convergence} targets this problem and studies a primal-dual REINFORCE-flavored algorithm. Their overall sample complexity is $\Ocal(\delta^{-8})$, which we significantly improve over.
Finally, we note that the sample complexity here matches that of Algorithm~\ref{Alg:MT-PDNAC} for the tabular setting.




\begin{algorithm}[!h]
\caption{Multi-Task Primal-Dual Natural Actor-Critic Algorithm under Linear Function Approximation (Decentralized)}
\label{Alg:MT-PDNAC_nestedloop}
\begin{algorithmic}[1]
\STATE{\textbf{Initialization:} For each task $i$, initialize $\theta_i^0\in\mathbb{R}^{|\Scal||\Acal|}=0$, dual variables $\lambda_i^0,\nu_i^0\in\mathbb{R}_+=0$, critic parameters $\widehat{Q}_i^0 = 0\in\mathbb{R}^{d}$, and uniform behavior policy $\widehat{\pi}_i^0$. For each task $i$, draw the initial sample $s_i^{0,0}$ and $a_i^{0,0}\sim\widehat{\pi}_i^{0,0}(\cdot\mid s_i^{0,0})$}\\
\FOR{$k=0,1,\cdots,K-1$}
\FOR{Each task $i=1,\cdots,N$}
\FOR{$t=0,1,\cdots,T-1$}
\STATE{1) Observe $s_i^{k,t+1}\sim P(\cdot\mid s_i^{k,t+1}, a_i^{k,t+1})$ and take action $a_i^{k,t+1}\sim\widehat{\pi}_i^{k}(\cdot\mid s_i^{k,t+1})$}
\STATE{2) Critic parameter update:
    \begin{gather}
        \widehat{\omega}_{i}^{k,t+1} =\omega_i^{k,t}+\beta\phi(s_i^{k,t},a_i^{k,t})\Big(r_i(s_i^{k,t},a_i^{k,t})+(\gamma \phi(s_i^{k,t+1},a_i^{k,t+1})-\phi(s_i^{k,t},a_i^{k,t}))^{\top}\omega_i^{k,t}\Big)\notag\\
        \omega_i^{k,t+1} = \Pi_{B_{\omega}}\left(\widehat{\omega}_i^{k,t+1}\right)\label{Alg:MT-PDNAC_nestedloop:critic_update}
    \end{gather}}
\ENDFOR
\STATE{3) Policy (actor) update: $\forall s\in\Scal, a\in\Acal$
    \begin{align}
        &\theta_i^{k+1} = \sum_{j\in\mathcal{N}_i}W_{i,j}\theta_{j}^{k} + \alpha (\frac{1}{N}+\lambda_i^k-\nu_i^k)\omega_i^{k,T}\notag\\
        &\pi_i^{k+1}(a\mid s)=\frac{\exp \left(\theta_i^{k+1}(s,a)^{\top}\omega_i^{k,T}\right)}{\sum_{a' \in \Acal} \exp\left(\theta_i^{k+1}(s,a')^{\top}\omega_i^{k,T}\right)}
        \label{Alg:MT-PDNAC_nestedloop:policy_update}
    \end{align}
    }
\STATE{4) Behavior policy update: $\forall s\in\Scal,a\in\Acal$
    \begin{align}
        \widehat{\pi}_i^{k+1}(a\mid s)=\frac{\epsilon}{|\Acal|}+\left(1-\epsilon\right) \pi_i^{k+1}(a\mid s).
        \label{Alg:MT-PDNAC_nestedloop:behaviorpolicy_update}
    \end{align}}\\
\STATE{5) Dual variable update:
    \begin{align}
        \begin{aligned}
        &\lambda_i^{k+1}=\Pi_{[0,B_{\lambda}]}\Big(\lambda_i^k-\eta\big(\sum_{s,a}\rho(s)\pi_i^k(a\mid s)\phi(s,a)^{\top}\omega_i^{k}-\ell_i\big)\hspace{-1pt}\Big)\\
        &\nu_i^{k+1}=\Pi_{[0,B_{\lambda}]}\Big(\nu_i^k +\eta\big(\sum_{s,a}\rho(s)\pi_i^k(a\mid s)\phi(s,a)^{\top}\omega_i^{k}-u_i\big)\hspace{-1pt}\Big)
        \end{aligned}\label{Alg:MT-PDNAC_nestedloop:dual_update}
    \end{align}}\\
\STATE{6) Set $s_i^{k+1,0}=s_i^{k,T}$ and $a_i^{k+1,0}=a_i^{k,T}$}
\ENDFOR
\ENDFOR
\end{algorithmic}
\end{algorithm}

\section{Numerical Simulations}\label{sec:simulations}

\begin{figure}[ht]
\centering
\includegraphics[width=\linewidth]{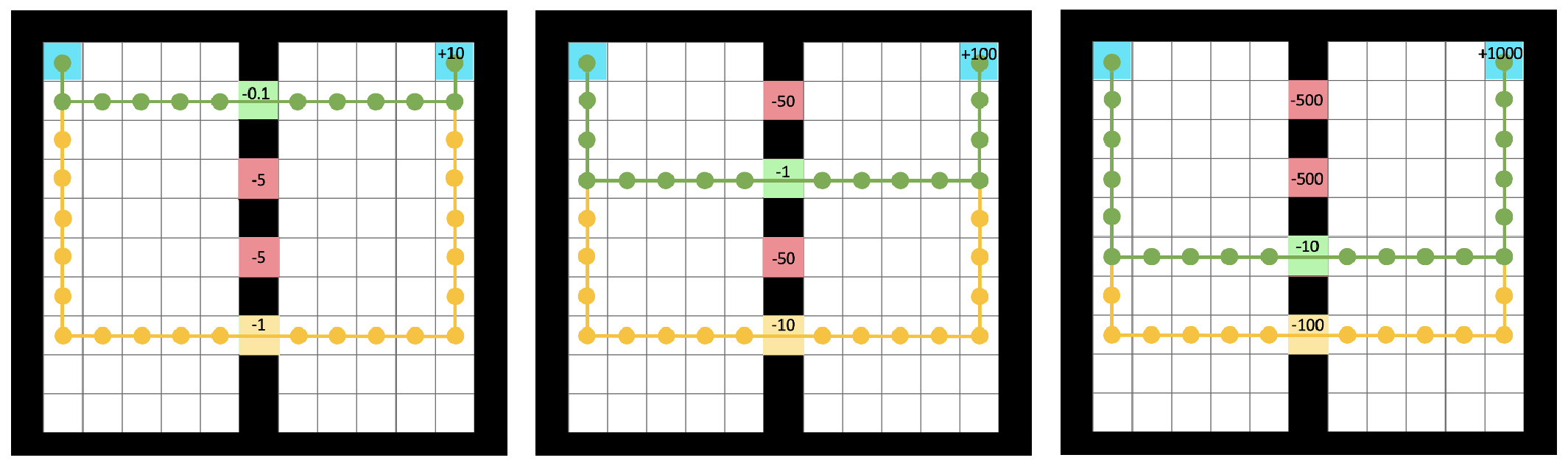}
\vspace{-.5cm}
\caption{\centering Reward in maze 1 (left): 10 for reaching target, -5 for crossing the second and third bridge, -1 for crossing the fourth bridge, and -0.1 for any other move. Reward in maze 2 (middle): 100 for reaching target, -50 for crossing the first and third bridge, -10 for crossing the fourth bridge, and -1 for any other move. Reward in maze 3 (right): 1000 for reaching target, -500 for crossing the second and the third bridge, -100 for crossing the fourth bridge, and -10 for any other move. Green dotted lines indicate optimal paths for local tasks. The yellow dotted line indicates a sub-optimal but acceptable policy for each local task, which is also the globally optimal policy of the constrained multi-task problem under $\ell_1=5$, $\ell_2=50$, $\ell_3=500$.}
\label{fig:maze}
\end{figure}

In this section, we use a small-scale GridWorld experiment to show how the constrained formulation allows us to control the multi-task policy in a fine-grained manner. We consider a three-task RL problem.
Associated with each task is a $10 \times 10$ maze in which an agent aims to navigate to a goal position by crossing bridges that charge different ``prices'' (i.e. incur different negative rewards). Shown in Figure~\ref{fig:maze}, the starting position for all tasks is the top left corner and the target is the top right corner, both marked in blue. There is a positive reward for reaching the target, which varies in magnitude across tasks.

In the first task, the agent receives the least negative reward for using the first bridge, which is -0.1. The reward of the fourth bridge, -1, is more negative, and the reward of the second and third bridges are the most negative. The reward of any other move in the maze is -0.1. If we are only interested in solving task 1, the optimal policies (one of which is drawn as the green path in the figure) obviously want to use the first bridge to reach the target as soon as possible.\looseness=-1

The second task can be regarded as a reward-magnified version of task 1. With the rewards marked in Figure~\ref{fig:maze}, the best bridge becomes the second one, which charges much lower prices than the other bridges. Making any other move not labeled incurs a reward of $-1$. One of the optimal policies for this task is drawn as the green dotted path, which uses the second bridge.

The rewards for all moves in the third task are further scaled up. The bonus of reaching the target increases to 1000, while the costs of crossing the bridges becomes 500, 100, and 10. The reward of making any other move also changes to -10. Since the magnitude of the rewards in this task is dominant over those in the other two tasks, it is easy to verify that the optimal policies of the third task coincide with the globally optimal policy for the unconstrained multi-task RL objective in Eq.~\eqref{eq:unconstrained_MTRL}, one of which is again shown in green in the figure. This optimal path takes the third bridge, which ensures that the greatest reward is collected in the third tasks, at the cost of completely unsatisfactory performance in the other two tasks.

\begin{figure}[ht]
\centering
\includegraphics[width=\linewidth]{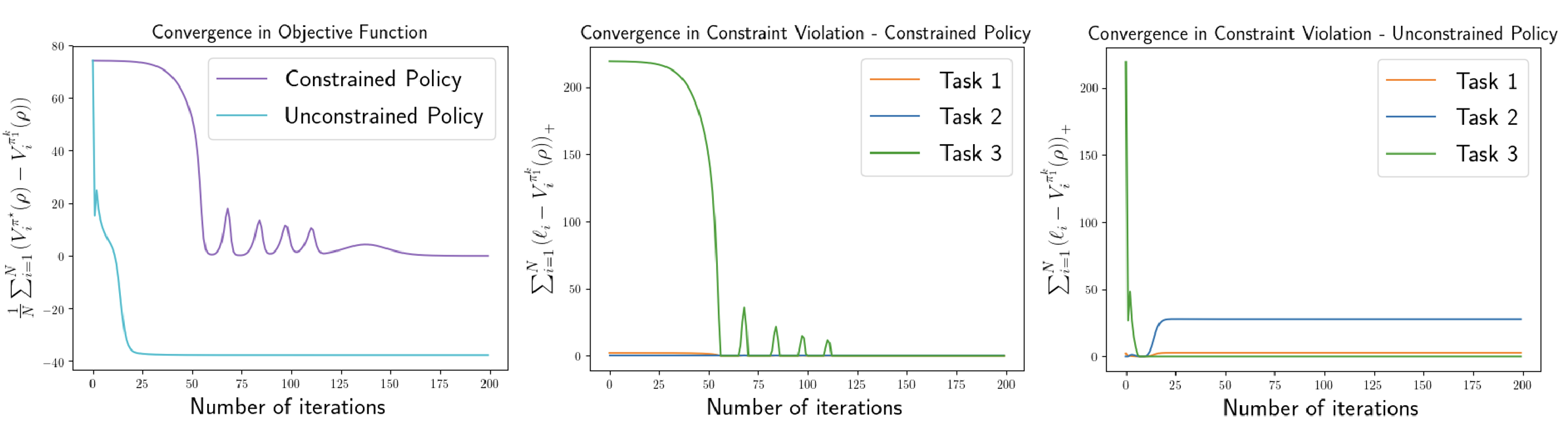}
\vspace{-.5cm}
\caption{\centering Left -- convergence of Algorithm~\ref{Alg:MT-PDNPG} in objective function with $\ell_1=5$, $\ell_2=50$, $\ell_3=500$ (constrained policy), and with $\ell_1=\ell_2=\ell_3=-\infty$ (unconstrained policy). Middle -- convergence of Algorithm~\ref{Alg:MT-PDNPG} in constraint violation with $\ell_1=5$, $\ell_2=50$, $\ell_3=500$. Right -- convergence of Algorithm~\ref{Alg:MT-PDNPG} in constraint violation with $\ell_1=\ell_2=\ell_3=-\infty$.}
\vspace{-.2cm}
\label{fig:simulations}
\end{figure}

Comparing with solving Eq.~\eqref{eq:unconstrained_MTRL}, our formulation in Eq.~\eqref{eq:objective} allows us to trade off the policy performance in the three tasks by properly specifying the performance lower bounds. In particular, we set $\ell_1=5$,$\ell_2=50$, and $\ell_3=500$. Calculations show that the optimal policies under this set of constraints will switch to take the fourth bridge, which means a slight but acceptable compromise of the policy performance in all tasks. Numerically, we apply our proposed primal-dual natural gradient algorithm\footnote{The small dimension and known transition kernel of the environments allow us to derive the exact gradient in this case.}, and verify that the local policy at every agent indeed converges to the optimal constrained policy in both objective function and constraint violation (see the first and second plots in Figure~\ref{fig:simulations}). In contrast, we also apply the same algorithm with $\ell_1=\ell_2=\ell_3=-\infty$, which essentially means that we run a decentralized natural policy gradient algorithm to solve the unconstrained multi-task problem in Eq.~\eqref{eq:unconstrained_MTRL}. The policy iterates of this unconstrained problem obviously achieves better objective function value, but fails to satisfy the constraints.

\section{Conclusion}

This paper considers a constrained multi-task RL objective where the goal is to find a policy with the maximum average performance across all tasks and guaranteed minimum performance in each environment. We show that three policy-gradient-based algorithms provably solve the problem, under different policy parameterization and information oracle. In the tabular setting with the access to the exact gradients, we study a primal-dual natural-gradient-descent-ascent algorithm that drives each agent's local policy to the globally optimal solution with finite-time complexity $\Ocal(K^{-1/2})$. When the exact gradient information is not available, we take a completely sampled-based actor-critic approach and show that it converges with rate $\Ocal(K^{-1/6})$. Finally, we extend the results to the linear function approximation setting and establish the same sample complexity as in the tabular case. To our best knowledge, a finite-time and finite-sample analysis for the on-policy natural actor-critic algorithm in the linear function approximation setting (even for standard non-constrained MDPs) has been missing from the existing literature, and our work fills in this gap of knowledge.

\section*{Disclaimer}
This paper was prepared for informational purposes in part by
the Artificial Intelligence Research group of JP Morgan Chase \& Co and its affiliates (``JP Morgan''),
and is not a product of the Research Department of JP Morgan.
JP Morgan makes no representation and warranty whatsoever and disclaims all liability,
for the completeness, accuracy or reliability of the information contained herein.
This document is not intended as investment research or investment advice, or a recommendation,
offer or solicitation for the purchase or sale of any security, financial instrument, financial product or service,
or to be used in any way for evaluating the merits of participating in any transaction,
and shall not constitute a solicitation under any jurisdiction or to any person,
if such solicitation under such jurisdiction or to such person would be unlawful.

\bibliography{references}

\begin{thebibliography}{57}
\providecommand{\natexlab}[1]{#1}
\providecommand{\url}[1]{\texttt{#1}}
\expandafter\ifx\csname urlstyle\endcsname\relax
  \providecommand{\doi}[1]{doi: #1}\else
  \providecommand{\doi}{doi: \begingroup \urlstyle{rm}\Url}\fi

\bibitem[Agarwal et~al.(2020)Agarwal, Kakade, Lee, and Mahajan]{agarwal2020optimality}
Alekh Agarwal, Sham~M Kakade, Jason~D Lee, and Gaurav Mahajan.
\newblock Optimality and approximation with policy gradient methods in markov decision processes.
\newblock In \emph{Conference on Learning Theory}, pp.\  64--66. PMLR, 2020.

\bibitem[Agarwal et~al.(2022)Agarwal, Bai, and Aggarwal]{agarwal2022regret}
Mridul Agarwal, Qinbo Bai, and Vaneet Aggarwal.
\newblock Regret guarantees for model-based reinforcement learning with long-term average constraints.
\newblock In \emph{Uncertainty in Artificial Intelligence}, pp.\  22--31. PMLR, 2022.

\bibitem[Altman(1999)]{altman1999constrained}
Eitan Altman.
\newblock \emph{Constrained {M}arkov decision processes}, volume~7.
\newblock Chapman and Hall/CRC Press, 1999.

\bibitem[Bai et~al.(2022)Bai, Bedi, Agarwal, Koppel, and Aggarwal]{bai2022achieving}
Qinbo Bai, Amrit~Singh Bedi, Mridul Agarwal, Alec Koppel, and Vaneet Aggarwal.
\newblock Achieving zero constraint violation for constrained reinforcement learning via primal-dual approach.
\newblock In \emph{Proceedings of the AAAI Conference on Artificial Intelligence}, volume~36, pp.\  3682--3689, 2022.

\bibitem[Barakat et~al.(2022)Barakat, Bianchi, and Lehmann]{barakat2022analysis}
Anas Barakat, Pascal Bianchi, and Julien Lehmann.
\newblock Analysis of a target-based actor-critic algorithm with linear function approximation.
\newblock In \emph{International Conference on Artificial Intelligence and Statistics}, pp.\  991--1040. PMLR, 2022.

\bibitem[Borkar(2005)]{borkar2005actor}
Vivek~S Borkar.
\newblock An actor-critic algorithm for constrained markov decision processes.
\newblock \emph{Systems \& control letters}, 54\penalty0 (3):\penalty0 207--213, 2005.

\bibitem[Bullins et~al.(2021)Bullins, Patel, Shamir, Srebro, and Woodworth]{bullins2021stochastic}
Brian Bullins, Kshitij Patel, Ohad Shamir, Nathan Srebro, and Blake~E Woodworth.
\newblock A stochastic newton algorithm for distributed convex optimization.
\newblock \emph{Advances in Neural Information Processing Systems}, 34:\penalty0 26818--26830, 2021.

\bibitem[Chang et~al.(2014)Chang, Nedi{\'c}, and Scaglione]{chang2014distributed}
Tsung-Hui Chang, Angelia Nedi{\'c}, and Anna Scaglione.
\newblock Distributed constrained optimization by consensus-based primal-dual perturbation method.
\newblock \emph{IEEE Transactions on Automatic Control}, 59\penalty0 (6):\penalty0 1524--1538, 2014.

\bibitem[Chen et~al.(2022{\natexlab{a}})Chen, Feng, Gao, and Wei]{chen2022decentralized}
Jinchi Chen, Jie Feng, Weiguo Gao, and Ke~Wei.
\newblock Decentralized natural policy gradient with variance reduction for collaborative multi-agent reinforcement learning.
\newblock \emph{arXiv preprint arXiv:2209.02179}, 2022{\natexlab{a}}.

\bibitem[Chen et~al.(2022{\natexlab{b}})Chen, Khodadadian, and Maguluri]{chen2022finite}
Zaiwei Chen, Sajad Khodadadian, and Siva~Theja Maguluri.
\newblock Finite-sample analysis of off-policy natural actor--critic with linear function approximation.
\newblock \emph{IEEE Control Systems Letters}, 6:\penalty0 2611--2616, 2022{\natexlab{b}}.

\bibitem[Chow et~al.(2017)Chow, Ghavamzadeh, Janson, and Pavone]{chow2017risk}
Yinlam Chow, Mohammad Ghavamzadeh, Lucas Janson, and Marco Pavone.
\newblock Risk-constrained reinforcement learning with percentile risk criteria.
\newblock \emph{The Journal of Machine Learning Research}, 18\penalty0 (1):\penalty0 6070--6120, 2017.

\bibitem[Chow et~al.(2018)Chow, Nachum, Duenez-Guzman, and Ghavamzadeh]{chow2018lyapunov}
Yinlam Chow, Ofir Nachum, Edgar Duenez-Guzman, and Mohammad Ghavamzadeh.
\newblock A lyapunov-based approach to safe reinforcement learning.
\newblock \emph{arXiv preprint arXiv:1805.07708}, 2018.

\bibitem[D'Eramo et~al.()D'Eramo, Tateo, Bonarini, Restelli, and Peters]{dsharing}
Carlo D'Eramo, Davide Tateo, Andrea Bonarini, Marcello Restelli, and Jan Peters.
\newblock Sharing knowledge in multi-task deep reinforcement learning.
\newblock In \emph{International Conference on Learning Representations}.

\bibitem[Ding et~al.(2020)Ding, Zhang, Basar, and Jovanovic]{ding2020natural}
Dongsheng Ding, Kaiqing Zhang, Tamer Basar, and Mihailo Jovanovic.
\newblock Natural policy gradient primal-dual method for constrained markov decision processes.
\newblock \emph{Advances in Neural Information Processing Systems}, 33:\penalty0 8378--8390, 2020.

\bibitem[Ding et~al.(2021)Ding, Wei, Yang, Wang, and Jovanovic]{ding2021provably}
Dongsheng Ding, Xiaohan Wei, Zhuoran Yang, Zhaoran Wang, and Mihailo Jovanovic.
\newblock Provably efficient safe exploration via primal-dual policy optimization.
\newblock In \emph{International Conference on Artificial Intelligence and Statistics}, pp.\  3304--3312. PMLR, 2021.

\bibitem[Ding et~al.(2022)Ding, Zhang, Ba{\c{s}}ar, and Jovanovi{\'c}]{ding2022convergence}
Dongsheng Ding, Kaiqing Zhang, Tamer Ba{\c{s}}ar, and Mihailo~R Jovanovi{\'c}.
\newblock Convergence and optimality of policy gradient primal-dual method for constrained markov decision processes.
\newblock In \emph{2022 American Control Conference (ACC)}, pp.\  2851--2856. IEEE, 2022.

\bibitem[Doan et~al.(2019)Doan, Maguluri, and Romberg]{doan2019finite}
Thinh Doan, Siva Maguluri, and Justin Romberg.
\newblock Finite-time analysis of distributed td (0) with linear function approximation on multi-agent reinforcement learning.
\newblock In \emph{International Conference on Machine Learning}, pp.\  1626--1635. PMLR, 2019.

\bibitem[Finn et~al.(2017)Finn, Abbeel, and Levine]{finn2017model}
Chelsea Finn, Pieter Abbeel, and Sergey Levine.
\newblock Model-agnostic meta-learning for fast adaptation of deep networks.
\newblock In \emph{International conference on machine learning}, pp.\  1126--1135. PMLR, 2017.

\bibitem[Guo et~al.(2022)Guo, Wu, and Lee]{guo2022learning}
Yijie Guo, Qiucheng Wu, and Honglak Lee.
\newblock Learning action translator for meta reinforcement learning on sparse-reward tasks.
\newblock In \emph{Proceedings of the AAAI Conference on Artificial Intelligence}, volume~36, pp.\  6792--6800, 2022.

\bibitem[Gupta et~al.(2017)Gupta, Devin, Liu, Abbeel, and Levine]{gupta2017learning}
Abhishek Gupta, Coline Devin, YuXuan Liu, Pieter Abbeel, and Sergey Levine.
\newblock Learning invariant feature spaces to transfer skills with reinforcement learning.
\newblock \emph{arXiv preprint arXiv:1703.02949}, 2017.

\bibitem[HasanzadeZonuzy et~al.(2021)HasanzadeZonuzy, Bura, Kalathil, and Shakkottai]{hasanzadezonuzy2021learning}
Aria HasanzadeZonuzy, Archana Bura, Dileep Kalathil, and Srinivas Shakkottai.
\newblock Learning with safety constraints: Sample complexity of reinforcement learning for constrained mdps.
\newblock In \emph{Proceedings of the AAAI Conference on Artificial Intelligence}, volume~35, pp.\  7667--7674, 2021.

\bibitem[Hayes et~al.(2022)Hayes, R{\u{a}}dulescu, Bargiacchi, K{\"a}llstr{\"o}m, Macfarlane, Reymond, Verstraeten, Zintgraf, Dazeley, Heintz, et~al.]{hayes2022practical}
Conor~F Hayes, Roxana R{\u{a}}dulescu, Eugenio Bargiacchi, Johan K{\"a}llstr{\"o}m, Matthew Macfarlane, Mathieu Reymond, Timothy Verstraeten, Luisa~M Zintgraf, Richard Dazeley, Fredrik Heintz, et~al.
\newblock A practical guide to multi-objective reinforcement learning and planning.
\newblock \emph{Autonomous Agents and Multi-Agent Systems}, 36\penalty0 (1):\penalty0 26, 2022.

\bibitem[Hessel et~al.(2019)Hessel, Soyer, Espeholt, Czarnecki, Schmitt, and van Hasselt]{hessel2019multi}
Matteo Hessel, Hubert Soyer, Lasse Espeholt, Wojciech Czarnecki, Simon Schmitt, and Hado van Hasselt.
\newblock Multi-task deep reinforcement learning with popart.
\newblock In \emph{Proceedings of the AAAI Conference on Artificial Intelligence}, volume~33, pp.\  3796--3803, 2019.

\bibitem[Hong et~al.(2021)Hong, Yoon, and Kim]{hong2021structure}
Sunghoon Hong, Deunsol Yoon, and Kee-Eung Kim.
\newblock Structure-aware transformer policy for inhomogeneous multi-task reinforcement learning.
\newblock In \emph{International Conference on Learning Representations}, 2021.

\bibitem[Islamov et~al.(2021)Islamov, Qian, and Richt{\'a}rik]{islamov2021distributed}
Rustem Islamov, Xun Qian, and Peter Richt{\'a}rik.
\newblock Distributed second order methods with fast rates and compressed communication.
\newblock In \emph{International conference on machine learning}, pp.\  4617--4628. PMLR, 2021.

\bibitem[Jiang et~al.(2022)Jiang, Lee, Tan, Tan, Balu, Lee, Hegde, and Sarkar]{jiang2022mdpgt}
Zhanhong Jiang, Xian~Yeow Lee, Sin~Yong Tan, Kai~Liang Tan, Aditya Balu, Young~M Lee, Chinmay Hegde, and Soumik Sarkar.
\newblock Mdpgt: momentum-based decentralized policy gradient tracking.
\newblock In \emph{Proceedings of the AAAI Conference on Artificial Intelligence}, volume~36, pp.\  9377--9385, 2022.

\bibitem[Ju et~al.(2022)Ju, Kotsalis, and Lan]{ju2022model}
Caleb Ju, Georgios Kotsalis, and Guanghui Lan.
\newblock A model-free first-order method for linear quadratic regulator with $\tilde {O}(1/\varepsilon)$ sampling complexity.
\newblock \emph{arXiv preprint arXiv:2212.00084}, 2022.

\bibitem[Junru et~al.(2022)Junru, Qiong, Muhua, Zhihang, Ruijuan, and Qingtao]{junru2022decentralized}
Shi Junru, Wang Qiong, Liu Muhua, Ji~Zhihang, Zheng Ruijuan, and Wu~Qingtao.
\newblock Decentralized multi-task reinforcement learning policy gradient method with momentum over networks.
\newblock \emph{Applied Intelligence}, pp.\  1--15, 2022.

\bibitem[Kalashnikov et~al.(2021)Kalashnikov, Varley, Chebotar, Swanson, Jonschkowski, Finn, Levine, and Hausman]{kalashnikov2021mt}
Dmitry Kalashnikov, Jacob Varley, Yevgen Chebotar, Benjamin Swanson, Rico Jonschkowski, Chelsea Finn, Sergey Levine, and Karol Hausman.
\newblock Mt-opt: Continuous multi-task robotic reinforcement learning at scale.
\newblock \emph{arXiv preprint arXiv:2104.08212}, 2021.

\bibitem[Khodadadian et~al.(2022)Khodadadian, Doan, Romberg, and Maguluri]{khodadadian2022finite}
Sajad Khodadadian, Thinh~T Doan, Justin Romberg, and Siva~Theja Maguluri.
\newblock Finite sample analysis of two-time-scale natural actor-critic algorithm.
\newblock \emph{IEEE Transactions on Automatic Control}, 2022.

\bibitem[Lei et~al.(2016)Lei, Chen, and Fang]{lei2016primal}
Jinlong Lei, Han-Fu Chen, and Hai-Tao Fang.
\newblock Primal--dual algorithm for distributed constrained optimization.
\newblock \emph{Systems \& Control Letters}, 96:\penalty0 110--117, 2016.

\bibitem[Liu et~al.(2021{\natexlab{a}})Liu, Zhou, Kalathil, Kumar, and Tian]{liu2021learning}
Tao Liu, Ruida Zhou, Dileep Kalathil, Panganamala Kumar, and Chao Tian.
\newblock Learning policies with zero or bounded constraint violation for constrained mdps.
\newblock \emph{Advances in Neural Information Processing Systems}, 34:\penalty0 17183--17193, 2021{\natexlab{a}}.

\bibitem[Liu et~al.(2021{\natexlab{b}})Liu, Zhou, Kalathil, Kumar, and Tian]{liu2021policy}
Tao Liu, Ruida Zhou, Dileep Kalathil, PR~Kumar, and Chao Tian.
\newblock Policy optimization for constrained mdps with provable fast global convergence.
\newblock \emph{arXiv preprint arXiv:2111.00552}, 2021{\natexlab{b}}.

\bibitem[Nedic(2020)]{nedic2020distributed}
Angelia Nedic.
\newblock Distributed gradient methods for convex machine learning problems in networks: Distributed optimization.
\newblock \emph{IEEE Signal Processing Magazine}, 37\penalty0 (3):\penalty0 92--101, 2020.

\bibitem[Nedic \& Ozdaglar(2009)Nedic and Ozdaglar]{nedic2009distributed}
Angelia Nedic and Asuman Ozdaglar.
\newblock Distributed subgradient methods for multi-agent optimization.
\newblock \emph{IEEE Transactions on Automatic Control}, 54\penalty0 (1):\penalty0 48--61, 2009.

\bibitem[Olshevsky(2015)]{olshevsky2015linear}
Alex Olshevsky.
\newblock Linear time average consensus on fixed graphs.
\newblock \emph{IFAC-PapersOnLine}, 48\penalty0 (22):\penalty0 94--99, 2015.

\bibitem[Paternain et~al.(2019)Paternain, Chamon, Calvo-Fullana, and Ribeiro]{paternain2019constrained}
Santiago Paternain, Luiz~FO Chamon, Miguel Calvo-Fullana, and Alejandro Ribeiro.
\newblock Constrained reinforcement learning has zero duality gap.
\newblock \emph{arXiv preprint arXiv:1910.13393}, 2019.

\bibitem[Prashanth \& Ghavamzadeh(2016)Prashanth and Ghavamzadeh]{prashanth2016variance}
LA~Prashanth and Mohammad Ghavamzadeh.
\newblock Variance-constrained actor-critic algorithms for discounted and average reward mdps.
\newblock \emph{Machine Learning}, 105\penalty0 (3):\penalty0 367--417, 2016.

\bibitem[Raghu et~al.()Raghu, Raghu, Bengio, and Vinyals]{raghurapid}
Aniruddh Raghu, Maithra Raghu, Samy Bengio, and Oriol Vinyals.
\newblock Rapid learning or feature reuse? towards understanding the effectiveness of maml.
\newblock In \emph{International Conference on Learning Representations}.

\bibitem[Rusu et~al.(2015)Rusu, Colmenarejo, Gulcehre, Desjardins, Kirkpatrick, Pascanu, Mnih, Kavukcuoglu, and Hadsell]{rusu2015policy}
Andrei~A Rusu, Sergio~Gomez Colmenarejo, Caglar Gulcehre, Guillaume Desjardins, James Kirkpatrick, Razvan Pascanu, Volodymyr Mnih, Koray Kavukcuoglu, and Raia Hadsell.
\newblock Policy distillation.
\newblock \emph{arXiv preprint arXiv:1511.06295}, 2015.

\bibitem[Tessler et~al.(2018)Tessler, Mankowitz, and Mannor]{tessler2018reward}
Chen Tessler, Daniel~J Mankowitz, and Shie Mannor.
\newblock Reward constrained policy optimization.
\newblock \emph{arXiv preprint arXiv:1805.11074}, 2018.

\bibitem[Traor{\'e} et~al.(2019)Traor{\'e}, Caselles-Dupr{\'e}, Lesort, Sun, Cai, D{\'\i}az-Rodr{\'\i}guez, and Filliat]{traore2019discorl}
Ren{\'e} Traor{\'e}, Hugo Caselles-Dupr{\'e}, Timoth{\'e}e Lesort, Te~Sun, Guanghang Cai, Natalia D{\'\i}az-Rodr{\'\i}guez, and David Filliat.
\newblock Discorl: Continual reinforcement learning via policy distillation.
\newblock \emph{arXiv preprint arXiv:1907.05855}, 2019.

\bibitem[Tsitsiklis \& Vanroy(1997)Tsitsiklis and Vanroy]{tsitsiklis1997analysis}
JN~Tsitsiklis and B~Vanroy.
\newblock An analysis of temporal-difference learning with function approximation.
\newblock \emph{IEEE Transactions on Automatic Control}, 42\penalty0 (5):\penalty0 674--690, 1997.

\bibitem[Wadhwania et~al.(2019)Wadhwania, Kim, Omidshafiei, and How]{wadhwania2019policy}
Samir Wadhwania, Dong-Ki Kim, Shayegan Omidshafiei, and Jonathan~P How.
\newblock Policy distillation and value matching in multiagent reinforcement learning.
\newblock In \emph{2019 IEEE/RSJ International Conference on Intelligent Robots and Systems (IROS)}, pp.\  8193--8200. IEEE, 2019.

\bibitem[Wu et~al.(2020)Wu, Zhang, Xu, and Gu]{wu2020finite}
Yue~Frank Wu, Weitong Zhang, Pan Xu, and Quanquan Gu.
\newblock A finite-time analysis of two time-scale actor-critic methods.
\newblock \emph{Advances in Neural Information Processing Systems}, 33:\penalty0 17617--17628, 2020.

\bibitem[Xu et~al.(2020)Xu, Wang, and Liang]{xu2020non}
Tengyu Xu, Zhe Wang, and Yingbin Liang.
\newblock Non-asymptotic convergence analysis of two time-scale (natural) actor-critic algorithms.
\newblock \emph{arXiv preprint arXiv:2005.03557}, 2020.

\bibitem[Yang et~al.(2020)Yang, Xu, Wu, and Wang]{yang2020multi}
Ruihan Yang, Huazhe Xu, Yi~Wu, and Xiaolong Wang.
\newblock Multi-task reinforcement learning with soft modularization.
\newblock \emph{Advances in Neural Information Processing Systems}, 33:\penalty0 4767--4777, 2020.

\bibitem[Yang et~al.(2019)Yang, Chen, Hong, and Wang]{yang2019provably}
Zhuoran Yang, Yongxin Chen, Mingyi Hong, and Zhaoran Wang.
\newblock Provably global convergence of actor-critic: A case for linear quadratic regulator with ergodic cost.
\newblock \emph{Advances in neural information processing systems}, 32, 2019.

\bibitem[Ying et~al.(2022)Ying, Ding, and Lavaei]{ying2022dual}
Donghao Ying, Yuhao Ding, and Javad Lavaei.
\newblock A dual approach to constrained markov decision processes with entropy regularization.
\newblock In \emph{International Conference on Artificial Intelligence and Statistics}, pp.\  1887--1909. PMLR, 2022.

\bibitem[Yu et~al.(2019)Yu, Yang, Kolar, and Wang]{yu2019convergent}
Ming Yu, Zhuoran Yang, Mladen Kolar, and Zhaoran Wang.
\newblock Convergent policy optimization for safe reinforcement learning.
\newblock \emph{Advances in Neural Information Processing Systems}, 32:\penalty0 3127--3139, 2019.

\bibitem[Yuan et~al.(2016)Yuan, Ling, and Yin]{yuan2016convergence}
Kun Yuan, Qing Ling, and Wotao Yin.
\newblock On the convergence of decentralized gradient descent.
\newblock \emph{SIAM Journal on Optimization}, 26\penalty0 (3):\penalty0 1835--1854, 2016.

\bibitem[Zeng et~al.(2021)Zeng, Anwar, Doan, Raychowdhury, and Romberg]{zeng2021decentralized}
Sihan Zeng, Malik~Aqeel Anwar, Thinh~T Doan, Arijit Raychowdhury, and Justin Romberg.
\newblock A decentralized policy gradient approach to multi-task reinforcement learning.
\newblock In \emph{Uncertainty in Artificial Intelligence}, pp.\  1002--1012. PMLR, 2021.

\bibitem[Zeng et~al.(2022)Zeng, Doan, and Romberg]{zeng2022finite}
Sihan Zeng, Thinh~T Doan, and Justin Romberg.
\newblock Finite-time complexity of online primal-dual natural actor-critic algorithm for constrained markov decision processes.
\newblock In \emph{2022 IEEE 61st Conference on Decision and Control (CDC)}, pp.\  4028--4033. IEEE, 2022.

\bibitem[Zeng et~al.(2023)Zeng, Doan, and Romberg]{zeng2022finite2}
Sihan Zeng, Thinh~T Doan, and Justin Romberg.
\newblock Finite-time convergence rates of decentralized stochastic approximation with applications in multi-agent and multi-task learning.
\newblock \emph{IEEE Transactions on Automatic Control}, 68:\penalty0 2758--2773, 2023.

\bibitem[Zeng et~al.(2024)Zeng, Doan, and Romberg]{zeng2021two}
Sihan Zeng, Thinh~T Doan, and Justin Romberg.
\newblock A two-time-scale stochastic optimization framework with applications in control and reinforcement learning.
\newblock \emph{SIAM Journal on Optimization}, 34\penalty0 (1):\penalty0 946--976, 2024.

\bibitem[Zhang et~al.(2018)Zhang, Yang, Liu, Zhang, and Basar]{zhang2018fully}
Kaiqing Zhang, Zhuoran Yang, Han Liu, Tong Zhang, and Tamer Basar.
\newblock Fully decentralized multi-agent reinforcement learning with networked agents.
\newblock In \emph{International Conference on Machine Learning}, pp.\  5872--5881. PMLR, 2018.

\bibitem[Zheng \& Ratliff(2020)Zheng and Ratliff]{zheng2020constrained}
Liyuan Zheng and Lillian Ratliff.
\newblock Constrained upper confidence reinforcement learning.
\newblock In \emph{Learning for Dynamics and Control}, pp.\  620--629. PMLR, 2020.

\end{thebibliography}
\bibliographystyle{tmlr}

\clearpage
\appendix

\section{Additional Notations}\label{sec:notations}
We introduce some additional shorthand notations frequently used in the appendix. First, we define $\widehat{V}_{i}^k\in\mathbb{R}^{|\Scal|}$ and $\widehat{A}_{i}^k\in\mathbb{R}^{|\Scal|\times|\Acal|}$ as the estimated value function and advantage function based on $\widehat{Q}_i^k$ such that
\begin{align*}
    \widehat{V}_{i}^k(s)&\triangleq\sum_{a'\in\Acal}\pi_i^k(a\mid s)\widehat{Q}_{i}^k(s,a),\,\,\forall s\in\Scal,\\
    \widehat{A}_{i}^k(s,a) &\triangleq \widehat{Q}_{i}^k(s,a)\hspace{-1pt}-\hspace{-1pt}\widehat{V}_{i}^k(s),\,\,\forall s\in\Scal,a\in\Acal.
\end{align*}

We use the subscript $g$ in the value/Q/advantage function to denote the vector formed by stacking the value/Q/advantage functions across agents $i=1,\dots,N$. Specifically, given a policy $\pi\in\Delta_{\Acal}^{\Scal}$, we write 
\begin{align*}
    V_{g}^{\pi}(s) &\triangleq [V_{1}^{\pi}(s),...,V_{N}^{\pi}(s)]^{\top}\in\mathbb{R}^{N}\\
    Q_{g}^{\pi}(s,a) &\triangleq [Q_{1}^{\pi}(s,a),...,Q_{N}^{\pi}(s,a)]^{\top}\in\mathbb{R}^{N}\\
    A_{g}^{\pi}(s,a) &\triangleq [A_{1}^{\pi}(s,a),...,A_{N}^{\pi}(s,a)]^{\top}\in\mathbb{R}^{N}
\end{align*}

\section{Proof of Theorems}

\subsection{Proof of Theorem \ref{thm:main_pg}}\label{sec:thm:main_pg}

The proof of Theorem \ref{thm:main_pg} relies on the proposition below, which is an intermediate result that characterizes the convergence of primal and dual variables under a general update rule. This proposition will also be used in the proofs of Theorem~\ref{thm:main_ac} and \ref{thm:LFA}.

\begin{prop}\label{prop:actor_conv}
Consider any decentralized algorithm where each agent in iteration $k$ maintains dual variables $\lambda_i^k,\nu_i^k\in\mathbb{R}$, critic variables $Q_i^k\in\mathbb{R}^{|\Scal||\Acal|}$, and a primal variable $\theta_i^k\in\mathbb{R}^d$ (for any positive integer $d$) which parameterizes a policy $\pi_i^{k} = g_i(\theta_i^{k})\in\Delta_{\Acal}^{\Scal}$ through some mapping $g_i$. Suppose the dual variable updates satisfy
\begin{align}
    \notag\\
    &\lambda_{i}^{k+1} =\Pi_{[0,B_{\lambda}]}\Big(\lambda_{i}^{k}-\eta\big(\sum_{s,a}\rho(s)\pi_i^k(a\mid s) Q_i^{k}(s,a)-\ell_i\big)\Big),\notag\\
    &\nu_{i}^{k+1} =\Pi_{[0,B_{\lambda}]}\Big(\lambda_{i}^{k}+\eta\big(\sum_{s,a}\rho(s)\pi_i^k(a\mid s) Q_i^{k}(s,a)-u_i\big)\Big).\label{prop:actor_conv:updates}
\end{align}
Also suppose that there exists a constant $B>0$ such that $\|Q_i^{k}\|_{\infty}\leq B$ for all $i$ and $k$, and that there exists a mapping $f:\mathbb{R}^{N\times d}\rightarrow\Delta_{\Acal}^{\Scal}$ such that $\widebar{\pi}^{k+1}\triangleq f(\{\theta_i^{k+1}\}_i)$ observes the recursive rule
\begin{align}
\widebar{\pi}^{k+1}(a\mid s)\propto\widebar{\pi}^{k}(a\mid s)\exp(\frac{\alpha}{N}\sum_{i=1}^{N}(\frac{1}{N}+\lambda_i^k-\nu_i^k)Q_i^{k}(s,a)),\quad\forall k.\label{prop:actor_conv:updates_pibar}
\end{align}
Then, the parameters $\{\pi_i^k\}_{i,k}$ after $K$ update iterations satisfy
\begin{align*}
    &\max\{\frac{1}{K}\sum_{k=0}^{K-1}\left(V_{0}^{\pi^{\star}}(\rho)-V_{0}^{\widebar{\pi}^{k}}(\rho)\right),\frac{1}{K}\sum_{k=0}^{K-1}\sum_{i=1}^{N}\left(\big[\ell_i-V_{i}^{\widebar{\pi}^{k}}(\rho)\big]_{+}+\big[V_{i}^{\widebar{\pi}^{k}}(\rho)-u_i\big]_{+}\right)\}\\
    &\leq \Ocal\left(\frac{N}{K\alpha}+N\eta+\frac{N}{K\eta}+\frac{1}{K}\sum_{k=0}^{K-1}\sum_{i=1}^{N}\|\widebar{\pi}^k-\pi_i^k\|+\frac{1}{K}\sum_{k=0}^{K-1}\sum_{i=1}^{N} \|Q_{i}^{\pi_i^k}-Q_{i}^{k}\|\right).
\end{align*}
\end{prop}

We also introduce the following lemma which establishes the Lipschitz continuity of the value function and Q function.
\begin{lem}[Lemma 8 of \cite{khodadadian2022finite}]\label{lem:misc}

For any policy $\pi_1,\pi_2$ and $i=1,\dots,N$
\begin{gather*}
    \|Q_i^{\pi_1}-Q_i^{\pi_2}\|\leq\frac{|\Scal|| \Acal|}{(1-\gamma)^2}\|\pi_1-\pi_2\|,\quad \|V_i^{\pi_1}-V_i^{\pi_2}\|\leq\frac{|\Scal|| \Acal|}{(1-\gamma)^2}\|\pi_1-\pi_2\|,\\
    |Q_i^{\pi_1}(s,a)-Q_i^{\pi_2}(s,a)|\leq\frac{\sqrt{|\Scal|| \Acal|}}{(1-\gamma)^2}\|\pi_1-\pi_2\|,\quad |V_i^{\pi_1}(s)-V_i^{\pi_2}(s)|\leq\frac{\sqrt{|\Scal|| \Acal|}}{(1-\gamma)^2}\|\pi_1-\pi_2\|.
\end{gather*}
\end{lem}

It is easy to verify that Algorithm~\ref{Alg:MT-PDNPG} observes the update rule in \ref{prop:actor_conv:updates} and \ref{prop:actor_conv:updates_pibar}, with $Q_i^k=Q_i^{\pi_i^k}$ and $\widebar{\pi}^{k+1}=f(\{\theta_i^{k+1}\}_i)$ defined as
\begin{align*}
    \widebar{\theta}^{k+1}=\frac{1}{N}\sum_{i=1}^{N}\theta_i^{k+1},\quad\widebar{\pi}^{k+1}=\frac{\exp \left(\widebar{\theta}^{k+1}(s,a)\right)}{\sum_{a' \in \Acal} \exp\left(\widebar{\theta}^{k+1}(s,a')\right)}
\end{align*}

Due to the bounded update, we can control the distance $\|\widebar{\pi}^k-\pi_i^k\|$ by the step size $\alpha$.

\begin{lem}\label{lem:consensus_error_pg}
The policy iterates $\{\pi_i^k\}$ generated by Algorithm~\ref{Alg:MT-PDNPG} satisfy
\begin{align*}
    \|\widebar{\pi}^k-\pi_i^k\|\leq\Ocal\left(\frac{\sqrt{N}\alpha}{1-\sigma_2(W)}\right),\quad\text{for all $k=0,\dots,K-1$ and $i=1,\dots,N$}.
\end{align*}
\end{lem}

As a result of the proposition and lemma above,
\begin{align*}
    &\max\{\frac{1}{K}\sum_{k=0}^{K-1}\left(V_{0}^{\pi^{\star}}(\rho)-V_{0}^{\widebar{\pi}^{k}}(\rho)\right),\frac{1}{K}\sum_{k=0}^{K-1}\sum_{i=1}^{N}\left(\big[\ell_i-V_{i}^{\widebar{\pi}^{k}}(\rho)\big]_{+}+\big[V_{i}^{\widebar{\pi}^{k}}(\rho)-u_i\big]_{+}\right)\}\\
    &\leq \Ocal\left(\frac{N}{K\alpha}+N\eta+\frac{N}{K\eta}+\frac{1}{K}\sum_{k=0}^{K-1}\sum_{i=1}^{N}\|\widebar{\pi}^k-\pi_i^k\|\right)\notag\\
    &=\Ocal\left(\frac{N}{K\alpha}+N\eta+\frac{N}{K\eta}+\frac{N^{3/2}\alpha}{1-\sigma_2(W)}\right).
\end{align*}

By the bound on consensus error in Lemma~\ref{lem:consensus_error_pg} and the Lipschitz continuity of the value function in Lemma~\ref{lem:misc}, this implies for any agent $j=1,\cdots,N$
\begin{align*}
    &\max\{\frac{1}{K}\sum_{k=0}^{K-1}\left(V_{0}^{\pi^{\star}}(\rho)-V_{0}^{\pi_j^k}(\rho)\right),\frac{1}{K}\sum_{k=0}^{K-1}\sum_{i=1}^{N}\left(\big[\ell_i-V_{i}^{\pi_j^k}(\rho)\big]_{+}+\big[V_{i}^{\pi_j^k}(\rho)-u_i\big]_{+}\right)\}\\
    &\leq\Ocal\left(\frac{N}{K\alpha}+N\eta+\frac{N}{K\eta}+\frac{N^{3/2}\alpha}{1-\sigma_2(W)}\right).
\end{align*}

Choosing the step sizes as $\alpha=\Ocal(\frac{\sqrt{1-\sigma_2(W)}}{N^{1/4}\sqrt{K}})$ and $\eta=\Ocal(1/\sqrt{K})$ leads to the claimed result.

\subsection{Proof of Theorem \ref{thm:main_ac}}\label{sec:thm_ac:proof}

Since each agent maintains and learns the local value function, the analysis of the critic is the same as in the single agent setting \cite{zeng2022finite}. Specifically, we define the critic error 
\begin{align}
    z_i^k= \widehat{Q}_{i}^k-Q_i^{\widehat{\pi}_i^k}\in\mathbb{R}^{|\Scal||\Acal|}.\label{eq:def_z}
\end{align}
and present its convergence rate in the following proposition.

\begin{prop}\label{prop:conv_critic}
Let Assumption \ref{assump:markov-chain} hold and
$\frac{(1-\gamma)\underline{\mu}\epsilon\beta}{|\Acal|}\leq 1$,
then the iterates $\{\widehat{Q}_i^k\}$ and $\{\widehat{\pi}_i^k\}$ satisfy
\begin{align*}
    \frac{1}{K}\sum_{k=0}^{K-1}\|z_i^k\|\leq
    \Ocal\Big(\frac{\tau^{1/2}}{\epsilon^{1/2}\beta^{1/2} K^{1/2}}+\frac{N\beta^{1/2}\tau}{\epsilon^{1/2}}+\frac{N\alpha}{\epsilon\beta}\Big).
\end{align*}
\end{prop}

\begin{lem}\label{lem:consensus_error_ac}
The policy iterates $\{\pi_i^k\}$ generated by Algorithm~\ref{Alg:MT-PDNAC} satisfy
\begin{align*}
    \|\widebar{\pi}^k-\pi_i^k\|\leq\Ocal\left(\frac{\sqrt{N}\alpha}{1-\sigma_2(W)}\right),\quad\text{for all $k=0,\dots,K-1$ and $i=1,\dots,N$}.
\end{align*}
\end{lem}

Similar to the proof of Theorem~\ref{thm:main_pg}, we can verify that Algorithm~\ref{Alg:MT-PDNAC} observes the update rule in \ref{prop:actor_conv:updates} and \ref{prop:actor_conv:updates_pibar}, with $Q_i^k=\widehat{Q}_i^{k}$ and $\widebar{\pi}^{k+1}=f(\{\theta_i^{k+1}\}_i)$ being defined as
\begin{align*}
    \widebar{\theta}^{k+1}=\frac{1}{N}\sum_{i=1}^{N}\theta_i^{k+1},\quad\widebar{\pi}^{k+1}=\frac{\exp \left(\widebar{\theta}^{k+1}(s,a)\right)}{\sum_{a' \in \Acal} \exp\left(\widebar{\theta}^{k+1}(s,a')\right)}
\end{align*}

The bounds in Proposition \ref{prop:actor_conv} depend on the error $\frac{1}{K}\sum_{k=0}^{K-1}\sum_{i=0}^{N}\left\|\widehat{Q}_{i}^k-Q_{i}^{\pi_i^k}\right\|$, which can be decomposed. From the definition of $z_i^k$ in Eq.~\eqref{eq:def_z}, we have for any $i=0,1,\cdots,N$
\begin{align}
    \frac{1}{K}\sum_{k=0}^{K-1}\|\widehat{Q}_{i}^k-Q_{i}^{\pi_i^k}\|&\leq \frac{1}{K}\sum_{k=0}^{K-1}\|Q_{i}^{\widehat{\pi}_i^k}-Q_{i}^{\pi_i^k}\|+\frac{1}{K}\sum_{k=0}^{K-1}\|z_i^k\|\notag\\
    &\leq \frac{1}{K}\sum_{k=0}^{K-1}\frac{|\Scal||\Acal|}{(1-\gamma)^2}\|\widehat{\pi}_i^k-\pi_i^k\|+\frac{1}{K}\sum_{k=0}^{K-1}\|z_i^k\|\notag\\
    &\leq \frac{2|\Scal|^{3/2}|\Acal| \epsilon}{(1-\gamma)^2}+\frac{1}{K}\sum_{k=0}^{K-1}\|z_i^k\|\notag\\
    &=\Ocal\left(\epsilon+\frac{\tau^{1/2}}{\epsilon^{1/2}\beta^{1/2} K^{1/2}}+\frac{N\beta^{1/2}\tau}{\epsilon^{1/2}}+\frac{N\alpha}{\epsilon\beta}\right),
    \label{thm:main:eq3}
\end{align}
where the second inequality uses Lemma \ref{lem:misc} and the third inequality follows from
\begin{align*}
    &\|\widehat{\pi}_i^k-\pi_i^k\|=\|\frac{\epsilon}{|\Acal|}\1_{|\Scal||\Acal|}+(1-\epsilon)\pi_i^k-\pi_i^k\|\notag\\
    &\leq\epsilon\|\frac{1}{|\Acal|}\1_{|\Scal||\Acal|}\|+\epsilon\|\pi_i^k\|\leq \epsilon\frac{|\Scal|^{1/2}}{|\Acal|^{1/2}}+\epsilon|\Scal|^{1/2}\hspace{-2pt}\leq\hspace{-2pt}2\epsilon|\Scal|^{1/2}.
\end{align*}

Using Eq.~\eqref{thm:main:eq3} and Lemma~\ref{lem:consensus_error_ac} in Proposition \ref{prop:actor_conv}, we obtain the following bound on the optimality gap
\begin{align*}
    &\max\{\frac{1}{K}\sum_{k=0}^{K-1}\left(V_{0}^{\pi^{\star}}(\rho)-V_{0}^{\widebar{\pi}^{k}}(\rho)\right),\frac{1}{K}\sum_{k=0}^{K-1}\sum_{i=1}^{N}\left(\big[\ell_i-V_{i}^{\widebar{\pi}^k}(\rho)\big]_{+}+\big[V_{i}^{\widebar{\pi}^k}(\rho)-u_i\big]_{+}\right)\}\\
    &\leq \Ocal\left(\frac{N}{K\alpha}+N\eta+\frac{1}{K}\sum_{k=0}^{K-1}\sum_{i=1}^{N}\|\widebar{\pi}^k-\pi_i^k\|+\frac{1}{K}\sum_{k=0}^{K-1}\sum_{i=1}^{N} \|Q_{i}^{\pi_i^k}-\widehat{Q}_{i}^{k}\|\right)\\
    &\leq \Ocal\left(\frac{N}{K\alpha}+N\eta+\frac{N}{K\eta}+\frac{N^{3/2}\alpha}{1-\sigma_2(W)}+\epsilon+\frac{\tau^{1/2}}{\epsilon^{1/2}\beta^{1/2} K^{1/2}}+\frac{N\beta^{1/2}\tau}{\epsilon^{1/2}}+\frac{N\alpha}{\epsilon\beta}\right).
\end{align*}

By the bound on consensus error in Lemma~\ref{lem:consensus_error_ac} and the Lipschitz continuity of the value function in Lemma~\ref{lem:misc}, this implies for any agent $j=1,\cdots,N$
\begin{align*}
    &\max\{\frac{1}{K}\sum_{k=0}^{K-1}\left(V_{0}^{\pi^{\star}}(\rho)-V_{0}^{\pi_j^k}(\rho)\right),\frac{1}{K}\sum_{k=0}^{K-1}\sum_{i=1}^{N}\left(\big[\ell_i-V_{i}^{\pi_j^{k}}(\rho)\big]_{+}+\big[V_{i}^{\pi_j^k}(\rho)-u_i\big]_{+}\right)\}\\
    &\leq \Ocal\left(\frac{N}{K\alpha}+N\eta+\frac{N}{K\eta}+\frac{N^{3/2}\alpha}{1-\sigma_2(W)}+\epsilon+\frac{\tau^{1/2}}{\epsilon^{1/2}\beta^{1/2} K^{1/2}}+\frac{N\beta^{1/2}\tau}{\epsilon^{1/2}}+\frac{N\alpha}{\epsilon\beta}\right).
\end{align*}

Plugging in the step sizes to the two inequalities above and recognizing from Eq.~\eqref{eq:mixing:tau} that
\begin{align*}
    \tau^{1/2} \leq \tau\leq D\log(1/\alpha)=D\log(\frac{K^{5/6}}{\alpha_0})=\Ocal(\log(K))
\end{align*}
completes the proof.

\qed

\subsection{Proof of Theorem~\ref{thm:LFA}}

In the context of linear function approximation, we denote
\begin{align}
\widehat{Q}_i^k=\Phi\omega_i^{k,T}, \quad \widehat{V}_{i}^k(s),\label{thm:LFA:eq1}
\end{align}
and adopt the rest of the notations from Section~\ref{sec:notations}.

We note that Algorithm~\ref{Alg:MT-PDNAC_nestedloop} observes the update rule in \ref{prop:actor_conv:updates} and \ref{prop:actor_conv:updates_pibar}, with $Q_i^k=\widehat{Q}_i^k$ and $\widebar{\pi}^{k+1}=f(\{\theta_i^{k+1}\}_i)$ being defined as
\begin{align*}
    \widebar{\theta}^{k+1}=\frac{1}{N}\sum_{i=1}^{N}\theta_i^{k+1},\quad\widebar{\pi}^{k+1}=\frac{\exp \left(\phi(s,a)^{\top}\widebar{\theta}^{k+1}(s,a)\right)}{\sum_{a' \in \Acal} \exp\left(\phi(s,a)^{\top}\widebar{\theta}^{k+1}(s,a')\right)}.
\end{align*}
As a result, we can apply Proposition~\ref{prop:actor_conv}, which implies
\begin{align}
    &\max\{\frac{1}{K}\sum_{k=0}^{K-1}\left(V_{0}^{\pi^{\star}}(\rho)-V_{0}^{\widebar{\pi}^{k}}(\rho)\right),\frac{1}{K}\sum_{k=0}^{K-1}\sum_{i=1}^{N}\left(\big[\ell_i-V_{i}^{\widebar{\pi}^{k}}(\rho)\big]_{+}+\big[V_{i}^{\widebar{\pi}^{k}}(\rho)-u_i\big]_{+}\right)\}\notag\\
    &\leq \Ocal\left(\frac{N}{K\alpha}+N\eta+\frac{N}{K\eta}+\frac{1}{K}\sum_{k=0}^{K-1}\sum_{i=1}^{N}\|\widebar{\pi}^k-\pi_i^k\|+\frac{1}{K}\sum_{k=0}^{K-1}\sum_{i=1}^{N} \|Q_{i}^{\pi_i^k}-\widehat{Q}_{i}^{k}\|\right).\label{thm:LFA:eq1.5}
\end{align}

A straightforward consequence of the assumption is the boundedness of $\widehat{Q}_{i}^k, \widehat{V}_{i}^k, Q_{i}^{\widehat{\pi}_i^k}$, and $\omega_{i}^{\star}(\widehat{\pi}_i^k)$, which we state in the lemma below.
\begin{lem}\label{lem:omega_bound}
Recall the definition of $B_{\omega}$ in Eq.~\eqref{eq:def_Bomega}
For all $i=1,\cdots,N$ and $k\geq 0$, we have 
\begin{gather*}
    \|\omega_i^{\star}(\widehat{\pi}_i^k)\|\leq B_{\omega},\quad \max\{\|\widehat{Q}_{i}^k\|,\|\widehat{V}_{i}^k\|,\|Q_{i}^{\widehat{\pi}_i^k}\|\}\leq Q_{\max},
\end{gather*}
where $Q_{\max}=\sigma_{\max}(\Phi)B_{\omega}$.
\end{lem}

We treat the policy space consensus error in the following lemma.
\begin{lem}\label{lem:consensus_error_LFA}
The policy iterates $\{\pi_i^k\}$ generated by Algorithm~\ref{Alg:MT-PDNAC_nestedloop} satisfy
\begin{align*}
    \|\widebar{\pi}^k-\pi_i^k\|\leq\Ocal\left(\frac{\sqrt{N}\alpha}{1-\sigma_2(W)}\right),\quad\text{for all $k=0,\dots,K-1$ and $i=1,\dots,N$}.
\end{align*}
\end{lem}

We will later decompose the error $\frac{1}{K}\sum_{k=0}^{K-1}\sum_{i=1}^{N} \|Q_{i}^{\pi_i^k}-\widehat{Q}_{i}^{k}\|$ and bound a component with the following proposition.
\begin{prop}\label{prop:critic:LFA}
Define $z_i^{k,t}=\omega_i^{k,t}-\omega_i^{\star}(\widehat{\pi}_i^k)$. Under Assumptions~\ref{assump:Slater}-\ref{assump:epsilon_max}, we have
\begin{align*}
    \mathbb{E}[\|z_i^{k,T}\|^2]\leq4B_{\omega}^2\left(1-\frac{2(1-\gamma)\underline{\mu}\sigma_{\min}(\Phi)\epsilon\beta}{|\Acal|}\right)^{T-\tau}+\frac{4(1+18B_{\omega})^2|\Acal|\beta\tau}{(1-\gamma)\underline{\mu}\sigma_{\min}(\Phi)\epsilon}.
\end{align*}
\end{prop}

By Jensen's inequality, Proposition~\ref{prop:critic:LFA} implies for all $i$ and $k$
\begin{align}
    \mathbb{E}[\|z_i^{k,T}\|]
    &\leq\sqrt{\mathbb{E}[\|z_i^{k,T}\|^2]}\notag\\
    &\leq 2B_{\omega}\left(1-\frac{2(1-\gamma)\underline{\mu}\sigma_{\min}(\Phi)\epsilon\beta}{|\Acal|}\right)^{\frac{T-\tau}{2}}+\frac{2(1+18B_{\omega})\sqrt{|\Acal|\beta\tau}}{\sqrt{(1-\gamma)\underline{\mu}\sigma_{\min}(\Phi)\epsilon}}.\label{thm:LFA:eq2}
\end{align}

By the definition of $\widehat{Q}_i^k$ in Eq.~\eqref{thm:LFA:eq1}, we have
\begin{align}
    \left\|\widehat{Q}_i^k-Q_{i}^{\pi_i^k}\right\| &\leq\left\|\Phi(\omega_i^{k,T}-\omega_i^{\star}(\widehat{\pi}_i^k))\right\|+\left\|\Phi\omega_i^{\star}(\widehat{\pi}_i^k)-Q_{i}^{\widehat{\pi}_i^k}\right\|+\left\|Q_{i}^{\widehat{\pi}_i^k}-Q_{i}^{\pi_i^k}\right\|\notag\\
    &\leq\sigma_{\max}(\Phi)\left\|\omega_i^{k,T}-\omega_i^{\star}(\widehat{\pi}_i^k)\right\|+\varepsilon_{\max}+\frac{\gamma|\Scal||\Acal|}{(1-\gamma)^2}\|\widehat{\pi}_i^k-\pi_i^k\|\notag\\
    &\leq\sigma_{\max}(\Phi)\left\|z_i^{k,T}\right\|+\varepsilon_{\max}+\frac{2\gamma|\Scal|^{3/2}|\Acal|\epsilon}{(1-\gamma)^2},\label{thm:LFA:eq3}
\end{align}
where the second inequality employs Assumption~\ref{assump:epsilon_max} and the last inequality follows from
\begin{align*}
    \|\widehat{\pi}_i^k-\pi_i^k\|&=\|\frac{\epsilon}{|\Acal|}\1_{|\Scal||\Acal|}+(1-\epsilon)\pi_i^k-\pi_i^k\|\leq\epsilon\|\frac{1}{|\Acal|}\1_{|\Scal||\Acal|}\|+\epsilon\|\pi_i^k\|\\
    &\leq \epsilon\frac{|\Scal|^{1/2}}{|\Acal|^{1/2}}+\epsilon|\Scal|^{1/2}\leq2\epsilon|\Scal|^{1/2}.
\end{align*}

As a result of Eqs.~\eqref{thm:LFA:eq2} and \eqref{thm:LFA:eq3},
\begin{align}
    \frac{1}{K}\sum_{k=0}^{K-1}\sum_{i=0}^{N}\mathbb{E}[\|Q_{i}^{\pi_i^k}-\widehat{Q}_{i}^k\|]&\leq\frac{\sigma_{\max}(\Phi)}{K} \sum_{k=0}^{K-1}\sum_{i=0}^{N}\mathbb{E}[\|z_i^{k,T}\|]+N\varepsilon_{\max}+\frac{2\gamma|\Scal|^{3/2}|\Acal|N\epsilon}{(1-\gamma)^2}\notag\\
    &\leq \Ocal\Big(N(1-\frac{2(1-\gamma)\underline{\mu}\sigma_{\min}(\Phi)\epsilon\beta}{|\Acal|})^{\frac{T-\tau}{2}}+\frac{\sqrt{\beta\tau}}{\sqrt{\epsilon}}+N\varepsilon_{\max}+N\epsilon\Big).\label{thm:LFA:eq4}
\end{align}

Plugging Eq.~\eqref{thm:LFA:eq4} and the result of Lemma~\ref{lem:consensus_error_LFA} into Eq.~\eqref{thm:LFA:eq1.5},
\begin{align*}
    &\max\left\{\frac{1}{K}\sum_{k=0}^{K-1}\left(V_{0}^{\pi^{\star}}(\rho)-V_{0}^{\widebar{\pi}^{k}}(\rho)\right),\frac{1}{K}\sum_{k=0}^{K-1}\sum_{i=1}^{N}\left(\big[\ell_i-V_{i}^{\widebar{\pi}^{k}}(\rho)\big]_{+}+\big[V_{i}^{\widebar{\pi}^{k}}(\rho)-u_i\big]_{+}\right)\right\}\notag\\
    &\leq \Ocal\hspace{-2pt}\left(\frac{N}{K\alpha}+N(\eta\hspace{-2pt}+\hspace{-2pt}\epsilon\hspace{-2pt}+\hspace{-2pt}\varepsilon_{\max})\hspace{-2pt}+\hspace{-2pt}\frac{N}{K\eta}\hspace{-2pt}+\hspace{-2pt}\frac{N^{3/2}\alpha}{1-\sigma_2(W)}+N(1-\frac{2(1-\gamma)\underline{\mu}\sigma_{\min}(\Phi)\epsilon\beta}{|\Acal|})^{\frac{T-\tau}{2}}+\frac{\sqrt{\beta\tau}}{\sqrt{\epsilon}}\right).
\end{align*}

By the bound on consensus error in Lemma~\ref{lem:consensus_error_LFA} and the Lipschitz continuity of the value function in Lemma~\ref{lem:misc}, this implies for any agent $j=1,\cdots,N$
\begin{align*}
    &\max\left\{\frac{1}{K}\sum_{k=0}^{K-1}\left(V_{0}^{\pi^{\star}}(\rho)-V_{0}^{\pi_j^{k}}(\rho)\right),\frac{1}{K}\sum_{k=0}^{K-1}\sum_{i=1}^{N}\left(\big[\ell_i-V_{i}^{\pi_j^{k}}(\rho)\big]_{+}+\big[V_{i}^{\pi_j^{k}}(\rho)-u_i\big]_{+}\right)\right\}\notag\\
    &\leq \Ocal\left(\frac{N}{K\alpha}\hspace{-2pt}+\hspace{-2pt}N(\eta+\epsilon+\varepsilon_{\max})\hspace{-2pt}+\hspace{-2pt}\frac{N}{K\eta}\hspace{-2pt}+\hspace{-2pt}\frac{N^{3/2}\alpha}{1-\sigma_2(W)}\hspace{-2pt}+\hspace{-2pt}N(1-\frac{2(1-\gamma)\underline{\mu}\sigma_{\min}(\Phi)\epsilon\beta}{|\Acal|})^{\frac{T-\tau}{2}}+\frac{\sqrt{\beta\tau}}{\sqrt{\epsilon}}\right).
\end{align*}

In order to have
\begin{align*}
    \max\left\{\frac{1}{K}\sum_{k=0}^{K-1}\left(V_{0}^{\pi^{\star}}(\rho)-V_{0}^{\pi_j^{k}}(\rho)\right),\frac{1}{K}\sum_{k=0}^{K-1}\sum_{i=1}^{N}\left(\big[\ell_i-V_{i}^{\pi_j^{k}}(\rho)\big]_{+}+\big[V_{i}^{\pi_j^{k}}(\rho)-u_i\big]_{+}\right)\right\}\\
    \leq \frac{N^{5/4}\delta}{\sqrt{1-\sigma_2(W)}}+\Ocal(N\varepsilon_{\max}),
\end{align*}
we can choose
\begin{align*}
    \alpha& \sim \Ocal(\frac{(\sqrt{1-\sigma_2(W)}\delta}{N^{1/4}}),\quad
    \beta \sim \Ocal(\frac{\delta^3}{\log(1/\delta)}),\quad\epsilon \sim\Ocal(\delta),\quad
    \eta\sim\Ocal(\delta), 
\end{align*}
which implies
\begin{gather*}
    K\sim\Ocal(\delta^{-2}),\quad T\sim\Ocal(\frac{\log(N/\delta)}{\beta\epsilon})=\Ocal(\frac{\log(N/\delta)}{\delta^4}),\quad TK\sim\Ocal(\frac{\log(1/\delta)}{\delta^6}).
\end{gather*}

\qed

\section{Proof of Propositions}

\subsection{Proof of Proposition~\ref{prop:actor_conv}}

Without loss of generality, we assume that $B\geq\frac{1}{1-\gamma}$, since otherwise we can safely set $B$ to be $\frac{1}{1-\gamma}$ in Proposition~\ref{prop:actor_conv}. This helps us simplify notations later in the analysis.

We define the following notation on the aggregate policy over the network.
\begin{align}
    \vpi^k=[\pi_1^k,\dots,\pi_N^k]\label{eq:def_aggregate_notation}
\end{align}

We have from Eq.~\eqref{prop:actor_conv:updates_pibar}
\begin{align}
\widebar{\pi}^{k+1}(a\mid s)=\widebar{\pi}^{k}(a\mid s)\frac{\exp(\frac{\alpha}{N}\sum_{i=1}^{N}(\frac{1}{N}+\lambda_i^k-\nu_i^k)Q_i^{k}(s,a))}{Z^k(s)},\label{prop:actor_conv:updates_rewrite}
\end{align}
where $Z^k(s)=\sum_{a'\in\Acal}\widebar{\pi}^{k}(a'\mid s)\exp(\frac{\alpha}{N}\sum_{i=1}^{N}(\frac{1}{N}+\lambda_i^k-\nu_i^k)Q_i^{k}(s,a'))$.

Let $V_i^k(s)=\sum_{a}\pi_i^k(a\mid s)Q_i^k(s,a)$.

Define $Q_{g}^{\vpi^k}=[Q_1^{\pi_1^k},\dots,Q_N^{\pi_N^k}]$ and $V_{g}^{\vpi^k}=[V_1^{\pi_1^k},\dots,V_N^{\pi_N^k}]$.
Define $Q_{g}^{k}=[Q_1^{k},\dots,Q_N^{k}]$ and $V_{g}^{k}=[V_1^{k},\dots,V_N^{k}]$.
Define $Q_{L,k}^{k}=\sum_{i=1}^{N}(\frac{1}{N}+\lambda_i^k-\nu_i^k)Q_i^{k}$ and $V_{L,k}^{k}=\sum_{i=1}^{N}(\frac{1}{N}+\lambda_i^k-\nu_i^k)V_i^{k}$.

\noindent\textbf{Objective function convergence.}
From the dual update in Eq.~\eqref{Alg:MT-PDNPG:dual_update}, we have
\begin{align}
    0\leq\|\lambda^K\|^2&=\sum_{k=0}^{K-1}\left(\|\lambda^{k+1}\|^2-\|\lambda^k\|^2\right)\notag\\
    &=\sum_{k=0}^{K-1}\left(\left\|\Pi_{[0,B_{\lambda}]}\left(\lambda^k-\eta\left(\sum_{s,a}\rho(s)\diag(\vpi^k(a\mid s))Q_{g}^{k}(s,a)-\ell\right)\right)\right\|^2-\|\lambda^k\|^2\right)\notag\\
    &\leq\sum_{k=0}^{K-1}\left(\left\|\lambda^k-\eta\left(\sum_{s,a}\rho(s)\diag(\vpi^k(a\mid s))Q_{g}^{k}(s,a)-\ell\right)\right\|^2-\|\lambda^k\|^2\right)\notag\\
    &=-2\eta\sum_{k=0}^{K-1}(\lambda^k)^{\top}\left(\sum_{s,a}\rho(s)\diag(\vpi^k(a\mid s))Q_{g}^{\vpi^k}(s,a)-\ell\right)\notag\\
    &\hspace{20pt}+2\eta\sum_{k=0}^{K-1}(\lambda^k)^{\top}\left(\sum_{s,a}\rho(s)\diag(\vpi^k(a\mid s))\left(Q_{g}^{\vpi^k}(s,a)-Q_{g}^{k}(s,a)\right)\right)\notag\\
    &\hspace{20pt}+\eta^2\sum_{k=0}^{K-1}\left\|\sum_{s,a}\rho(s)\diag(\vpi^k(a\mid s))Q_{g}^{k}(s,a)-\ell\right\|^2\notag\\
    &=-2\eta\sum_{k=0}^{K-1}(\lambda^k)^{\top}\left(V_{g}^{\vpi^k}(\rho)-\ell\right)+\eta^2\sum_{k=0}^{K-1}\left\|\sum_{s,a}\rho(s)\diag(\vpi^k(a\mid s))Q_{g}^{k}(s,a)-\ell\right\|^2\notag\\
    &\hspace{20pt}+2\eta\sum_{k=0}^{K-1}(\lambda^k)^{\top}\left(\sum_{s,a}\rho(s)\diag(\vpi^k(a\mid s))\left(Q_{g}^{\vpi^k}(s,a)-Q_{g}^{k}(s,a)\right)\right).\label{prop:convergence_obj_constraintviolation:eq0}
\end{align}

Since the value function and constant $\ell_i$ are within $[0,\frac{1}{1-\gamma}]$, the second term of Eq.~\eqref{prop:convergence_obj_constraintviolation:eq0} obeys 
\begin{align}
    \sum_{k=0}^{K-1}\left\|\sum_{s,a}\rho(s)\diag(\vpi^k(a\mid s))Q_{g}^{k}(s,a)-\ell\right\|^2&= \sum_{k=0}^{K-1}\sum_{i=1}^{N}\left(\sum_{s,a}\rho(s)\pi_i^k(a\mid s)Q_{i}^{k}(s,a)-\ell_i\right)^2\notag\\
    &\leq 2\sum_{k=0}^{K-1}\sum_{i=1}^{N}\left(\left(\sum_{s,a}\rho(s)\pi_i^k(a\mid s)Q_{i}^{k}(s,a)\right)^2+\left(\ell_i\right)^2\right)\notag\\
    &\leq 2\sum_{k=0}^{K-1}\sum_{i=1}^{N}\left(B^2+\frac{1}{(1-\gamma)^2}\right)\leq\frac{4KN}{(1-\gamma)^2}.\label{prop:convergence_obj_constraintviolation:eq0.5}
\end{align}

The third term of Eq.~\eqref{prop:convergence_obj_constraintviolation:eq0} can be treated as
\begin{align}
    &2\eta\sum_{k=0}^{K-1}(\lambda^k)^{\top}\left(\sum_{s,a}\rho(s)\diag(\vpi^k(a\mid s))\left(Q_{g}^{\vpi^k}(s,a)-Q_{g}^{k}(s,a)\right)\right)\notag\\
    &=2\eta\sum_{k=0}^{K-1}\sum_{i=1}^{N}\lambda_i^k\left(\sum_{s,a}\rho(s)\pi_i^k(a\mid s)\left(Q_{i}^{\pi_i^k}(s,a)-Q_{i}^{k}(s,a)\right)\right)\notag\\
    &\leq 2B_{\lambda}\eta\sum_{k=0}^{K-1}\sum_{i=1}^{N} \left(\sum_{s,a}\rho(s)^2\pi_i^k(a\mid s)^2\right)^{1/2}\|Q_{i}^{\pi_i^k}-Q_{i}^{k}\|\notag\\
    &\leq 2B_{\lambda}\eta\sum_{k=0}^{K-1}\sum_{i=1}^{N} \|Q_{i}^{\pi_i^k}-Q_{i}^{k}\|,\label{prop:convergence_obj_constraintviolation:eq0.6}
\end{align}
where the first inequality follows from the Cauchy-Schwarz inequality, and the second inequality follows from the fact that the $\ell_2$ norm of a vector is upper bounded by its $\ell_1$ norm.

Using Eqs.~\eqref{prop:convergence_obj_constraintviolation:eq0.5} and \eqref{prop:convergence_obj_constraintviolation:eq0.6} in Eq.~\eqref{prop:convergence_obj_constraintviolation:eq0}, we get
\begin{align*}
    &0\leq-2\eta\sum_{k=0}^{K-1}(\lambda^k)^{\top}\left(V_{g}^{\vpi^k}(\rho)-\ell\right)+\eta^2\sum_{k=0}^{K-1}\left\|\sum_{s,a}\rho(s)\diag(\vpi^k(a\mid s))Q_{g}^{k}(s,a)-b\right\|^2\notag\\
    &\hspace{20pt}+2\eta\sum_{k=0}^{K-1} (\lambda^k)^{\top}\left(\sum_{s,a}\rho(s)\diag(\vpi^k(a\mid s))\left(Q_{g}^{\vpi^k}(s,a)-Q_{g}^{k}(s,a)\right)\right)\notag\\
    &\leq 2\eta\sum_{k=0}^{K-1}(\lambda^k)^{\top} \left(V_{g}^{\pi^{\star}}(\rho)\hspace{-2pt}-\hspace{-2pt}V_{g}^{\widebar{\pi}^k}(\rho)\right)+2\eta\sum_{k=0}^{K-1}(\lambda^k)^{\top} \left(V_{g}^{\widebar{\pi}^k}(\rho)\hspace{-2pt}-\hspace{-2pt}V_{g}^{\vpi^k}(\rho)\right)\\
    &\hspace{20pt}+\frac{4KN\eta^2}{(1-\gamma)^2}+2B_{\lambda}\eta\sum_{k=0}^{K-1}\sum_{i=1}^{N} \|Q_{i}^{\pi_i^k}\hspace{-2pt}-\hspace{-2pt}Q_{i}^{k}\|\notag\\
    &\leq 2\eta\sum_{k=0}^{K-1}\hspace{-2pt}(\lambda^k)^{\top} \hspace{-2pt}\left(V_{g}^{\pi^{\star}}(\rho)\hspace{-2pt}-\hspace{-2pt}V_{g}^{\widebar{\pi}^k}(\rho)\right)\hspace{-2pt}+\hspace{-2pt}\frac{2\sqrt{|\Scal|| \Acal|}B_{\lambda}\eta}{(1-\gamma)^2}\sum_{k=0}^{K-1}\sum_{i=1}^{N}\|\widebar{\pi}^k-\pi_i^k\| \\
    &\hspace{20pt}+\frac{4KN\eta^2}{(1-\gamma)^2}+2B_{\lambda}\eta\sum_{k=0}^{K-1}\sum_{i=1}^{N} \|Q_{i}^{\pi_i^k}-Q_{i}^{k}\|,
\end{align*}
where the second inequality follows from the fact that the optimal policy satisfies the constraints, i.e. $V_i^{\pi^{\star}}(\rho)\geq \ell_i$ for all $i=1,\cdots,N$, and the third inequality is applies Lemma~\ref{lem:misc}.

Re-arranging this inequality and dividing by $2K\eta$ lead to
\begin{align}
    &\frac{1}{K}\sum_{k=0}^{K-1}(\lambda^k)^{\top} \left(V_{g}^{\pi^{\star}}(\rho)-V_{g}^{\widebar{\pi}^k}(\rho)\right)\notag\\
    &\hspace{50pt}\geq-\frac{\sqrt{|\Scal||\Acal|}B_{\lambda}}{(1-\gamma)^2 K}\sum_{k=0}^{K-1}\sum_{i=1}^{N}\|\widebar{\pi}^k-\pi_i^k\| -\frac{2N\eta}{(1-\gamma)^2}-\frac{B_{\lambda}}{K}\sum_{k=0}^{K-1}\sum_{i=1}^{N} \|Q_{i}^{\pi_i^k}-Q_{i}^{k}\|.\label{prop:convergence_obj_constraintviolation:eq1}
\end{align}

A similar analysis on $\nu^k$ implies
\begin{align}
    &-\frac{1}{K}\sum_{k=0}^{K-1}(\nu^k)^{\top} \left(V_{g}^{\pi^{\star}}(\rho)-V_{g}^{\widebar{\pi}^k}(\rho)\right)\notag\\
    &\hspace{50pt}\geq-\frac{\sqrt{|\Scal||\Acal|}B_{\lambda}}{(1-\gamma)^2 K}\sum_{k=0}^{K-1}\sum_{i=1}^{N}\|\widebar{\pi}^k-\pi_i^k\| -\frac{2N\eta}{(1-\gamma)^2}-\frac{B_{\lambda}}{K}\sum_{k=0}^{K-1}\sum_{i=1}^{N} \|Q_{i}^{\pi_i^k}-Q_{i}^{k}\|.\label{prop:convergence_obj_constraintviolation:eq1.5}
\end{align}

\begin{lem}\label{lem:step_improvement}
The iterates of Eq.~\eqref{prop:actor_conv:updates} satisfy for all $k=0,\dots,K-1$
\begin{align*}
    &V_{L,k}^{\widebar{\pi}^{k+1}}(\zeta)-V_{L,k}^{\widebar{\pi}^{k}}(\zeta)\geq \frac{N}{\alpha}\mathbb{E}_{s \sim \zeta}\left[\log Z_k(s)-\frac{\alpha}{N}V_{L,k}^{\widebar{\pi}^k}(s)\right]\notag\\
    &\hspace{80pt}-\frac{2(B_{\lambda}+1/N)}{1-\gamma}\sum_{i=1}^{N}\|Q_i^{\pi_i^k}-Q_i^k\|-\frac{2\sqrt{|\Scal|| \Acal|}(B_{\lambda}+1/N)}{(1-\gamma)^3}\sum_{i=1}^{N}\|\widebar{\pi}^k-\pi_i^k\|.
\end{align*}
\end{lem}

\begin{lem}\label{lem:bounded_Lagrangian}
The iterates of Eq.~\eqref{prop:actor_conv:updates} satisfy for all $k=0,\dots,K-1$
\begin{align*}
    \frac{1}{K}\sum_{k=0}^{K-1}\left(V_{L,k}^{\pi^{\star}}(\rho)-V_{L,k}^{\widebar{\pi}^{k}}(\rho)\right)&\leq \frac{N\log|\Acal|}{(1-\gamma)K\alpha}+\frac{3(B_{\lambda}+1/N)}{(1-\gamma)^2 K} \sum_{k=0}^{K-1}\sum_{i=1}^{N}\|Q_i^{\pi_i^k}-Q_i^k\|\notag\\
    &\hspace{-20pt}+\frac{3\sqrt{|\Scal|| \Acal|}(B_{\lambda}+1/N)}{(1-\gamma)^4 K}\sum_{k=0}^{K-1}\sum_{i=1}^{N}\|\widebar{\pi}^k-\pi_i^k\|+\frac{2N B_{\lambda}}{(1-\gamma)^2 K}+\frac{4N\eta}{(1-\gamma)^3 K}.
\end{align*}
\end{lem}

Combining Eq.~\eqref{prop:convergence_obj_constraintviolation:eq1}, Eq.~\eqref{prop:convergence_obj_constraintviolation:eq1.5}, and Lemma~\ref{lem:bounded_Lagrangian},
\begin{align*}
    &\frac{1}{K}\sum_{k=0}^{K-1}\left(V_{0}^{\pi^{\star}}(\rho)-V_{0}^{\widebar{\pi}^{k}}(\rho)\right)\\
    &\leq\frac{N\log|\Acal|}{(1-\gamma)K\alpha}+\frac{3(B_{\lambda}+1/N)}{(1-\gamma)^2 K} \sum_{k=0}^{K-1}\sum_{i=1}^{N}\|Q_i^{\pi_i^k}-Q_i^k\|\notag\\
    &\hspace{20pt}+\frac{3\sqrt{|\Scal|| \Acal|}(B_{\lambda}+1/N)}{(1-\gamma)^4 K}\sum_{k=0}^{K-1}\sum_{i=1}^{N}\|\widebar{\pi}^k-\pi_i^k\|+\frac{2N B_{\lambda}}{(1-\gamma)^2 K}+\frac{4N\eta}{(1-\gamma)^3 K}\notag\\
    &\hspace{20pt}+\frac{2\sqrt{|\Scal||\Acal|}B_{\lambda}}{(1-\gamma)^2 K}\sum_{k=0}^{K-1}\sum_{i=1}^{N}\|\widebar{\pi}^k-\pi_i^k\| +\frac{4N\eta}{(1-\gamma)^2}+\frac{2B_{\lambda}}{K}\sum_{k=0}^{K-1}\sum_{i=1}^{N} \|Q_{i}^{\pi_i^k}-Q_{i}^{k}\|\notag\\
    &\leq \frac{N\log|\Acal|}{(1-\gamma)K\alpha}\hspace{-2pt}+\hspace{-2pt}\frac{2N B_{\lambda}}{(1-\gamma)^2 K}\hspace{-2pt}+\hspace{-2pt}\frac{8N\eta}{(1-\gamma)^3}\hspace{-2pt}+\hspace{-2pt}\frac{5\sqrt{|\Scal|| \Acal|}(B_{\lambda}+1/N)}{(1-\gamma)^4 K}\sum_{k=0}^{K-1}\sum_{i=1}^{N}\|\widebar{\pi}^k-\pi_i^k\|\notag\\
    &\hspace{20pt}+\frac{8(B_{\lambda}+1/N)}{(1-\gamma)^2 K}\sum_{k=0}^{K-1}\sum_{i=1}^{N} \|Q_{i}^{\pi_i^k}-Q_{i}^{k}\|.
\end{align*}

\textbf{Constraint violation convergence.}
For any $\lambda\in[0,B_{\lambda}]^{N}$, since the projection operator $\Pi_{[0,B_{\lambda}]}$ is non-expansive, we have
\begin{align*}
    \|\lambda^{k+1}-\lambda\|^2 &= \|\Pi_{[0,B_{\lambda}]}(\lambda^k-\eta(V_{g}^{k}(\rho)-\ell))-\lambda\|^2\notag\\
    &\leq \|\lambda^k-\eta(V_{g}^{k}(\rho)-\ell)-\lambda\|^2\notag\\
    &= \|\lambda^k-\lambda\|^2-2\eta(\lambda^k-\lambda)^{\top}(V_{g}^{k}(\rho)-\ell)+\eta^2\sum_{i=1}^{N}\|\sum_{s,a}\rho(s)\pi_i^k(a\mid s)Q_{i}^{k}(\rho)-\ell_i\|^2\notag\\
    &\leq \|\lambda^k-\lambda\|^2-2\eta(\lambda^k-\lambda)^{\top}(V_{g}^{k}(\rho)-\ell)+\frac{4N\eta^2}{(1-\gamma)^2},
\end{align*}
where the last inequality bounds the quadratic term using an approach similar to Eq.~\eqref{prop:convergence_obj_constraintviolation:eq0.5}.

Re-arranging the terms and summing up from $k=0$ to $k=K-1$, we get
\begin{align*}
    \frac{1}{K}\sum_{k=0}^{K-1}(\lambda^k-\lambda)^{\top}(V_{g}^{k}(\rho)-b)
    &\leq\frac{1}{K}\left(\|\lambda^0-\lambda\|^2-\|\lambda^K-\lambda\|^2\right)+\frac{2N\eta}{(1-\gamma)^2}\notag\\
    &\leq \frac{1}{2K\eta}\|\lambda^0-\lambda\|^2+\frac{2N\eta}{(1-\gamma)^2},
\end{align*}
which implies
\begin{align*}
    &\frac{1}{K}\sum_{k=0}^{K-1}(\lambda^k-\lambda)^{\top}(V_{g}^{\widebar{\pi}^k}(\rho)-\ell) \notag\\
    &=\frac{1}{K}\sum_{k=0}^{K-1}(\lambda^k-\lambda)^{\top}(V_{g}^{k}(\rho)-\ell)+\frac{1}{K}\sum_{k=0}^{K-1}(\lambda^k-\lambda)^{\top}(V_{g}^{\widebar{\pi}^k}(\rho)-V_{g}^{\vpi_k}(\rho))\notag\\
    &\hspace{20pt}+\frac{1}{K}\sum_{k=0}^{K-1}(\lambda^k-\lambda)^{\top}(V_{g}^{\vpi^k}(\rho)-V_{g}^{k}(\rho))\notag\\
    &\leq \frac{1}{2K\eta}\|\lambda^0-\lambda\|^2+\frac{2N\eta}{(1-\gamma)^2}+\frac{2\sqrt{|\Scal|| \Acal|}B_{\lambda}}{(1-\gamma)^2 K}\sum_{k=0}^{K-1}\sum_{i=1}^{N}\|\widebar{\pi}^k-\pi_i^k\|+\frac{2B_{\lambda}}{K}\sum_{k=0}^{K-1}\sum_{i=1}^{N}\|Q_i^k-Q_i^{\pi_i^k}\|.
\end{align*}

Similarly, we can show for any $\nu\in[0,B_{\lambda}]^{N}$
\begin{align*}
    \frac{1}{K}\sum_{k=0}^{K-1}(\lambda^k-\lambda)^{\top}(u-V_{g}^{\widebar{\pi}^k}(\rho))
    &\leq \frac{1}{2K\eta}\|\nu^0-\nu\|^2+\frac{2N\eta}{(1-\gamma)^2}+\frac{2\sqrt{|\Scal|| \Acal|}B_{\lambda}}{(1-\gamma)^2 K}\sum_{k=0}^{K-1}\sum_{i=1}^{N}\|\widebar{\pi}^k-\pi_i^k\|\notag\\
    &\hspace{20pt}+\frac{2B_{\lambda}}{K}\sum_{k=0}^{K-1}\sum_{i=1}^{N}\|Q_i^k-Q_i^{\pi_i^k}\|.
\end{align*}

Since 
$\lambda^k$, $\nu^k$ are non-negative, we have 
from Lemma \ref{lem:bounded_Lagrangian} and the two inequalities above
\begin{align}
    &\frac{1}{K}\sum_{k=0}^{K-1}\left(V_{0}^{\pi^{\star}}(\rho)-V_{0}^{\widebar{\pi}^{k}}(\rho)+\lambda^{\top}(\ell-V_g^{\widebar{\pi}^k}(\rho))+\nu^{\top}(V_g^{\widebar{\pi}^k}(\rho)-u)\right)\notag\\
    &\leq\frac{1}{K}\hspace{-2pt}\sum_{k=0}^{K-1}\hspace{-3pt}\left(V_{0}^{\pi^{\star}}\hspace{-2pt}(\rho)\hspace{-3pt}-\hspace{-3pt}V_{0}^{\widebar{\pi}^{k}}\hspace{-2pt}(\rho)\hspace{-3pt}+\hspace{-3pt}(\lambda^k)^{\top}\hspace{-3pt}(V_{g}^{\pi^{\star}}\hspace{-2pt}(\rho)\hspace{-3pt}-\hspace{-3pt}\ell)\hspace{-3pt}+\hspace{-3pt}\lambda^{\top}\hspace{-2pt}(\ell\hspace{-3pt}-\hspace{-3pt}V_g^{\widebar{\pi}^k}\hspace{-2pt}(\rho))\hspace{-3pt}+\hspace{-3pt}(\nu^k)^{\top}\left(u\hspace{-3pt}-\hspace{-3pt}V_{g}^{\pi^{\star}}\hspace{-2pt}(\rho)\right)\hspace{-3pt}+\hspace{-3pt}\nu^{\top}\hspace{-2pt}(V_g^{\widebar{\pi}^k}\hspace{-3pt}(\rho)\hspace{-3pt}-\hspace{-3pt}u)\right)\notag\\
    &=\frac{1}{K}\sum_{k=0}^{K-1}\hspace{-2pt}\Big(V_{0}^{\pi^{\star}}(\rho)-V_{0}^{\widebar{\pi}^{k}}(\rho)+(\lambda^k-\nu^k)^{\top}\left(V_{g}^{\pi^{\star}}(\rho)-V_{g}^{\widebar{\pi}^k}(\rho)\right)\notag\\
    &\hspace{120pt}+(\lambda^k-\lambda)^{\top}(V_g^{\widebar{\pi}^k}(\rho)-\ell)+(\nu^k-\nu)^{\top}(u-V_g^{\widebar{\pi}^k}(\rho))\Big)\notag\\
    &\leq \frac{N\log|\Acal|}{(1-\gamma)K\alpha}+\frac{3(B_{\lambda}+1/N)}{(1-\gamma)^2 K} \sum_{k=0}^{K-1}\sum_{i=1}^{N}\|Q_i^{\pi_i^k}-Q_i^k\|\notag\\
    &\hspace{20pt}+\frac{3\sqrt{|\Scal|| \Acal|}(B_{\lambda}+1/N)}{(1-\gamma)^4 K}\sum_{k=0}^{K-1}\sum_{i=1}^{N}\|\widebar{\pi}^k-\pi_i^k\|+\frac{2N B_{\lambda}}{(1-\gamma)^2 K}+\frac{4N\eta}{(1-\gamma)^3 K}\notag\\
    &\hspace{20pt}+\frac{1}{2K\eta}\|\lambda^0-\lambda\|^2+\frac{2N\eta}{(1-\gamma)^2}+\frac{2\sqrt{|\Scal|| \Acal|}B_{\lambda}}{(1-\gamma)^2 K}\sum_{k=0}^{K-1}\sum_{i=1}^{N}\|\widebar{\pi}^k-\pi_i^k\|+\frac{2B_{\lambda}}{K}\sum_{k=0}^{K-1}\sum_{i=1}^{N}\|Q_i^k-Q_i^{\pi_i^k}\|\notag\\
    &\hspace{20pt}+\frac{1}{2K\eta}\|\nu^0-\nu\|^2+\frac{2N\eta}{(1-\gamma)^2}+\frac{2\sqrt{|\Scal|| \Acal|}B_{\lambda}}{(1-\gamma)^2 K}\sum_{k=0}^{K-1}\sum_{i=1}^{N}\|\widebar{\pi}^k-\pi_i^k\|+\frac{2B_{\lambda}}{K}\sum_{k=0}^{K-1}\sum_{i=1}^{N}\|Q_i^k-Q_i^{\pi_i^k}\|\notag\\
    &\leq \frac{N\log|\Acal|}{(1-\gamma)K\alpha}+\frac{2N B_{\lambda}}{(1-\gamma)^2 K}+\frac{8N\eta}{(1-\gamma)^3}+\frac{\|\lambda^0-\lambda\|^2+\|\nu^0-\nu\|^2}{2K\eta}\label{eq:constraint_violation_4}\\
    &\hspace{20pt}+\frac{7\sqrt{|\Scal|| \Acal|}(B_{\lambda}+1/N)}{(1-\gamma)^4 K}\sum_{k=0}^{K-1}\sum_{i=1}^{N}\|\widebar{\pi}^k-\pi_i^k\|+\frac{7(B_{\lambda}+1/N)}{(1-\gamma)^2 K} \sum_{k=0}^{K-1}\sum_{i=1}^{N}\|Q_i^{\pi_i^k}-Q_i^k\|.\notag
\end{align}

Now, choosing $\lambda$ and $\nu$ such that 
\begin{align*}
    \lambda_i = \begin{cases}B_{\lambda}, & \text { if } \ell_i-V_i^{\pi_k}(\rho)\geq 0 \\ 0, & \text { else }\end{cases}\quad \nu_i = \begin{cases}B_{\lambda}, & \text { if } V_i^{\pi_k}(\rho)-u_i\geq 0 \\ 0, & \text { else }\end{cases}
\end{align*}

Then, Eq.~\eqref{eq:constraint_violation_4} leads to
\begin{align}
    &\frac{1}{K}\sum_{k=0}^{K-1}\left(V_{0}^{\pi^{\star}}(\rho)-V_{0}^{\widebar{\pi}^k}(\rho)\right)+\frac{1}{K}\sum_{k=0}^{K-1}\sum_{i=1}^{N}B_{\lambda}\left(\big[\ell_i-V_{i}^{\widebar{\pi}^k}(\rho)\big]_{+}+\big[V_{i}^{\widebar{\pi}^k}(\rho)-u_i\big]_{+}\right)\notag\\
    &\leq \frac{N\log|\Acal|}{(1-\gamma)K\alpha}+\frac{2N B_{\lambda}}{(1-\gamma)^2 K}+\frac{8N\eta}{(1-\gamma)^3}+\frac{NB_{\lambda}^2}{K\eta}\notag\\
    &\hspace{20pt} +\frac{7\sqrt{|\Scal|| \Acal|}(B_{\lambda}+1/N)}{(1-\gamma)^4 K}\sum_{k=0}^{K-1}\sum_{i=1}^{N}\|\widebar{\pi}^k-\pi_i^k\|+\frac{7(B_{\lambda}+1/N)}{(1-\gamma)^2 K} \sum_{k=0}^{K-1}\sum_{i=1}^{N}\|Q_i^{\pi_i^k}-Q_i^k\|.\label{eq:constraint_violation_5}
\end{align}

Note that there always exists a policy $\widetilde{\pi}^K$ such that $d_{\rho}^{\widetilde{\pi}^K}=\frac{1}{K}\sum_{k=0}^{K-1}d_{\rho}^{\widebar{\pi}^k}$, which implies
\begin{align*}
    V_i^{\widetilde{\pi}^K}=\frac{1}{K}\sum_{k=0}^{K-1} V_i^{\widebar{\pi}^k}\quad\forall i=0,1,\cdots,N.
\end{align*}

As a result, Eq.~\eqref{eq:constraint_violation_5} becomes
\begin{align}
    &\left(V_{0}^{\pi^{\star}}(\rho)-V_{0}^{\widetilde{\pi}^K}(\rho)\right)+B_{\lambda}\sum_{i=1}^{N}\left(\big[\ell_i-V_{i}^{\widetilde{\pi}^K}(\rho)\big]_{+}+\big[V_{i}^{\widetilde{\pi}^K}(\rho)-u_i\big]_{+} \right)\notag\\
    &\leq \frac{N\log|\Acal|}{(1-\gamma)K\alpha}+\frac{2N B_{\lambda}}{(1-\gamma)^2 K}+\frac{8N\eta}{(1-\gamma)^3}+\frac{NB_{\lambda}^2}{K\eta}\notag\\
    &\hspace{20pt} +\frac{7\sqrt{|\Scal|| \Acal|}(B_{\lambda}+1/N)}{(1-\gamma)^4 K}\sum_{k=0}^{K-1}\sum_{i=1}^{N}\|\widebar{\pi}^k-\pi_i^k\|+\frac{7(B_{\lambda}+1/N)}{(1-\gamma)^2 K} \sum_{k=0}^{K-1}\sum_{i=1}^{N}\|Q_i^{\pi_i^k}-Q_i^k\|.\label{eq:constraint_violation_6}
\end{align}

\begin{lem}[Theorem 6 of \cite{ding2020natural}]\label{lem:constraint_violation}
Suppose that Assumption \ref{assump:Slater} holds.
Let the constant $C$ obey $C\geq 2\|\lambda^{\star}\|_{\infty}$ and $C\geq 2\|\nu^{\star}\|_{\infty}$.
Then, given a policy $\pi$, if there exists a constant $\delta>0$ such that
\begin{align*}
    V_{0}^{\pi^{\star}}(\rho)-V_{0}^{\pi}(\rho)+C\sum_{i=1}^{N}\left([\ell_i-V_i^{\pi}(\rho)]_{+}+[V_i^{\pi}(\rho)-u_i]_{+}\right)\leq\delta,
\end{align*}
then we have
\begin{align*}
    \sum_{i=1}^{N}\left([\ell_i-V_i^{\pi}(\rho)]_{+}+[V_i^{\pi}(\rho)-u_i]_{+}\right)\leq \frac{2\delta}{C}.
\end{align*}
\end{lem}

Recall that Lemma \ref{lem:bounded_lambdastar} states that $2\|\lambda^{\star}\|_{\infty}\leq B_{\lambda}$ and $2\|\nu^{\star}\|_{\infty}\leq B_{\lambda}$. Applying Lemma \ref{lem:constraint_violation} with $C=B_{\lambda}$ and $\delta$ being the terms on the left hand side of Eq.~\eqref{eq:constraint_violation_6}, we have
\begin{align*}
    \sum_{i=1}^{N}\left(\big[\ell_i-V_{i}^{\widetilde{\pi}^{K}}(\rho)\big]_{+}+\big[V_{i}^{\widetilde{\pi}^{K}}(\rho)-u_i\big]_{+}\right)\leq\hspace{-2pt}\frac{2}{B_{\lambda}}\Big(\frac{N\log|\Acal|}{(1-\gamma)K\alpha}+\frac{2N B_{\lambda}}{(1-\gamma)^2 K}+\frac{8N\eta}{(1-\gamma)^3}+\frac{NB_{\lambda}^2}{K\eta}&\notag\\
    +\frac{7\sqrt{|\Scal|| \Acal|}(B_{\lambda}+1/N)}{(1-\gamma)^4 K}\sum_{k=0}^{K-1}\sum_{i=1}^{N}\|\widebar{\pi}^k-\pi_i^k\|+\frac{7(B_{\lambda}+1/N)}{(1-\gamma)^2 K} \sum_{k=0}^{K-1}\sum_{i=1}^{N}\|Q_i^{\pi_i^k}-Q_i^k\|\Big).&
\end{align*}

\qed

\subsection{Proof of Proposition~\ref{prop:conv_critic}}

We omit the proof and note that this result is adapted from \cite{zeng2022finite}[Proposition 2].

\subsection{Proof of Proposition~\ref{prop:critic:LFA}}

This section presents the proof of Proposition~\ref{prop:conv_critic}.
We use $O_i^{k,t}$ to denote the data observation used for variable updates in iteration $(t,k)$, i.e.
\[O_i^{k,t}=(s_i^{k,t},a_i^{k,t},s_i^{k,t+1},a_i^{k,t+1}).\]

We use $R_i:|\Scal|\times|\Acal|\times|\Scal|\times|\Acal|\rightarrow\mathbb{R}^{d}$ to denote the feature-reward composite operator such that
\begin{align}
    R_i(s,a,s',a')= r_i(s,a)\phi(s,a).
    \label{eq:def_vectorR}
\end{align}

We also define the operator $H:|\Scal|\times|\Acal|\times|\Scal|\times|\Acal|\rightarrow\mathbb{R}^{d\times d}$ such that
\begin{align}
    H(s,a,s',a')=\phi(s,a)(\gamma\phi(s',a')-\phi(s,a))^{\top},
    \label{eq:def_matrixA}
\end{align}
which means that $\widebar{H}^{\pi}$ defined in Eq.~\eqref{eq:def_Ab} satisfies
\[\widebar{H}^{\pi}=\mathbb{E}_{s\sim\mu_{\pi},a\sim\pi(\cdot\mid s),s'\sim P(\cdot\mid s,a),a'\sim\pi(\cdot\mid s')}[H(s,a,s',a')].\]
Finally, we denote
\begin{align}
    &\Gamma_i(\pi,z,O)\triangleq z^{\top}(R_i(O)\hspace{-2pt}+\hspace{-2pt}H(O)\omega_i^{\star}(\pi))\hspace{-2pt}+\hspace{-2pt}z^{\top}(H(O)\hspace{-2pt}-\hspace{-2pt}\widebar{H}^{\pi})z,\notag\\
    &p_i^{k,t}\triangleq\omega_i^{k,t+1}-\widehat{\omega}_i^{k,t+1}
    \label{eq:def_Gamma_p}
\end{align}

We introduce a few technical lemmas to support the analysis.

\begin{lem}\label{lem:p_bound}
For all $k\geq 0$, we have
\begin{gather*}
    \|p_i^{k,t}\|\leq2(1+2B_{\omega})\beta,\notag\\
    \|z_i^{k,t+1}-z_i^{k,t}\|=\|\omega_i^{k,t+1}-\omega_i^{k,t}\|\leq(1+2B_{\omega})\beta.
\end{gather*}
\end{lem}

\begin{lem}\label{lem:A_negdef}
The matrix $\widebar{H}^{\widehat{\pi}_i^{k}}$ is negative definite
\begin{align*}
    \omega^{\top}\widebar{H}^{\widehat{\pi}_i^{k}}\omega\leq -\frac{(1-\gamma)\underline{\mu}\sigma_{\min}(\Phi)\epsilon}{|\Acal|}\|\omega\|^2, \quad\forall \omega\in\mathbb{R}^{d},k\geq 0,
\end{align*}
where $\underline{\mu}=\min_{\pi,s}\mu_{\pi}(s)$ is a positive constant due to the uniform ergodicity of the Markov chain under any policy.
\end{lem}

\begin{lem}\label{lem:bound_Gamma}
Under Assumption \ref{assump:markov-chain}, we have for all $t\geq\tau$ and $i=0,\cdots,N$
\begin{align*}
    \mathbb{E}[\Gamma_i(\widehat{\pi}_i^k,z_i^{k,t},O_i^{k,t})]\leq2(1+18B_{\omega})^2\beta\tau.
\end{align*}
\end{lem}

With the definition of $R_i$ and $H$ in Eqs.~\eqref{eq:def_vectorR} and \eqref{eq:def_matrixA}, the critic update in Eq.~\eqref{Alg:MT-PDNAC_nestedloop:critic_update} can be re-expressed as
\begin{align*}
    \omega_i^{k,t+1}=\omega_i^{k,t}+\beta(R_i(O_i^{k,t})+H(O_i^{k,t})\omega_i^{k,t})+p_i^{k,t},
\end{align*}
which implies
\begin{align}
    z_i^{k,t+1}-z_i^{k,t}&=\omega_i^{k,t+1}-\omega_i^{k,t}\notag\\
    &=\beta(R_i(O_i^{k,t})+H(O_i^{k,t})\omega_i^{k,t})+p_i^{k,t}\notag\\
    &=\beta(R_i(O_i^{k,t})+H(O_i^{k,t})\omega_i^{\star}(\widehat{\pi}_i^k) + H(O_i^{k,t})z_i^{k,t}) + p_i^{k,t}.\label{prop:conv_critic:eq0}
\end{align}

Then, straightforward manipulations yield
\begin{align}
    \|z_i^{k,t+1}\|^2-\|z_i^{k,t}\|^2&=2(z_i^{k,t})^{\top}(z_i^{k,t+1}-z_i^{k,t})+\|z_i^{k,t+1}-z_i^{k,t}\|^2\notag\\
    &=2(z_i^{k,t})^{\top}(z_i^{k,t+1}-z_i^{k,t}-\beta\widebar{H}^{\widehat{\pi}_i^k}z_i^{k,t})+\|z_i^{k,t+1}-z_i^{k,t}\|^2+2\beta (z_i^{k,t})^{\top}\widebar{H}^{\widehat{\pi}_i^k}z_i^{k,t}\notag\\
    &=2\beta (z_i^{k,t})^{\top}(R_i(O_i^{k,t})+H(O_i^{k,t})(\omega_i^{\star}(\widehat{\pi}_i^k)+z_i^{k,t})-\widebar{H}^{\widehat{\pi}_i^k} z_i^{k,t})\notag\\
    &\hspace{20pt}+2(z_i^{k,t})^{\top}p_i^{k,t}+\|z_i^{k,t+1}-z_i^{k,t}\|^2+2\beta (z_i^{k,t})^{\top}\widebar{H}^{\widehat{\pi}_i^k}z_i^{k,t}\notag\\
    &\leq 2\beta\Gamma_i(\widehat{\pi}_i^k,z_i^{k,t},O_i^{k,t})+2(z_i^{k,t})^{\top}p_i^{k,t}\notag\\
    &\hspace{20pt}+\|z_i^{k,t+1}-z_i^{k,t}\|^2+2\beta (z_i^{k,t})^{\top}\widebar{H}^{\widehat{\pi}_i^k}z_i^{k,t},
    \label{prop:conv_critic:eq1}
\end{align}
where the third equality applies Eq.~\eqref{prop:conv_critic:eq0}.

To bound $(z_i^{k,t})^{\top}p_i^{k,t}$, 
\begin{align}
    (z_i^{k,t})^{\top}p_i^{k,t}&=\langle\omega_i^{k,t}-\omega_i^{\star}(\widehat{\pi}_i^k),p_i^{k,t}\rangle\notag\\
    &\leq\langle\omega_i^{k,t+1}-\omega_i^{\star}(\widehat{\pi}_i^k),\omega_i^{k,t+1}-\widehat{\omega}_i^{k,t+1}\rangle\notag\\
    &\hspace{20pt}+\|z_i^{k,t+1}-z_i^{k,t}\|\|p_i^{k,t}\|\notag\\
    &=\langle\widehat{\omega}_i^{k,t+1} - \omega_i^{\star}(\widehat{\pi}_i^k),\omega_i^{k,t+1}-\widehat{\omega}_i^{k,t+1}\rangle\notag\\
    &\hspace{20pt}+\|\omega_i^{k,t+1}-\widehat{\omega}_i^{k,t+1}\|^2+\|z_i^{k,t+1}-z_i^{k,t}\|\|p_i^{k,t}\|\notag\\
    &\leq 2(1+2B_{\omega})^2\beta^2,\label{prop:conv_critic:eq2}
\end{align}
where the second inequality is due to \cite{doan2019finite}[Lemma 3(a)].

Taking the expectation in Eq.~\eqref{prop:conv_critic:eq1} and applying Eq.~\eqref{prop:conv_critic:eq2} and Lemmas~\ref{lem:A_negdef} and \ref{lem:bound_Gamma}, we have
\begin{align*}
    \mathbb{E}[\|z_i^{k,t+1}\|^2-\|z_i^{k,t}\|^2]&\leq 2\beta\mathbb{E}[\Gamma_i(\widehat{\pi}_i^k,z_i^{k,t},O_i^{k,t})]+2\mathbb{E}[(z_i^{k,t})^{\top}p_i^{k,t}]\notag\\
    &\hspace{20pt}+\mathbb{E}[\|z_i^{k,t+1} - z_i^{k,t}\|^2]+2\beta \mathbb{E}[(z_i^{k,t})^{\top}\widebar{H}^{\widehat{\pi}_i^k}z_i^{k,t}]\notag\\
    &\hspace{-75pt}\leq 4(1+18B_{\omega})^2\beta^2\tau+2(1+2B_{\omega})^2\beta^2+2(1+2B_{\omega})^2\beta^2-\frac{2(1-\gamma)\underline{\mu}\sigma_{\min}(\Phi)\epsilon\beta}{|\Acal|}\mathbb{E}[\|z_i^{k,t}\|^2].
\end{align*}

Re-arranging the terms,
\begin{align}
    \mathbb{E}[\|z_i^{k,t+1}\|^2]&\leq \big(1-\frac{2(1-\gamma)\underline{\mu}\sigma_{\min}(\Phi)\epsilon\beta}{|\Acal|}\big)\mathbb{E}[\|z_i^{k,t}\|^2]+8(1+18B_{\omega})^2\beta^2\tau.\label{prop:conv_critic:eq3}
\end{align}

Recursively applying Eq.~\eqref{prop:conv_critic:eq3}, we get
\begin{align*}
    &\mathbb{E}[\|z_i^{k,T}\|^2]\\
    &\leq\hspace{-2pt} \big(1-\frac{2(1-\gamma)\underline{\mu}\sigma_{\min}(\Phi)\epsilon\beta}{|\Acal|}\big)^{T-\tau}\mathbb{E}[\|z_i^{k,\tau}\|^2]+\hspace{-5pt}\sum_{t=0}^{T-\tau-1}\hspace{-5pt}8(1+18B_{\omega})^2\beta^2\tau\big(1\hspace{-2pt}-\hspace{-2pt}\frac{2(1-\gamma)\underline{\mu}\sigma_{\min}(\Phi)\epsilon\beta}{|\Acal|}\big)^t\notag\\
    &\leq \hspace{-2pt}\big(1\hspace{-2pt}-\hspace{-2pt}\frac{2(1-\gamma)\underline{\mu}\sigma_{\min}(\Phi)\epsilon\beta}{|\Acal|}\big)^{T-\tau}\mathbb{E}[\|z_i^{k,\tau}\|^2]+8(1+18B_{\omega})^2\beta^2\tau\sum_{t=0}^{\infty}\big(1-\frac{2(1-\gamma)\underline{\mu}\sigma_{\min}(\Phi)\epsilon\beta}{|\Acal|}\big)^t\notag\\
    &=4B_{\omega}^2\left(1-\frac{2(1-\gamma)\underline{\mu}\sigma_{\min}(\Phi)\epsilon\beta}{|\Acal|}\right)^{T-\tau}+\frac{4(1+18B_{\omega})^2|\Acal|\beta\tau}{(1-\gamma)\underline{\mu}\sigma_{\min}(\Phi)\epsilon}.
\end{align*}

\qed

\section{Proof of Technical Lemmas}

\subsection{Proof of Lemma~\ref{lem:bounded_lambdastar}}

We skip the proof as it is a simple extension of \cite{ding2020natural}[Lemma 1].

\subsection{Proof of Lemma~\ref{lem:consensus_error_pg}}

We denote $g_i^k=(\frac{1}{N}+\lambda_i^k-\nu_i^k)Q_i^{\pi_i^k}\in\mathbb{R}^{|\Scal||\Acal|}$ and $g^k=[(g_1^k)^{\top},\dots,(g_N^k)^{\top}]^{\top}\in\mathbb{R}^{N|\Scal||\Acal|}$.
It is easy to see
\begin{align*}
    \|g_i^k\|\leq \left|\frac{1}{N}+\lambda_i^k-\nu_i^k\right|\|Q_i^{\pi_i^k}\|\leq\frac{(B_{\lambda}+\frac{1}{N})\sqrt{|\Scal||\Acal|}}{1-\gamma},
\end{align*}
which implies $\|g^k\|\leq\frac{(B_{\lambda}+\frac{1}{N})\sqrt{N|\Scal||\Acal|}}{1-\gamma}$ for all $k$.
Then, using an argument similar to the one in \cite{yuan2016convergence}[Lemma 1], we can get
\begin{align}
    \|\widebar{\theta}^k-\theta_i^k\|\leq\frac{(B_{\lambda}+\frac{1}{N})\sqrt{N|\Scal||\Acal|}\alpha}{(1-\gamma)(1-\sigma_2(W))}.\label{lem:consensus_error:eq1}
\end{align}

The softmax function is Lipschitz with constant $1$, i.e.
\begin{align*}
    \|\pi_{\theta}-\pi_{\theta'}\|\leq\|\theta-\theta'\|,\quad\forall \theta,\theta',
\end{align*}

Recall the definition of $\widebar{\pi}^k$ in Eq.~\eqref{eq:def_aggregate_notation}.
The Lipschitz continuity and Eq.~\eqref{lem:consensus_error:eq1} imply the claimed result.

\qed

\subsection{Proof of Lemma~\ref{lem:consensus_error_ac}}

We skip the proof and note that it is almost identical to the proof of Lemma~\ref{lem:consensus_error_pg}.

\subsection{Proof of Lemma~\ref{lem:omega_bound}}

Due to the boundedness of the reward function, it is easy to see that $|Q_i^{\pi}(s,a)|\leq\frac{1}{1-\gamma}$, which implies
\[\|Q_i^{\pi}\|\leq\sqrt{\frac{|\Scal||\Acal|}{1-\gamma}}.\]

Since $\|\Phi\omega_i^{\star}(\widehat{\pi}_i^k)-Q_{i}^{\widehat{\pi}_i^k}\|\leq\varepsilon_{\max}$ due to Assumption~\ref{assump:epsilon_max}, we have
\begin{align}
    \|\Phi\omega_i^{\star}(\widehat{\pi}_i^k)\|\leq\|Q_{i}^{\widehat{\pi}_i^k}\|+\varepsilon_{\max}\leq\sqrt{\frac{|\Scal||\Acal|}{1-\gamma}}+\varepsilon_{\max},
\end{align}
which implies
\begin{align}
    \|\omega_i^{\star}(\widehat{\pi}_i^k)\|\leq\sigma_{\min}^{-1}(\Phi)\Big(\sqrt{\frac{|\Scal||\Acal|}{1-\gamma}}+\varepsilon_{\max}\Big)=B_{\omega}.
\end{align}

To show the bound on $\widehat{Q}_{i}^k$, note that $\|\omega_i^{k,T}\|\leq B_{\omega}$ due to the projection in Eq.~\eqref{Alg:MT-PDNAC_nestedloop:critic_update}. As a result,
\[\|\widehat{Q}_{i}^k\|=\|\Phi\omega_i^{k,T}\|\leq \sigma_{\max}(\Phi) B_{\omega}.\]
The bound on $\|\widehat{V}_{i,t}\|$ easily follows from Jensen's inequality
\begin{align*}
\|\widehat{V}_{i}^k\|^2&=\sum_{s}\left(\sum_{a} \pi_i^k(a\mid s) \widehat{Q}_{i}^k(s,a)\right)^2\leq \sum_{s}\sum_{a}\pi_i^k(a\mid s)\left(\widehat{Q}_{i}^k(s,a)\right)^2\leq \sum_{s,a}\left(\widehat{Q}_{i}^k(s,a)\right)^2=\|\widehat{Q}_{i}^k\|^2.
\end{align*}

\qed

\subsection{Proof of Lemma~\ref{lem:consensus_error_LFA}}

We denote 
\[g_i^k=(\frac{1}{N}+\lambda_i^k-\nu_i^k)\omega_i^{k,T}\in\mathbb{R}^{d},\quad\text{and}\quad g^k=[(g_1^k)^{\top},\dots,(g_N^k)^{\top}]^{\top}\in\mathbb{R}^{Nd}.\]
Then,
\begin{align*}
    \|g_i^k\|\leq \left|\frac{1}{N}+\lambda_i^k-\nu_i^k\right|\|\omega_i^{k,T}\|\leq (B_{\lambda}+\frac{1}{N})B_{\omega},
\end{align*}
which implies $\|g^k\|\leq(B_{\lambda}+\frac{1}{N})B_{\omega}\sqrt{N}$ for all $k$.
Then, using an argument similar to the one in \cite{yuan2016convergence}[Lemma 1], we can get
\begin{align}
    \|\widebar{\theta}^k-\theta_i^k\|\leq\frac{(B_{\lambda}+\frac{1}{N})B_{\omega}\sqrt{N}\alpha}{1-\sigma_2(W)}.\label{lem:consensus_error_LFA:eq1}
\end{align}

The softmax function is Lipschitz continuous with constant $1$, which implies
\begin{align*}
    \|\widebar{\pi}^k-\pi_i^k\|\leq\|\Phi(\widebar{\theta}^k-\theta_i^k)\|\leq\sigma_{\max}(\Phi)\|\widebar{\theta}^k-\theta_i^k\|\leq\Ocal\left(\frac{\sqrt{N}\alpha}{1-\sigma_2(W)}\right).
\end{align*}

\qed

\subsection{Proof of Lemma \ref{lem:step_improvement}}

The performance difference lemma states that for any policies $\pi_1$, $\pi_2$, initial distribution $\zeta$, and $i=0,\cdots,N$
\begin{align}
    V_i^{\pi_1}(\zeta)-V_i^{\pi_2}(\zeta)=\frac{1}{1-\gamma}\mathbb{E}_{s \sim d_{\zeta}^{\pi^{\star}}, a \sim \pi^{\star}(\cdot \mid s)}\left[A_{0}^{\pi_k}(s, a)\right].
    \label{eq:def_lemma_performancediff}
\end{align}

By this lemma,
\begin{align*}
&V_{0}^{\widebar{\pi}^{k+1}}(\zeta)-V_{0}^{\widebar{\pi}^{k}}(\zeta)\notag\\
&= \frac{1}{1-\gamma}\mathbb{E}_{s \sim d_{\zeta}^{\widebar{\pi}^{k+1}}, a \sim \widebar{\pi}^{k+1}(\cdot \mid s)}\left[A_{0}^{\widebar{\pi}^k}(s, a)\right] \notag\\
&= \frac{1}{1-\gamma}\mathbb{E}_{s \sim d_{\zeta}^{\widebar{\pi}^{k+1}}, a \sim \widebar{\pi}^{k+1}(\cdot \mid s)}\left[Q_{0}^{\widebar{\pi}^k}(s, a)\right] -\frac{1}{1-\gamma}\mathbb{E}_{s \sim d_{\zeta}^{\widebar{\pi}^{k+1}}}\left[V_{0}^{\widebar{\pi}^k}(s)\right] \notag\\
&=\frac{1}{1-\gamma}\mathbb{E}_{s \sim d_{\zeta}^{\widebar{\pi}^{k+1}}, a \sim \widebar{\pi}^{k+1}(\cdot \mid s)}\left[Q_{L,k}^{k}(s, a)\right]+\frac{1}{1-\gamma}\mathbb{E}_{s \sim d_{\zeta}^{\widebar{\pi}^{k+1}}, a \sim \widebar{\pi}^{k+1}(\cdot \mid s)}\left[Q_{L,k}^{\vpi^k}(s, a)-Q_{L,k}^{k}(s, a)\right]\notag\\
&\hspace{20pt}+\frac{1}{1-\gamma}\mathbb{E}_{s \sim d_{\zeta}^{\widebar{\pi}^{k+1}}, a \sim \widebar{\pi}^{k+1}(\cdot \mid s)}\left[Q_{L,k}^{\widebar{\pi}^k}(s, a)-Q_{L,k}^{\vpi^k}(s, a)\right]\notag\\
&\hspace{20pt}-\frac{(\lambda^k-\nu^k)^{\top}}{1-\gamma}\mathbb{E}_{s \sim d_{\zeta}^{\widebar{\pi}^{k+1}}, a \sim \widebar{\pi}^{k+1}(\cdot \mid s)}\left[Q_{g}^{\widebar{\pi}^k}(s, a)\right]-\frac{1}{1-\gamma}\mathbb{E}_{s \sim d_{\zeta}^{\widebar{\pi}^{k+1}}}\left[V_{0}^{\widebar{\pi}^k}(s)\right].
\end{align*}

Note that the actor update rule in Eq.~\eqref{prop:actor_conv:updates_rewrite} implies
\begin{align*}
    Q_{L,k}^{k}(s, a)=\frac{N}{\alpha}\log\left(\frac{\widebar{\pi}^{k+1}(a\mid s)}{\widebar{\pi}^k(a\mid s)}Z_k(s)\right).
\end{align*}

Combining the two equalities above, we have
\begin{align*}
    &V_{0}^{\widebar{\pi}^{k+1}}(\zeta)-V_{0}^{\widebar{\pi}^{k}}(\zeta)\notag\\
    &=\frac{N}{\alpha(1-\gamma)}\mathbb{E}_{s \sim d_{\zeta}^{\widebar{\pi}^{k+1}}\hspace{-2pt}, a \sim \widebar{\pi}^{k+1}(\cdot \mid s)}\left[\log\left(\frac{\widebar{\pi}^{k+1}(a\mid s)}{\widebar{\pi}^k(a\mid s)}Z_k(s)\right)\right]\\
    &\hspace{20pt}+\frac{1}{1-\gamma}\mathbb{E}_{s \sim d_{\zeta}^{\widebar{\pi}^{k+1}}, a \sim \widebar{\pi}^{k+1}(\cdot \mid s)} \left[Q_{L,k}^{\vpi^k}(s, a)-Q_{L,k}^{k}(s, a)\right]\notag\\
    &\hspace{20pt}+\frac{1}{1-\gamma}\mathbb{E}_{s \sim d_{\zeta}^{\widebar{\pi}^{k+1}}, a \sim \widebar{\pi}^{k+1}(\cdot \mid s)}\left[Q_{L,k}^{\widebar{\pi}^k}(s, a)-Q_{L,k}^{\vpi^k}(s, a)\right]\notag\\
    &\hspace{20pt}-\frac{(\lambda^k-\nu^k)^{\top}}{1-\gamma}\mathbb{E}_{s \sim d_{\zeta}^{\widebar{\pi}^{k+1}}, a \sim \widebar{\pi}^{k+1}(\cdot \mid s)}\left[Q_{g}^{\widebar{\pi}^k}(s, a)\right]-\frac{1}{1-\gamma}\mathbb{E}_{s \sim d_{\zeta}^{\widebar{\pi}^{k+1}}}\left[V_{0}^{\widebar{\pi}^k}(s)\right]\notag\\
    &\geq\frac{N}{\alpha(1-\gamma)} \mathbb{E}_{s \sim d_{\zeta}^{\widebar{\pi}^{k+1}}}\left[D_{KL}(\widebar{\pi}^{k+1}(\cdot\mid s)||\widebar{\pi}^{k}(\cdot\mid s))\right]+\frac{N}{\alpha(1-\gamma)}\mathbb{E}_{s \sim d_{\zeta}^{\widebar{\pi}^{k+1}}}\left[\log Z_k(s)\right]\notag\\
    &\hspace{20pt}-\frac{B_{\lambda}+1/N}{1-\gamma}\sum_{i=1}^{N}\|Q_i^{\pi_i^k}-Q_i^k\|-\frac{\sqrt{|\Scal|| \Acal|}(B_{\lambda}+1/N)}{(1-\gamma)^3} \sum_{i=1}^{N}\|\widebar{\pi}^k-\pi_i^k\|\notag\\
    &\hspace{20pt}-\frac{(\lambda^k-\nu^k)^{\top}}{1-\gamma}\mathbb{E}_{s \sim d_{\zeta}^{\widebar{\pi}^{k+1}}\hspace{-2pt}, a \sim \widebar{\pi}^{k+1}(\cdot \mid s)}\left[A_{g}^{\widebar{\pi}^k}(s, a)\right]-\frac{1}{1-\gamma}\mathbb{E}_{s \sim d_{\zeta}^{\widebar{\pi}^{k+1}}}\left[V_{L,k}^{\widebar{\pi}^k}(s)\right]\notag\\
    &\geq\frac{N}{\alpha(1\hspace{-2pt}-\hspace{-2pt}\gamma)} \mathbb{E}_{s \sim d_{\zeta}^{\widebar{\pi}^{k+1}}}\left[\log Z_k(s)\right]\hspace{-2pt}-\hspace{-2pt}\frac{B_{\lambda}\hspace{-2pt}+\hspace{-2pt}1/N}{1-\gamma}\sum_{i=1}^{N}\|Q_i^{\pi_i^k}\hspace{-2pt}-\hspace{-2pt}Q_i^k\|\hspace{-2pt}-\hspace{-2pt}\frac{\sqrt{|\Scal|| \Acal|}(B_{\lambda}+1/N)}{(1-\gamma)^3}\sum_{i=1}^{N}\|\widebar{\pi}^k-\pi_i^k\|\notag\\
    &\hspace{20pt}-(\lambda^k-\nu^k)^{\top}\left(V_{g}^{\widebar{\pi}^{k+1}}(\zeta)-V_{g}^{\widebar{\pi}^k}(\zeta)\right)-\frac{1}{1-\gamma}\mathbb{E}_{s \sim d_{\zeta}^{\widebar{\pi}^{k+1}}}\left[V_{L,k}^{\widebar{\pi}^k}(s)\right],
\end{align*}
where the last inequality applies the performance difference lemma. Rearranging this inequality leads to
\begin{align}
    V_{L,k}^{\widebar{\pi}^{k+1}}(\zeta)-V_{L,k}^{\widebar{\pi}^{k}}(\zeta)&\geq\frac{N}{\alpha(1-\gamma)}\mathbb{E}_{s \sim d_{\zeta}^{\widebar{\pi}^{k+1}}}\left[\log Z_k(s)\right]-\frac{\sqrt{|\Scal|| \Acal|}(B_{\lambda}+1/N)}{(1-\gamma)^3}\sum_{i=1}^{N}\|\widebar{\pi}^k-\pi_i^k\| \notag\\
    &\hspace{20pt}-\frac{B_{\lambda}+1/N}{1-\gamma}\sum_{i=1}^{N}\|Q_i^{\pi_i^k}-Q_i^k\|- \frac{1}{1-\gamma}\mathbb{E}_{s \sim d_{\zeta}^{\widebar{\pi}^{k+1}}}\left[V_{L,k}^{\widebar{\pi}^k}(s)\right].\label{lem:step_improvement:eq1}
\end{align}

From the definition of $Z^k$ and Jensen's inequality,
\begin{align*}
    \log Z^k(s)&=\log\left(\sum_{a'\in\Acal}\widebar{\pi}^{k}(a'\mid s)\exp(\frac{\alpha}{N}\sum_{i=1}^{N}(\frac{1}{N}+\lambda_i^k-\nu_i^k)Q_i^{k}(s,a'))\right)\notag\\
    &\geq \sum_{a'\in\Acal}\widebar{\pi}^{k}(a'\mid s)\log\left(\exp(\frac{\alpha}{N}\sum_{i=1}^{N}(\frac{1}{N}+\lambda_i^k-\nu_i^k)Q_i^{k}(s,a'))\right)\notag\\
    &= \frac{\alpha}{N}\sum_{a'\in\Acal}\widebar{\pi}^{k}(a'\mid s)\sum_{i=1}^{N}(\frac{1}{N}+\lambda_i^k-\nu_i^k)Q_i^{k}(s,a')\notag\\
    &= \frac{\alpha}{N}\sum_{a'\in\Acal}\widebar{\pi}^{k}(a'\mid s)\sum_{i=1}^{N}(\frac{1}{N}+\lambda_i^k-\nu_i^k)Q_i^{\widebar{\pi}^k}(s,a')\notag\\
    &\hspace{20pt}+\frac{\alpha}{N}\sum_{a'\in\Acal}\widebar{\pi}^{k}(a'\mid s)\sum_{i=1}^{N}(\frac{1}{N}+\lambda_i^k-\nu_i^k)(Q_i^{k}(s,a')-Q_i^{\pi_i^k}(s,a'))\notag\\
    &\hspace{20pt}+\frac{\alpha}{N}\sum_{a'\in\Acal}\widebar{\pi}^{k}(a'\mid s)\sum_{i=1}^{N}(\frac{1}{N}+\lambda_i^k-\nu_i^k)(Q_i^{\pi_i^k}(s,a')-Q_i^{\widebar{\pi}^k}(s,a'))\notag\\
    &\geq \frac{\alpha}{N}V_{L,k}^{\widebar{\pi}^k}(s)-\frac{(B_{\lambda}+1/N)\alpha}{N}\sum_{i=1}^{N}\|Q_i^{\pi_i^k}-Q_i^k\|-\frac{\sqrt{|\Scal|| \Acal|}(B_{\lambda}+1/N)\alpha}{N(1-\gamma)^2}\sum_{i=1}^{N}\|\widebar{\pi}^k-\pi_i^k\|.
\end{align*}

This bound on $\log Z_k(s)$ implies
\begin{align*}
    &\frac{N}{\alpha(1-\gamma)}\mathbb{E}_{s \sim d_{\zeta}^{\widebar{\pi}^{k+1}}}\left[\log Z_k(s)\right]-\frac{1}{1-\gamma}\mathbb{E}_{s \sim d_{\zeta}^{\widebar{\pi}^{k+1}}}\left[V_{L,k}^{\widebar{\pi}^k}(s)\right]\notag\\
    &=\frac{N}{\alpha(1-\gamma)}\sum_{s}d_{\zeta}^{\widebar{\pi}^{k+1}}(s)\left(\log Z_k(s)-\frac{\alpha}{N}V_{L,k}^{\widebar{\pi}^k}(s)\right)\notag\\
    &=\frac{N}{\alpha(1-\gamma)}\sum_{s}d_{\zeta}^{\widebar{\pi}^{k+1}}(s)\Big(\log Z_k(s)-\frac{\alpha}{N}V_{L,k}^{\widebar{\pi}^k}(s)+\frac{(B_{\lambda}+1/N)\alpha}{N}\sum_{i=1}^{N}\|Q_i^{\pi_i^k}-Q_i^k\|\notag\\
    &\hspace{120pt}+\frac{\sqrt{|\Scal|| \Acal|}(B_{\lambda}+1/N)\alpha}{N(1-\gamma)^2}\sum_{i=1}^{N}\|\widebar{\pi}^k-\pi_i^k\|\Big)\notag\\
    &\hspace{20pt}-\frac{B_{\lambda}+1/N}{1-\gamma}\sum_{i=1}^{N}\|Q_i^{\pi_i^k}-Q_i^k\|-\frac{\sqrt{|\Scal|| \Acal|}(B_{\lambda}+1/N)}{(1-\gamma)^3}\sum_{i=1}^{N}\|\widebar{\pi}^k-\pi_i^k\|\notag\\
    &\geq \frac{N}{\alpha}\sum_{s}\zeta(s)\left(\log Z_k(s)-\frac{\alpha}{N}V_{L,k}^{\widebar{\pi}^k}(s)\right)-\frac{B_{\lambda}+1/N}{1-\gamma}\sum_{i=1}^{N}\|Q_i^{\pi_i^k}-Q_i^k\| \notag\\
    &\hspace{20pt}-\frac{\sqrt{|\Scal|| \Acal|}(B_{\lambda}+1/N)}{(1-\gamma)^3}\sum_{i=1}^{N}\|\widebar{\pi}^k-\pi_i^k\|,
\end{align*}
where the inequality follows from the fact that $d_{\zeta}^{\pi} \geq(1-\gamma) \zeta$ elementwise for any policy $\pi$. Plugging this bound into Eq.~\eqref{lem:step_improvement:eq1}, we have
\begin{align*}
    V_{L,k}^{\widebar{\pi}^{k+1}}(\zeta)-V_{L,k}^{\widebar{\pi}^{k}}(\zeta)
    &\geq \frac{N}{\alpha}\mathbb{E}_{s \sim \zeta}\left[\log Z_k(s)-\frac{\alpha}{N}V_{L,k}^{\widebar{\pi}^k}(s)\right]-\frac{2(B_{\lambda}+1/N)}{1-\gamma}\sum_{i=1}^{N}\|Q_i^{\pi_i^k}-Q_i^k\|\notag\\
    &\hspace{20pt}-\frac{2\sqrt{|\Scal|| \Acal|}(B_{\lambda}+1/N)}{(1-\gamma)^3}\sum_{i=1}^{N}\|\widebar{\pi}^k-\pi_i^k\|.
\end{align*}

\qed

\subsection{Proof of Lemma~\ref{lem:bounded_Lagrangian}}

By the performance difference lemma in Eq.~\eqref{eq:def_lemma_performancediff},
\begin{align*}
V_{0}^{\pi^{\star}}(\rho)-V_0^{\widebar{\pi}^k}(\rho)
&=\frac{1}{N}\sum_{i=1}^{N}(V_{i}^{\pi^{\star}}(\rho)-V_i^{\widebar{\pi}^k}(\rho))\notag\\
&=\frac{1}{N}\sum_{i=1}^{N}(V_{i}^{\pi^{\star}}(\rho)-V_i^{\pi_i^k}(\rho))+\frac{1}{N}\sum_{i=1}^{N}(V_{i}^{\pi_i^k}(\rho)-V_i^{\widebar{\pi}^k}(\rho))\notag\\
&\leq \frac{1}{N(1-\gamma)}\sum_{i=1}^{N}\mathbb{E}_{s \sim d_{\rho}^{\pi^{\star}}, a \sim \pi^{\star}(\cdot \mid s)}\left[A_{i}^{\pi_i^k}(s, a)\right] + \frac{\sqrt{|\Scal||\Acal|}}{N(1-\gamma)^3}\sum_{i=1}^{N}\|\widebar{\pi}^k-\pi_i^k\|\notag\\
&=\frac{1}{N(1-\gamma)}\sum_{i=1}^{N}\mathbb{E}_{s \sim d_{\rho}^{\pi^{\star}},a\sim\pi^{\star}(\cdot\mid s)}\left[ Q_{i}^{\pi_i^k}(s, a)\right]-\frac{1}{N(1-\gamma)}\sum_{i=1}^{N}\mathbb{E}_{s \sim d_{\rho}^{\pi^{\star}}}\left[V_{i}^{\pi_i^k}(s)\right]\notag\\
&\hspace{20pt}+ \frac{\sqrt{|\Scal||\Acal|}}{N(1-\gamma)^3}\sum_{i=1}^{N}\|\widebar{\pi}^k-\pi_i^k\|.
\end{align*}

Plugging in the update rule of the policy,
\begin{align*}
&V_{0}^{\pi^{\star}}(\rho)-V_0^{\widebar{\pi}^k}(\rho)\notag\\
&\leq\frac{1}{N(1-\gamma)}\sum_{i=1}^{N}\mathbb{E}_{s \sim d_{\rho}^{\pi^{\star}},a\sim\pi^{\star}(\cdot\mid s)}\left[ Q_{i}^{\pi_i^k}(s, a)\right]-\frac{1}{N(1-\gamma)}\sum_{i=1}^{N}\mathbb{E}_{s \sim d_{\rho}^{\pi^{\star}}}\left[V_{i}^{\pi_i^k}(s)\right]\notag\\
&\hspace{20pt}+ \frac{\sqrt{|\Scal||\Acal|}}{N(1-\gamma)^3}\sum_{i=1}^{N}\|\widebar{\pi}^k-\pi_i^k\|\notag\\
&=\frac{1}{1\hspace{-2pt}-\hspace{-2pt}\gamma}\mathbb{E}_{s \sim d_{\rho}^{\pi^{\star}}, a \sim \pi^{\star}(\cdot \mid s)}\left[\sum_{i=1}^{N}(\frac{1}{N}\hspace{-2pt}+\hspace{-2pt}\lambda_i^k\hspace{-2pt}-\hspace{-2pt}\nu_i^k)Q_i^{\pi_i^k}(s, a)\right]-\frac{(\lambda^k\hspace{-2pt}-\hspace{-2pt}\nu^k)^{\top}}{1-\gamma}\mathbb{E}_{s \sim d_{\rho}^{\pi^{\star}}, a \sim \pi^{\star}(\cdot \mid s)}\left[Q_{g}^{\vpi^k}(s, a)\right]\notag\\
&\hspace{20pt}-\frac{1}{N(1-\gamma)}\sum_{i=1}^{N}\mathbb{E}_{s \sim d_{\rho}^{\pi^{\star}}}\left[V_{i}^{\pi_i^k}(s)\right]+ \frac{\sqrt{|\Scal||\Acal|}}{N(1-\gamma)^3}\sum_{i=1}^{N}\|\widebar{\pi}^k-\pi_i^k\| \notag\\
&=\frac{N}{\alpha(1-\gamma)}\mathbb{E}_{s \sim d_{\rho}^{\pi^{\star}}, a \sim \pi^{\star}(\cdot \mid s)}\left[\log\left(\frac{\widebar{\pi}^{k+1}(a\mid s)}{\widebar{\pi}^k(a\mid s)}Z_k(s)\right)\right]\notag\\
&\hspace{20pt}+\frac{1}{1-\gamma}\mathbb{E}_{s \sim d_{\rho}^{\pi^{\star}}, a \sim \pi^{\star}(\cdot \mid s)}\left[\sum_{i=1}^{N}(\frac{1}{N}+\lambda_i^k-\nu_i^k)(Q_i^{\pi_i^k}(s, a)-Q_i^{k}(s, a))\right]\notag\\
&\hspace{20pt}-\frac{(\lambda^k-\nu_k)^{\top}}{1-\gamma}\mathbb{E}_{s \sim d_{\rho}^{\pi^{\star}}, a \sim \pi^{\star}(\cdot \mid s)}\left[A_{g}^{\vpi^k}(s, a)\right]-\frac{(\lambda^k-\nu_k)^{\top}}{1-\gamma}\mathbb{E}_{s \sim d_{\rho}^{\pi^{\star}}}\left[V_{g}^{\vpi^k}(s)\right]\notag\\
&\hspace{20pt}-\frac{1}{N(1-\gamma)}\sum_{i=1}^{N}\mathbb{E}_{s \sim d_{\rho}^{\pi^{\star}}}\left[V_{i}^{\pi_i^k}(s)\right]+ \frac{\sqrt{|\Scal||\Acal|}}{N(1-\gamma)^3}\sum_{i=1}^{N}\|\widebar{\pi}^k-\pi_i^k\| \notag\\
&\leq\frac{N}{\alpha(1-\gamma)}\mathbb{E}_{s \sim d_{\rho}^{\pi^{\star}}}\left[D_{\text{KL}}(\pi^{\star}(\cdot\mid s)||\widebar{\pi}^k(\cdot\mid s))\right]-\frac{N}{\alpha(1-\gamma)}\mathbb{E}_{s \sim d_{\rho}^{\pi^{\star}}}\left[D_{\text{KL}}(\pi^{\star}(\cdot\mid s)||\widebar{\pi}^{k+1}(\cdot\mid s))\right]\notag\\
&\hspace{20pt}+\frac{N}{\alpha(1-\gamma)}\mathbb{E}_{s \sim d_{\rho}^{\pi^{\star}}}\left[\log Z^k(s)\right]-\frac{(\lambda^k-\nu^k)^{\top}}{1-\gamma}\mathbb{E}_{s \sim d_{\rho}^{\pi^{\star}}, a \sim \pi^{\star}(\cdot \mid s)}\left[A_{g}^{\vpi^k}(s, a)\right]\notag\\
&\hspace{20pt}-\frac{(\lambda^k-\nu^k)^{\top}}{1-\gamma}\mathbb{E}_{s \sim d_{\rho}^{\pi^{\star}}}\left[V_{g}^{\vpi^k}(s)\right]-\frac{1}{N(1-\gamma)}\sum_{i=1}^{N}\mathbb{E}_{s \sim d_{\rho}^{\pi^{\star}}}\left[V_{i}^{\pi_i^k}(s)\right]\notag\\
&\hspace{20pt}+\frac{B_{\lambda}+1/N}{1-\gamma}\sum_{i=1}^{N}\|Q_i^k-Q_i^{\pi_i^k}\|+ \frac{\sqrt{|\Scal||\Acal|}}{N(1-\gamma)^3}\sum_{i=1}^{N}\|\widebar{\pi}^k-\pi_i^k\|.
\end{align*}

Re-grouping the terms,
\begin{align}
&V_{0}^{\pi^{\star}}(\rho)-V_0^{\widebar{\pi}^k}(\rho)\notag\\
&\leq\frac{N}{\alpha(1-\gamma)}\mathbb{E}_{s \sim d_{\rho}^{\pi^{\star}}}\left[D_{\text{KL}}(\pi^{\star}(\cdot\mid s)||\widebar{\pi}^k(\cdot\mid s))\right]\hspace{-2pt}-\hspace{-2pt}\frac{N}{\alpha(1-\gamma)}\mathbb{E}_{s \sim d_{\rho}^{\pi^{\star}}}\left[D_{\text{KL}}(\pi^{\star}(\cdot\mid s)||\widebar{\pi}^{k+1}(\cdot\mid s))\right]\notag\\
&\hspace{20pt}+\frac{N}{\alpha(1-\gamma)}\mathbb{E}_{s \sim d_{\rho}^{\pi^{\star}}}\left[\log Z^k(s)\right]-\frac{(\lambda^k-\nu_k)^{\top}}{1-\gamma}\mathbb{E}_{s \sim d_{\rho}^{\pi^{\star}}, a \sim \pi^{\star}(\cdot \mid s)}\left[A_{g}^{\widebar{\pi}^k}(s, a)\right]\notag\\
&\hspace{20pt}+\frac{(\lambda^k-\nu^k)^{\top}}{1-\gamma}\mathbb{E}_{s \sim d_{\rho}^{\pi^{\star}}, a \sim \pi^{\star}(\cdot \mid s)}\left[A_{g}^{\widebar{\pi}^k}(s, a)-A_{g}^{\vpi^k}(s, a)\right] -\frac{1}{1-\gamma}\mathbb{E}_{s \sim d_{\rho}^{\pi^{\star}}}\left[V_{L,k}^{\vpi^k}(s)\right]\notag\\
&\hspace{20pt}+\frac{B_{\lambda}+1/N}{1-\gamma}\sum_{i=1}^{N}\|Q_i^k-Q_i^{\pi_i^k}\| + \frac{\sqrt{|\Scal||\Acal|}}{N(1-\gamma)^3}\sum_{i=1}^{N}\|\widebar{\pi}^k-\pi_i^k\|\notag\\
&\leq\frac{N}{\alpha(1-\gamma)}\mathbb{E}_{s \sim d_{\rho}^{\pi^{\star}}}\left[D_{\text{KL}}(\pi^{\star}(\cdot\mid s)||\widebar{\pi}^k(\cdot\mid s))\right]\hspace{-2pt}-\hspace{-2pt}\frac{N}{\alpha(1-\gamma)}\mathbb{E}_{s \sim d_{\rho}^{\pi^{\star}}}\left[D_{\text{KL}}(\pi^{\star}(\cdot\mid s)||\widebar{\pi}^{k+1}(\cdot\mid s))\right]\label{lem:bounded_Lagrangian:eq1}\\
&\hspace{20pt}+\frac{N}{\alpha(1-\gamma)}\mathbb{E}_{s \sim d_{\rho}^{\pi^{\star}}}\left[\log Z^k(s)\right]\hspace{-2pt}-\hspace{-2pt}\frac{(\lambda^k\hspace{-2pt}-\hspace{-2pt}\nu_k)^{\top}}{1-\gamma}\left(V_g^{\pi^{\star}}(s)-V_g^{\widebar{\pi}^k}(s)\right)\hspace{-2pt}-\hspace{-2pt}\frac{1}{1-\gamma}\mathbb{E}_{s \sim d_{\rho}^{\pi^{\star}}}\left[V_{L,k}^{\vpi^k}(s)\right] \notag\\
&\hspace{20pt}+\frac{B_{\lambda}+1/N}{1-\gamma}\sum_{i=1}^{N}\|Q_i^k-Q_i^{\pi_i^k}\|+ \frac{\sqrt{|\Scal||\Acal|}(B_{\lambda}+1/N)}{(1-\gamma)^3}\sum_{i=1}^{N}\|\widebar{\pi}^k-\pi_i^k\|,\notag
\end{align}
where the second inequality follows from the performance difference lemma and the Lipschitz continuity of the advantage.

Applying Lemma~\ref{lem:step_improvement} with $\zeta=d_{\rho}^{\pi^{\star}}$,
\begin{align}
    \frac{N}{\alpha}\mathbb{E}_{s \sim d_{\rho}^{\pi^{\star}}}\left[\log Z_k(s)-\frac{\alpha}{N}V_{L,k}^{\widebar{\pi}^k}(s)\right]&\leq V_{L,k}^{\widebar{\pi}^{k+1}}(d_{\rho}^{\pi^{\star}})-V_{L,k}^{\widebar{\pi}^{k}}(d_{\rho}^{\pi^{\star}})+\frac{2(B_{\lambda}+1/N)}{1-\gamma}\sum_{i=1}^{N}\|Q_i^{\pi_i^k}-Q_i^k\|\notag\\
    &\hspace{20pt}+ \frac{2\sqrt{|\Scal|| \Acal|}(B_{\lambda}+1/N)}{(1-\gamma)^3}\sum_{i=1}^{N}\|\widebar{\pi}^k-\pi_i^k\|.\label{lem:bounded_Lagrangian:eq2}
\end{align}

Combining Eqs.~\eqref{lem:bounded_Lagrangian:eq1} and \eqref{lem:bounded_Lagrangian:eq2},
\begin{align*}
&V_{0}^{\pi^{\star}}(\rho)-V_0^{\widebar{\pi}^k}(\rho)\\
&\leq\frac{N}{\alpha(1-\gamma)}\mathbb{E}_{s \sim d_{\rho}^{\pi^{\star}}}\left[D_{\text{KL}}(\pi^{\star}(\cdot\mid s)||\widebar{\pi}^k(\cdot\mid s))\right]\hspace{-2pt}-\hspace{-2pt}\frac{N}{\alpha(1-\gamma)}\mathbb{E}_{s \sim d_{\rho}^{\pi^{\star}}}\left[D_{\text{KL}}(\pi^{\star}(\cdot\mid s)||\widebar{\pi}^{k+1}(\cdot\mid s))\right]\notag\\
&\hspace{20pt}+\frac{1}{1-\gamma}\left(V_{L,k}^{\widebar{\pi}^{k+1}}(d_{\rho}^{\pi^{\star}})-V_{L,k}^{\widebar{\pi}^{k}}(d_{\rho}^{\pi^{\star}})\right)+\frac{2(B_{\lambda}+1/N)}{(1-\gamma)^2}\sum_{i=1}^{N}\|Q_i^{\pi_i^k}-Q_i^k\|\notag\\
&\hspace{20pt}+\frac{2\sqrt{|\Scal|| \Acal|}(B_{\lambda}+1/N)}{(1-\gamma)^4}\sum_{i=1}^{N}\|\widebar{\pi}^k-\pi_i^k\|-\frac{(\lambda^k-\nu^k)^{\top}}{1-\gamma}\left(V_g^{\pi^{\star}}(s)-V_g^{\widebar{\pi}^k}(s)\right)\notag\\
&\hspace{20pt} +\frac{B_{\lambda}+1/N}{1-\gamma}\sum_{i=1}^{N}\|Q_i^{\pi_i^k}-Q_i^k\|+ \frac{\sqrt{|\Scal||\Acal|}(B_{\lambda}+1/N)}{N(1-\gamma)^3}\sum_{i=1}^{N}\|\widebar{\pi}^k-\pi_i^k\|,
\end{align*}
which implies
\begin{align*}
    V_{L,k}^{\pi^{\star}}(\rho)-V_{L,k}^{\widebar{\pi}^k}(\rho)&\leq\frac{N}{\alpha(1-\gamma)}\mathbb{E}_{s \sim d_{\rho}^{\pi^{\star}}}\left[D_{\text{KL}}(\pi^{\star}(\cdot\mid s)||\widebar{\pi}^k(\cdot\mid s))\right]\notag\\
    &\hspace{-25pt}-\frac{N}{\alpha(1-\gamma)}\mathbb{E}_{s \sim d_{\rho}^{\pi^{\star}}}\left[D_{\text{KL}}(\pi^{\star}(\cdot\mid s)||\widebar{\pi}^{k+1}(\cdot\mid s))\right]+\frac{1}{1-\gamma}\left(V_{L,k}^{\widebar{\pi}^{k+1}}(d_{\rho}^{\pi^{\star}})-V_{L,k}^{\widebar{\pi}^{k}}(d_{\rho}^{\pi^{\star}})\right) \notag\\
    &\hspace{20pt} +\frac{3(B_{\lambda}+1/N)}{(1-\gamma)^2}\sum_{i=1}^{N}\|Q_i^{\pi_i^k}-Q_i^k\|+\frac{3\sqrt{|\Scal|| \Acal|}(B_{\lambda}+1/N)}{(1-\gamma)^4}\sum_{i=1}^{N}\|\widebar{\pi}^k-\pi_i^k\|.
\end{align*}

Taking the average from $k=0$ to $k=K-1$, we have
\begin{align}
    &\frac{1}{K}\sum_{k=0}^{K-1}\left(V_{L,k}^{\pi^{\star}}(\rho)-V_{L,k}^{\widebar{\pi}^{k}}(\rho)\right)\notag\\
    &\leq \frac{N}{(1-\gamma)K\alpha}\mathbb{E}_{s\sim d_{\rho}^{\pi^{\star}}}\left[D_{\text{KL}}(\pi^{\star}(\cdot\mid s)||\pi_0(\cdot\mid s))\right]+\frac{1}{(1-\gamma)K}\sum_{k=0}^{K-1}\left(V_{L,k}^{\widebar{\pi}^{k+1}}(d_{\rho}^{\pi^{\star}})-V_{L,k}^{\widebar{\pi}^{k}}(d_{\rho}^{\pi^{\star}})\right)\notag\\
    &\hspace{20pt} +\frac{3(B_{\lambda}+1/N)}{(1-\gamma)^2 K} \sum_{k=0}^{K-1}\sum_{i=1}^{N}\|Q_i^{\pi_i^k}-Q_i^k\| +\frac{3\sqrt{|\Scal|| \Acal|}(B_{\lambda}+1/N)}{(1-\gamma)^4 K}\sum_{k=0}^{K-1}\sum_{i=1}^{N}\|\widebar{\pi}^k-\pi_i^k\|.
    \label{lem:bounded_Lagrangian:eq3}
\end{align}

The second term on the right hand side can be decomposed as follows
\begin{align}
    &\frac{1}{(1-\gamma)K}\sum_{k=0}^{K-1}\left(V_{L,k}^{\widebar{\pi}^{k+1}}(d_{\rho}^{\pi^{\star}})-V_{L,k}^{\widebar{\pi}^{k}}(d_{\rho}^{\pi^{\star}})\right)\notag\\
    &\leq \frac{1}{(1-\gamma)K}\sum_{k=0}^{K-1}\left(V_{0}^{\widebar{\pi}^{k+1}}(d_{\rho}^{\pi^{\star}})-V_{0}^{\widebar{\pi}^{k}}(d_{\rho}^{\pi^{\star}})\right)\notag\\
    &\hspace{20pt}+\frac{1}{(1-\gamma)K}\sum_{k=0}^{K-1}(\lambda^k-\nu^k)^{\top}\left(V_{g}^{\widebar{\pi}^{k+1}}(d_{\rho}^{\pi^{\star}})-V_{g}^{\widebar{\pi}^{k}}(d_{\rho}^{\pi^{\star}})\right)\notag\\
    &= \frac{V_{0}^{\widebar{\pi}^{K}}(d_{\rho}^{\pi^{\star}})}{(1-\gamma)K}+\frac{1}{(1-\gamma)K}\sum_{k=0}^{K-1}\left((\lambda^{k+1}-\nu^{k+1})^{\top}V_{g}^{\widebar{\pi}^{k+1}}(d_{\rho}^{\pi^{\star}})-(\lambda^k-\nu^k)^{\top}V_{g}^{\widebar{\pi}^{k}}(d_{\rho}^{\pi^{\star}})\right)\notag\\
    &\hspace{20pt}+\frac{1}{(1-\gamma)K}\sum_{k=0}^{K-1}\left(\lambda^{k}-\nu^k-\lambda^{k+1}+\nu^{k+1}\right)^{\top}V_{g}^{\widebar{\pi}^{k+1}}(d_{\rho}^{\pi^{\star}})\notag\\
    &= \frac{V_{0}^{\widebar{\pi}^{K}}(d_{\rho}^{\pi^{\star}})}{(1-\gamma) K}+\frac{1}{(1-\gamma) K}\sum_{i=1}^{N}(\lambda_i^K-\nu_i^K)V_i^{\widebar{\pi}^K}(d_{\rho}^{\pi^{\star}})\label{lem:bounded_Lagrangian:eq4}\\
    &\hspace{20pt}+\frac{1}{(1-\gamma)K}\sum_{k=0}^{K-1}\sum_{i=1}^{N}(\lambda_i^k-\nu_i^k-\lambda_i^{k+1}+\nu_i^{k+1})V_{i}^{\widebar{\pi}^{k+1}}(d_{\rho}^{\pi^{\star}}).\notag
\end{align}

We know that the value functions are bounded between $[0,\frac{1}{1-\gamma}]$. The projection in the update of the dual variable in Eq.~\eqref{prop:actor_conv:updates} guarantees $\lambda_i^k\in[0,B_{\lambda}]$. It is also straightforward to see that 
\begin{align*}
    |\lambda_{i,k}\hspace{-2pt}-\hspace{-2pt}\lambda_{i,k+1}|\leq \frac{\eta}{1-\gamma}+B\eta,\,|\nu_{i,k}-\nu_{i,k+1}|\leq \frac{\eta}{1-\gamma}+B\eta,\,\forall i=1,2,\cdot,N, \,\,k=0,1,\cdots,K-1.
\end{align*}

Using these bounds in Eq.~\eqref{lem:bounded_Lagrangian:eq4}, we get
\begin{align}
    \frac{1}{(1-\gamma)K}\sum_{k=0}^{K-1}\left(V_{L,k}^{\widebar{\pi}^{k+1}}(d_{\rho}^{\pi^{\star}})-V_{L,k}^{\widebar{\pi}^{k}}(d_{\rho}^{\pi^{\star}})\right) &\leq \frac{1}{(1-\gamma)^2 K}+\frac{N B_{\lambda}}{(1-\gamma)^2 K}\hspace{-2pt}+\hspace{-2pt}\frac{2N\eta}{(1-\gamma)^3 K}\hspace{-2pt}+\hspace{-2pt}\frac{2NB\eta}{(1-\gamma)^2 K}\notag\\
    &\leq\frac{2N B_{\lambda}}{(1-\gamma)^2 K}+\frac{4N\eta}{(1-\gamma)^3 K}. \label{lem:bounded_Lagrangian:eq5}
\end{align}

Finally, combining Eqs.~\eqref{lem:bounded_Lagrangian:eq3} and \eqref{lem:bounded_Lagrangian:eq5} yields
\begin{align*}
    &\frac{1}{K}\sum_{k=0}^{K-1}\left(V_{L,k}^{\pi^{\star}}(\rho)-V_{L,k}^{\widebar{\pi}^{k}}(\rho)\right)\\
    &\leq \frac{N}{(1\hspace{-2pt}-\hspace{-2pt}\gamma)K\alpha}\mathbb{E}_{s\sim d_{\rho}^{\pi^{\star}}}\left[D_{\text{KL}}(\pi^{\star}(\cdot\mid s)||\pi_0(\cdot\mid s))\right]\hspace{-2pt}+\hspace{-2pt}\frac{3(B_{\lambda}+1/N)}{(1-\gamma)^2 K} \sum_{k=0}^{K-1}\sum_{i=1}^{N}\|Q_i^{\pi_i^k}-Q_i^k\| \notag\\
    &\hspace{20pt}+\frac{3\sqrt{|\Scal|| \Acal|}(B_{\lambda}+1/N)}{(1-\gamma)^4 K}\sum_{k=0}^{K-1}\sum_{i=1}^{N}\|\widebar{\pi}^k-\pi_i^k\|+\frac{2N B_{\lambda}}{(1-\gamma)^2 K}+\frac{4N\eta}{(1-\gamma)^3 K},
\end{align*}
which leads to the claimed result by recognizing the fact that for $D_{\text{KL}}(p_1||p_2)\leq\log|\Acal|$ for $p_1,p_2\in\Delta_{\Acal}$ if $p_2$ is a uniform distribution.

\qed

\subsection{Proof of Lemma~\ref{lem:p_bound}}
From the definition of $p_i^{k,t}$ in Eq.~\eqref{eq:def_Gamma_p},
\begin{align*}
    \|p_i^{k,t}\| &=\|\omega_i^{k,t+1}-\widehat{\omega}_i^{k,t+1}\|\notag\\
    &\leq\|\omega_i^{k,t+1}-\omega_i^{k,t}\|+\|\omega_i^{k,t}-\widehat{\omega}_i^{k,t+1}\|\notag\\
    &\leq 2\|\omega_i^{k,t}-\widehat{\omega}_i^{k,t+1}\|\notag\\
    &\leq2\beta\|r_i(s_i^{k,t},a_i^{k,t})+(\gamma \phi(s_i^{k,t+1},a_i^{k,t+1})-\phi(s_i^{k,t},a_i^{k,t}))^{\top}\omega_i^{k,t}\|\notag\\
    &\leq 2(1+2B_{\omega})\beta,
\end{align*}
where the second inequality is due to the fact that $\Pi_{B_{\omega}}$ is the projection to a convex set and $\Pi_{B_{\omega}}\omega_i^{k,t}=\omega_i^{k,t}$.

Similarly,
\begin{align*}
    \|z_i^{k,t+1}-z_i^{k,t}\|=\|\omega_i^{k,t+1}-\omega_i^{k,t}\|\leq\|\omega_i^{k,t}-\widehat{\omega}_i^{k,t+1}\|\leq (1+2B_{\omega})\beta.
\end{align*}

\qed

\subsection{Proof of Lemma \ref{lem:A_negdef}}
Recall the definition of matrix $M^{\pi}$ in Eq.~\eqref{eq:def_M}. Given any $\pi$, we define the matrix $\widetilde{P}^{\pi}\in\mathbb{R}^{|\Scal||\Acal|\times|\Scal||\Acal|}$ such that $\widetilde{P}^{\pi}(s',a'\mid s,a)=P(s'\mid s,a)\pi(a'\mid s')$. Then, we have for any two vectors $\omega_1,\omega_2\in\mathbb{R}^d$
\begin{align*}
    \mathbb{E}_{\widehat{\pi}_i^k}[\omega_1^{\top}\phi(s,a)\phi(s',a')^{\top}\omega_2]&=\sum_{s,s',a,a'}\widetilde{\mu}_{\widehat{\pi}_i^k}(s,a)\widetilde{P}^{\widehat{\pi}_i^k}(s',a'\mid s,a) \omega_1^{\top}\phi(s,a)\phi(s',a')^{\top}\omega_2\notag\\
    &=\omega_1^{\top}\Phi^{\top}M^{\widehat{\pi}_i^k} \widetilde{P}^{\widehat{\pi}_i^k}\Phi\omega_2.
\end{align*}
Since this is true for any $\omega_1,\omega_2$, we know
\begin{align*}
    \mathbb{E}_{\widehat{\pi}_i^k}[\phi(s,a)\phi(s',a')^{\top}]=\Phi^{\top}M^{\widehat{\pi}_i^k} \widetilde{P}^{\widehat{\pi}_i^k}\Phi.
\end{align*}

Similarly, we can show 
\begin{align*}
    \mathbb{E}_{\widehat{\pi}_i^k}[\phi(s,a)\phi(s,a)^{\top}] &=\Phi^{\top}M^{\widehat{\pi}_i^k}\Phi,\\
    \mathbb{E}_{\widehat{\pi}_i^k}[r(s,a)\phi(s,a)] &=\Phi^{\top}M^{\widehat{\pi}_i^k}r_i,
\end{align*}
where $r_i=[r_i(s_0,a_0),\cdots,r_i(s_{|\Scal|},a_{|\Acal|})]^{\top}\in\mathbb{R}^{|\Scal||\Acal|}$.

We define $e_i^{\widehat{\pi}_i^k}(\omega)=\widebar{H}^{\widehat{\pi}_i^k}\omega+\widebar{b}_i^{\widehat{\pi}_i^k}$ and note that $e_i^{\widehat{\pi}_i^k}(\omega_i^{\star}(\widehat{\pi}_i^k))=0$. From the equations above and the definition of $\widebar{H}^{\pi}$ and $\widebar{b}_i^{\pi}$,
\begin{align*}
    e_i^{\widehat{\pi}_i^k}(\omega)&=\widebar{H}^{\widehat{\pi}_i^k}\omega+\widebar{b}_i^{\widehat{\pi}_i^k}\notag\\
    &=\mathbb{E}_{\widehat{\pi}_i^k}\big[\phi(s,a)\big((\gamma\phi(s',s')^{\top}-\phi(s,a)^{\top})\,\omega-r_i(s,a)\big)\big]\notag\\
    &=(\gamma\Phi^{\top}M^{\widehat{\pi}_i^k}\widetilde{P}^{\widehat{\pi}_i^k} \Phi-\Phi^{\top}M^{\widehat{\pi}_i^k}\Phi)\omega+\Phi^{\top}M^{\widehat{\pi}_i^k}r_i\notag\\
    &=\Phi^{\top}M^{\widehat{\pi}_i^k}(\gamma \widetilde{P}^{\widehat{\pi}_i^k}\Phi \omega-\Phi\omega+r_i)\notag\\
    &=\Phi^{\top}M^{\widehat{\pi}_i^k}(T_i^{\widehat{\pi}_i^k}(\Phi \omega)-\Phi\omega),
\end{align*}
where $T_i^{\pi}:\mathbb{R}^{|\Scal||\Acal|}\rightarrow\mathbb{R}^{|\Scal||\Acal|}$ for any policy $\pi$ is the Bellman operator defined such that
\begin{align}
    (T_i^{\pi}Q)(s,a) = \mathbb{E}_{\pi}[r_i(s,a)+\gamma Q(s',a')],\,\,\forall Q\in\mathbb{R}^{|\Scal||\Acal|}.\label{lem:A_negdef:eq1}
\end{align}

The projection of a vector to the span of $\Phi$ under the weighted $M^{\widehat{\pi}_i^k}$ norm is carried out through the projection matrix $\Pi_{\Phi}^{\widehat{\pi}_i^k}$
\begin{align*}
    \Pi_{\Phi}^{\widehat{\pi}_i^k}=\Phi(\Phi^{\top} M^{\widehat{\pi}_i^k} \Phi)^{-1} \Phi^{\top} M^{\widehat{\pi}_i^k}.
\end{align*}
It is obvious that $\Phi^{\top} M^{\widehat{\pi}_i^k}\Pi_{\Phi}^{\widehat{\pi}_i^k}=\Phi^{\top} M^{\widehat{\pi}_i^k}$. 

Eq.~\eqref{lem:A_negdef:eq1} implies that for any $Q_1, Q_2\in\mathbb{R}^{|\Scal||\Acal|}$
\begin{align}
    \|T_i^{\widehat{\pi}_i^k}Q_1-T_i^{\widehat{\pi}_i^k}Q_2\|_{M^{\widehat{\pi}_i^k}}^2
    &=\gamma(Q_1-Q_2)^{\top}(\widetilde{P}^{\widehat{\pi}_i^k})^{\top}M^{\widehat{\pi}_i^k}\widetilde{P}^{\widehat{\pi}_i^k}(Q_1-Q_2)\notag\\
    &=\sum_{s,a}\widetilde{\mu}_{\widehat{\pi}_i^k}(s,a)\Big(\sum_{s',a'}P(s'\mid s,a)\widehat{\pi}_i^k(a'\mid s')(Q_1-Q_2)(s',a')\Big)^2\notag\\
    &\leq \sum_{s,a}\widetilde{\mu}_{\widehat{\pi}_i^k}(s,a)\sum_{s',a'}P(s'\mid s,a)\widehat{\pi}_i^k(a'\mid s')(Q_1(s',a')-Q_2(s',a'))^2\notag\\
    &\leq \sum_{s',a'}\widetilde{\mu}_{\widehat{\pi}_i^k}(s',a')(Q_1(s',a')-Q_2(s',a'))^2\notag\\
    &=\|Q_1-Q_2\|_{M^{\widehat{\pi}_i^k}}^2,\label{lem:A_negdef:eq2}
\end{align}
where the inequality follows from Jensen's inequality. Eq.~\eqref{lem:A_negdef:eq2} implies another property of $\Pi_{\Phi}^{\widehat{\pi}_i^k}$, which is the contraction of $\Pi_{\Phi}^{\widehat{\pi}_i^k} T_i^{\widehat{\pi}_i^k}$ under the weighted $M^{\widehat{\pi}_i^k}$ norm. Specifically, we have for any $\omega\in\mathbb{R}^d$
\begin{align*}
    \|\Pi_{\Phi}^{\widehat{\pi}_i^k} T_i^{\widehat{\pi}_i^k}(\Phi\omega)-\Phi\omega_i^{\star}(\widehat{\pi}_i^k)\|_{M^{\widehat{\pi}_i^k}} &=\|\Pi_{\Phi}^{\widehat{\pi}_i^k} T_i^{\widehat{\pi}_i^k}(\Phi\omega)-\Pi_{\Phi}^{\widehat{\pi}_i^k} T_i^{\widehat{\pi}_i^k}(\Phi\omega_i^{\star})\|_{M^{\widehat{\pi}_i^k}}\notag\\
    &\leq\|T_i^{\widehat{\pi}_i^k}(\Phi\omega) - T_i^{\widehat{\pi}_i^k}(\Phi\omega_i^{\star}(\widehat{\pi}_i^k))\|_{M^{\widehat{\pi}_i^k}}\notag\\
    &\leq\gamma\|\Phi\omega-\Phi\omega_i^{\star}(\widehat{\pi}_i^k)\|_{M^{\widehat{\pi}_i^k}}.
\end{align*}

Then, we have for any $\omega\in\mathbb{R}^d$,
\begin{align*}
    &(\omega-\omega_i^{\star}(\widehat{\pi}_i^k))^{\top}\widebar{H}^{\widehat{\pi}_i^k} (\omega-\omega_i^{\star}(\widehat{\pi}_i^k))\\
    &= (\omega-\omega_i^{\star}(\widehat{\pi}_i^k))^{\top}(e_i^{\widehat{\pi}_i^k}(\omega)-e_i^{\widehat{\pi}_i^k}(\omega_i^{\star}(\widehat{\pi}_i^k)))\notag\\
    &=(\omega-\omega_i^{\star}(\widehat{\pi}_i^k))^{\top}e_i^{\widehat{\pi}_i^k}(\omega)\notag\\
    &=(\omega-\omega_i^{\star}(\widehat{\pi}_i^k))^{\top}\Phi^{\top}M^{\widehat{\pi}_i^k}(T_i^{\widehat{\pi}_i^k}(\Phi \omega)-\Phi\omega)\notag\\
    &=(\omega-\omega_i^{\star}(\widehat{\pi}_i^k))^{\top}\Phi^{\top}M^{\widehat{\pi}_i^k}((I-\Pi_{\Phi}^{\widehat{\pi}_i^k})T_i^{\widehat{\pi}_i^k}(\Phi \omega)+\Pi_{\Phi}^{\widehat{\pi}_i^k} T_i^{\widehat{\pi}_i^k}(\Phi \omega)-\Phi\omega)\notag\\
    &=(\omega-\omega_i^{\star}(\widehat{\pi}_i^k))^{\top}\Phi^{\top}M^{\widehat{\pi}_i^k}(\Pi_{\Phi}^{\widehat{\pi}_i^k} T_i^{\widehat{\pi}_i^k}(\Phi \omega)-\Phi\omega),
\end{align*}
where the last equality follows from $\Phi^{\top} M^{\widehat{\pi}_i^k}\Pi_{\Phi}^{\widehat{\pi}_i^k}=\Phi^{\top} M^{\widehat{\pi}_i^k}$.

Using the contraction of $\Pi_{\Phi}^{\widehat{\pi}_i^k} T_i^{\widehat{\pi}_i^k}$ and the Cauchy-Schwarz inequality, we have for any $\omega\in\mathbb{R}^d$
\begin{align*}
    &(\omega-\omega_i^{\star}(\widehat{\pi}_i^k))^{\top}\widebar{H}^{\widehat{\pi}_i^k} (\omega-\omega_i^{\star}(\widehat{\pi}_i^k))\notag\\
    &=(\omega-\omega_i^{\star}(\widehat{\pi}_i^k))^{\top}\Phi^{\top}M^{\widehat{\pi}_i^k}(\Pi_{\Phi}^{\widehat{\pi}_i^k} T_i^{\widehat{\pi}_i^k}(\Phi \omega)-\Phi\omega)\notag\\
    &=(\Phi\omega-\Phi\omega_i^{\star}(\widehat{\pi}_i^k))^{\top}M^{\widehat{\pi}_i^k}(\Pi_{\Phi}^{\widehat{\pi}_i^k} T_i^{\widehat{\pi}_i^k}(\Phi \omega)-\Phi\omega_i^{\star}(\widehat{\pi}_i^k))+(\Phi\omega-\Phi\omega_i^{\star}(\widehat{\pi}_i^k))^{\top}M^{\widehat{\pi}_i^k}(\Phi\omega_i^{\star}(\widehat{\pi}_i^k)-\Phi\omega)\notag\\
    &\leq \|\Phi\omega-\Phi\omega_i^{\star}(\widehat{\pi}_i^k)\|_{M^{\widehat{\pi}_i^k}} \|\Pi_{\Phi}^{\widehat{\pi}_i^k} T_i^{\widehat{\pi}_i^k}(\Phi \omega)-\Phi\omega_i^{\star}(\widehat{\pi}_i^k)\|_{M^{\widehat{\pi}_t}}-\|\Phi\omega-\Phi\omega_i^{\star}(\widehat{\pi}_i^k)\|_{M^{\widehat{\pi}_i^k}}^2\notag\\
    &\leq (\gamma-1)\|\Phi\omega-\Phi\omega_i^{\star}(\widehat{\pi}_i^k)\|_{M^{\widehat{\pi}_i^k}}^2\notag\\
    &\leq (\gamma-1)\sigma_{\min}(\Phi)\sigma_{\min}(M^{\widehat{\pi}_i^k})\|\omega-\omega_i^{\star}(\widehat{\pi}_i^k)\|^2\notag\\
    &\leq (\gamma-1)\sigma_{\min}(\Phi)\min_{s,a}\mu_{\widehat{\pi}_i^k}(s)\widehat{\pi}_t(a\mid s)\|\omega-\omega_i^{\star}(\widehat{\pi}_i^k)\|^2,
\end{align*}
where $\sigma_{\min}$ denotes the singular value with the smallest magnitude, and the last inequality follows from the fact that the singular values of a diagonal matrix are the diagonal entries.

\qed

\subsection{Proof of Lemma \ref{lem:bound_Gamma}}
Recall that the Markovian samples generated by our algorithm are
\begin{align*}
{s}_i^{k,t-\tau} \stackrel{\widehat{\pi}_i^k}{\longrightarrow}  a_i^{k,t-\tau}\stackrel{P}{\longrightarrow}  s_i^{k,t-\tau+1} \stackrel{\widehat{\pi}_i^k}{\longrightarrow} a_i^{k,t-\tau+1}\stackrel{P}{\longrightarrow} \cdots\stackrel{P}{\longrightarrow} s_i^{k,t} \stackrel{\widehat{\pi}_i^k}{\longrightarrow}  a_i^{k,t} \stackrel{P}{\longrightarrow}  s_i^{k,t+1} \stackrel{\widehat{\pi}_i^k}{\longrightarrow}  a_i^{k,t+1}.
\end{align*}

Let $\widebar{O}_i^k=(\widebar{s}_i^k,\widebar{a}_i^k,\widebar{s}_i^{k\prime},\widebar{a}_i^{k\prime})$ where $\widebar{s}_i^k\sim\mu_{\widehat{\pi}_i^k}$, $\widebar{a}_i^k\sim\widehat{\pi}_i^k(\cdot\mid\widebar{s}_i^k)$, $\widebar{s}_i^{k\prime}\sim P(\cdot\mid\widebar{s}_i^k,\widebar{a}_i^k)$, and $\widebar{a}_i^{k\prime}\sim\widehat{\pi}_i^k(\cdot\mid\widebar{s}_i^{k\prime})$. 
We can decompose the term of interest as
\begin{align}
    \mathbb{E}[\Gamma_i(\widehat{\pi}_i^k,z_i^{k,t},O_i^{k,t})]&=\mathbb{E}[\Gamma_i(\widehat{\pi}_i^k,z_i^{k,t},O_i^{k,t})-\Gamma_i(\widehat{\pi}_i^k,z_i^{k,t-\tau},O_i^{k,t})]\notag\\
    &\hspace{20pt}+\mathbb{E}[\Gamma_i(\widehat{\pi}_i^k,z_i^{k,t-\tau},O_i^{k,t})-\Gamma_i(\widehat{\pi}_i^k,z_i^{k,t-\tau},\widebar{O}_i^{k})]+\mathbb{E}[\Gamma_i(\widehat{\pi}_i^k,z_i^{k,t-\tau},\widebar{O}_i^{k})].\label{lem:bound_Gamma_decomposition}
\end{align}

To bound the first term of Eq.~\eqref{lem:bound_Gamma_decomposition}, 
\begin{align}
    &\Gamma_i(\widehat{\pi}_i^k,z_i^{k,t},O_i^{k,t})-\Gamma_i(\widehat{\pi}_i^k,z_i^{k,t-\tau},O_i^{k,t})\notag\\
    &\leq (z_i^{k,t}-z_i^{k,t-\tau})^{\top}\left(R_i(O_i^{k,t})+H(O_i^{k,t})\omega_i^{\star}(\widehat{\pi}_i^k)\right)\notag\\
    &\hspace{20pt}+(z_i^{k,t})^{\top}(H(O_i^{k,t})-\widebar{H}^{\widehat{\pi}_i^k})z_i^{k,t}-(z_i^{k,t-\tau})^{\top}(H(O_i^{k,t})-\widebar{H}^{\widehat{\pi}_i^k})z_i^{k,t-\tau}\notag\\
    &\leq \|z_i^{k,t}-z_i^{k,t-\tau}\|\left(1+2B_{\omega}\right)+(z_i^{k,t}-z_i^{k,t-\tau})^{\top}(H(O_i^{k,t})-\widebar{H}^{\widehat{\pi}_i^k})z_i^{k,t-\tau}\notag\\
    &\hspace{20pt}+(z_i^{k,t})^{\top}(H(O_i^{k,t})-\widebar{H}^{\widehat{\pi}_i^k})(z_i^{k,t}-z_i^{k,t-\tau})\notag\\
    &\leq (1+2B_{\omega})\|z_i^{k,t}-z_i^{k,t-\tau}\| +16B_{\omega}\sum_{t'=t-\tau}^{t-1}\|z_i^{k,t'+1}-z_i^{k,t'}\|\notag\\
    &\leq (1+18B_{\omega})\cdot\tau(1+2B_{\omega})\beta\notag\\
    &\leq(1+18B_{\omega})^2\beta\tau.
    \label{lem:bound_Gamma:eq1}
\end{align}

Let $\Fcal_i^{k,t}$ denote the past randomness in outer loop iteration $k$ up to inner loop iteration $t$ at agent $i$, i.e. $\Fcal_i^{k,t}=\{O_i^{k,0},O_i^{k,1},\dots,O_i^{k,t}\}$. To treat the second term of Eq.~\eqref{lem:bound_Gamma_decomposition}, we have for all $t\geq\tau$
\begin{align}
    &\mathbb{E}[\Gamma_i(\widehat{\pi}_i^k,z_i^{k,t-\tau},O_i^{k,t})\mid\Fcal_i^{k,t-\tau}-\Gamma_i(\widehat{\pi}_i^k,z_i^{k,t-\tau},\widebar{O}_i^k)\mid\Fcal_i^{k,t-\tau}]\notag\\
    &=(z_i^{k,t-\tau})^{\top}\mathbb{E}[R_i(O_i^{k,t})-R_i(\widebar{O}_i^k)\mid\Fcal_i^{k,t-\tau}]+(z_i^{k,t-\tau_k})^{\top}\mathbb{E}[H(O_i^{k,t})-H(\widebar{O}_i^k)\mid\Fcal_i^{k,t-\tau}]\omega_i^{\star}(\widehat{\pi}_i^k)\notag\\
    &\hspace{20pt}+(z_i^{k,t-\tau})^{\top}\mathbb{E}[H(O_i^{k,t})-H(\widebar{O}_i^k)\mid\Fcal_i^{k,t-\tau}]z_i^{k,t-\tau}\notag\\
    &\leq \|z_i^{k,t-\tau}\|\|\mathbb{E}[R_i(O_i^{k,t})-R_i(\widebar{O}_i^k)\mid\Fcal_i^{k,t-\tau}]\|+\|z_i^{k,t-\tau}\|\mathbb{E}[H(O_i^{k,t})-H(\widebar{O}_i^k)\mid\Fcal_i^{k,t-\tau}]\|\|\omega_i^{\star}(\widehat{\pi}_i^k)\|\notag\\
    &\hspace{20pt}+\|z_i^{k,t-\tau}\|^2\|\mathbb{E}[H(O_i^{k,t})-H(\widebar{O}_i^k)\mid\Fcal_i^{k,t-\tau}]\|\notag\\
    &\leq B_{\omega}\cdot d_{TV}(O_i^{k,t},\widebar{O}_i^k\mid\Fcal_i^{k,t-\tau})+2B_{\omega}\cdot B_{\omega}\cdot 2d_{TV}(O_i^{k,t},\widebar{O}_i^k\mid\Fcal_i^{k,t-\tau})\notag\\
    &\hspace{20pt}+4B_{\omega}^2\cdot 2d_{TV}(O_i^{k,t},\widebar{O}_i^k\mid\Fcal_i^{k,t-\tau})\notag\\
    &\leq (B_{\omega}+12B_{\omega}^2)d_{TV}(O_i^{k,t},\widebar{O}_i^k\mid\Fcal_i^{k,t-\tau})\notag\\
    &\leq (B_{\omega}+12B_{\omega}^2)C_0\ell^k\notag\\
    &\leq (B_{\omega}+12B_{\omega}^2)\beta,
    \label{lem:bound_Gamma:eq3}
\end{align}
where the fourth inequality follows from an argument similar to the one in \cite{wu2020finite}[Lemma D.11], the last inequality uses Assumption~\ref{assump:markov-chain}, and the last inequality is a result of Eq.~\eqref{eq:tau_def}.

By the definition of $H(\widebar{O}_i^k),\widebar{H}^{\widehat{\pi}_i^k}$, and $\omega_i^{\star}(\widehat{\pi}_i^k)$, we have for the last term of Eq.~\eqref{lem:bound_Gamma_decomposition}
\begin{align}
    \mathbb{E}[\Gamma_i(\widehat{\pi}_i^k,z_i^{k,t-\tau},\widebar{O}_i^k)\mid\Fcal_i^{k,t-\tau}]&=\mathbb{E}[(z_i^{k,t-\tau})^{\top}(R_i(\widebar{O}_i^k)+H(\widebar{O}_i^k)\omega_i^{\star}(\widehat{\pi}_i^k))\mid\Fcal_i^{k,t-\tau}]\notag\\
    &\hspace{20pt}+\mathbb{E}[(z_i^{k,t-\tau})^{\top}(H(\widebar{O}_i^k)-\widebar{H}^{\widehat{\pi}_i^k})z_i^{k,t-\tau}\mid\Fcal_i^{k,t-\tau}]\notag\\
    &=0+0=0.\label{lem:bound_Gamma:eq2}
\end{align}

Plugging Eqs.~\eqref{lem:bound_Gamma:eq1} and \eqref{lem:bound_Gamma:eq3} into Eq.~\eqref{lem:bound_Gamma_decomposition}, we have
\begin{align*}
    \mathbb{E}[\Gamma_i(\widehat{\pi}_i^k,z_i^{k,t},O_i^{k,t})]&\leq (1+18B_{\omega})^2\beta\tau+(B_{\omega}+12B_{\omega}^2)\beta\notag\\
    &\leq2(1+18B_{\omega})^2\beta\tau.
\end{align*}

\qed

\end{document}